\documentclass[11pt,a4paper]{article}
\usepackage{ifthen,latexsym,amssymb,amsmath,bbm,fixmath}
\usepackage{hyperref}
\usepackage[nobysame,initials]{amsrefs}
\usepackage{bm,booktabs,float}
\usepackage{graphicx,tikz}
\usepackage{authblk}
\usepackage{enumitem}
\renewcommand{\ldots}{\hspace{0.9pt}.\hspace{0.3pt}.\hspace{0.3pt}.\hspace{1.5pt}}
\renewcommand{\ge}{\geqslant}

\newif\ifnotesw\noteswtrue

\newcommand{\hide}[1]{}

\newcommand{\beq}[1]{\begin{equation}\label{#1}}
\newcommand{\eeq}{\end{equation}}

\newtheorem{theorem}{Theorem}[section]
\newtheorem{lemma}[theorem]{Lemma}
\newtheorem{problem}[theorem]{Problem}
\newtheorem{conjecture}[theorem]{Conjecture}
\newtheorem{corollary}[theorem]{Corollary}
\newtheorem{definition}[theorem]{Definition}
\newtheorem{proposition}[theorem]{Proposition}

\newtheorem{remark}[theorem]{Remark}

\newcommand{\bpf}[1][Proof.]{\smallskip\noindent{\it #1} }
\newcommand{\qed}{\nolinebreak\mbox{\hspace{5 true pt}%
  \rule[-0.85 true pt]{3.9 true pt}{8.1 true pt}}}
\newcommand{\epf}{\qed \medskip}

\newtheorem{claim}{Claim}[theorem]
\newcommand{\cqed}{\nolinebreak\mbox{\hspace{5 true pt}%
  \rule[-0.85 true pt]{2.0 true pt}{8.1 true pt}}}
\newcommand{\bcpf}{\bpf[Proof of Claim.]}
\newcommand{\ecpf}{\cqed \medskip}

\begin{document}

%text specific macros 

\def\eps{{\varepsilon}}
\newcommand{\cP}{\mathcal{P}}
\newcommand{\cT}{\mathcal{T}}
\newcommand{\cL}{\mathcal{L}}
\newcommand{\ex}{\mathbb{E}}
\newcommand{\eul}{e}
\newcommand{\pr}{\mathbb{P}}
\newcommand{\supp}{\mathop{\mathrm{supp}}}
\newcommand{\Capa}{\mathop{\mathrm{Cap}}}
\newcommand{\opt}{\mathop{\textsc{opt}}}
\newcommand{\feas}{\mathop{\textsc{feas}}}
\newcommand{\pat}{\mathop{\textsc{pat}}}
\newcommand{\wt}{\mathop{\textsc{wt}}}
\newcommand{\coars}{\mathop{\mathrm{Coars}}}
\newcommand{\RL}{\mathop{\mathrm{RL}}}
\newcommand{\ext}{\mathop{\mathrm{ext}}}
\newcommand{\ba}{\bm{\alpha}}
\newcommand{\bb}{\bm{\beta}}
\newcommand{\perf}{\mathrm{perf}}
\newcommand{\ed}{d}

\title{Exact solutions to the Erd\H{o}s-Rothschild problem}
\author[1]{Oleg Pikhurko\footnote{Supported by ERC Advanced Grant 101020255 and Leverhulme Research Project Grant RPG-2018-424}}
\affil[1]{Mathematics Institute and DIMAP, University of Warwick, Coventry CV4 7AL, UK}
\author[2]{Katherine Staden\footnote{Supported by EPSRC Fellowship EP/V025953/1}}
\affil[2]{School of Mathematics and Statistics, The Open University, Milton Keynes, MK7 6AA, UK}
\maketitle

\begin{abstract}
Let $\bm{k} := (k_1,\ldots,k_s)$ be a sequence of natural numbers.
For a graph $G$, let $F(G;\bm{k})$ denote the number of colourings of the edges of $G$ with colours $1,\dots,s$ such that, for every  $c \in \{1,\dots,s\}$,
the edges of colour $c$ contain no clique of order $k_c$. 
Write $F(n;\bm{k})$ to denote the maximum of $F(G;\bm{k})$ over all graphs $G$ on $n$ vertices.
There are currently very few known exact (or asymptotic) results for this problem, posed by Erd\H{o}s and Rothschild in 1974.
We prove some new exact results for $n \to \infty$:
\begin{itemize}
\item[(i)] A sufficient condition on $\bm{k}$ which guarantees that every extremal graph is a complete multipartite graph, which systematically recovers all existing exact results.
\item[(ii)] Addressing the original question of Erd\H{o}s and Rothschild, in the case $\bm{k}=(3,\ldots,3)$ of length $7$, the unique extremal graph is the complete balanced $8$-partite graph, with colourings coming from Hadamard matrices of order $8$.
\item[(iii)] In the case $\bm{k}=(k+1,k)$, for which the sufficient condition in~(i) does not hold, for $3 \leq k \leq 10$, the unique extremal graph is complete $k$-partite with one part of size less than $k$ and the other parts as equal in size as possible.
\end{itemize}

\end{abstract}

\section{Introduction}

Let a non-increasing sequence $\bm{k} = (k_1,\ldots,k_s) \in \mathbb{N}^s$ of natural numbers be given. By an \emph{$s$-edge-colouring} (or \emph{colouring} for brevity)
of a graph $G=(V,E)$ we mean a function $\chi:E\to [s]$, where we denote $[s]:=\{1,\dots,s\}$. Note that colourings do not have to be proper, that is, adjacent edges can have the same colour. A colouring $\chi$ of $G$ is called $\bm{k}$-\emph{valid} if, for every $c\in [s]$, the colour-$c$ subgraph $\chi^{-1}(c)$ contains no copy of $K_{k_c}$, the complete graph of order $k_c$. 
Write $F(G; \bm{k})$ for the number of $\bm{k}$-valid colourings of $G$.
This paper concerns the parameter
$$
F(n;\bm{k}) := \max_{G: v(G)=n}F(G;\bm{k}),
$$
the maximum value of $F(G;\bm{k})$ for an $n$-vertex graph $G$,
and the \emph{$\bm{k}$-extremal} graphs (i.e.~those graphs which attain this maximum).
We always assume $s \geq 2$ and $k_s \geq 3$ for otherwise the problem is trivial
or reduces to one with shorter $\bm{k}$.
Determining $F(n;\bm{k})$ in the case $\bm{k}=(3,3)$ was originally studied by Erd\H{o}s and Rothschild~\cite{ER,ER2}, and hence it is called the \emph{Erd\H{o}s-Rothschild problem}.
It is in general wide open. 
The following summarises all of the cases where $F(n;\bm{k})$ is known (in one case only asymptotically).
Write $(k;s)$ for the tuple $(k,\ldots,k)$ of length $s$.

\begin{theorem}\label{knownthm}
There exists $n_0(k)>0$ such that the following hold for all integers $n \geq n_0$.

\medskip
\rm
\begin{tabular}{ll|l|l|l}
$k$ & $s$         &                     $F(n;(k;s))$        & extremal graph(s)                           & citation \\ 
\hline
any $k$ & $s=2$   & $2^{t_{k-1}(n)}$                     & $T_{k-1}(n)$                 &         \cite{abks}                 \\
 & $s=3$  & $3^{t_{k-1}(n)}$                            & $T_{k-1}(n)$              &     \cite{abks}      \\
$k=3$ & $s=4$           &  $(C_{3,4}+o(1))\cdot 18^{t_4(n)/3}$                        & $T_{4}(n)$      &         \cite{abks,PY}  \\
 & $s=5$           &      $6^{t_2(n)+o(n^2)}$                 & $T_{\alpha,\beta}(n)$, $0 \leq \alpha,\beta \leq \frac{1}{4}$  \hfill  (*)           &     \cite{skokan} \\      
  &          &                     & $T'_{\alpha,\beta}(n)$, $0 \leq \alpha,\beta, \alpha+\beta \leq \frac{1}{4}$  \hfill  (*)           &     \cite{skokan} \\  
 & $s=6$           &          $(C_{3,6}+o(1))4^{t_4(\lfloor n/2\rfloor)}4^{t_4(\lceil n/2\rceil)} 3^{t_2(n)}$          & $T_{8}(n)$      &    \cite{skokan}           \\      
$k=4$ & $s=4$           &      $(C_{4,4}+o(1))\cdot 3^{t_9(n)}$                       & $T_{9}(n)$                        &    \cite{abks,PY}                  
\end{tabular}

\medskip
(*) These graphs are known to be asymptotically extremal only: they achieve the right exponent in $F(n;\bm{k})$.
Here, $T_{\alpha,\beta}(n)$ denotes the complete partite graph with parts of size
$\alpha n, \alpha n, (\frac{1}{4}-\alpha)n, (\frac{1}{4}-\alpha)n, \beta n, \beta n, (\frac{1}{4}-\beta)n, (\frac{1}{4}-\beta)n$, and 
$T'_{\alpha,\beta}(n)$ denotes the complete partite graph with parts of size
$\frac{n}{4},\frac{n}{4},\alpha n, \alpha n, \beta n, \beta n, (\frac{1}{4}-\alpha-\beta)n, (\frac{1}{4}-\alpha-\beta)n$.

\end{theorem}
The constants $C_{3,4},C_{3,6},C_{4,4}$ can be determined, and generally depend on the remainder when $n$ is divided by some small integer; for example, $C_{3,4}$ equals $(2^{14}\cdot 3)^{1/3}$ if $n \equiv 2~(\!\!\!\!\mod 4)$ and equals $36$ otherwise.
Note that for $\bm{k} = (k;s)$, the trivial lower bound for $F(n;\bm{k})$ is $s^{t_{k-1}(n)}$, given by taking every $s$-edge-colouring of the largest $K_k$-free graph on $n$ vertices, namely the Tur\'an graph $T_{k-1}(n)$, the complete partite graph with $k-1$ parts of size as equal as possible (with $t_{k-1}(n) := e(T_{k-1}(n))$).
This trivial lower bound is in fact sharp for $s=2,3$, but $F(n;\bm{k})$ is exponentially larger for $s \geq 4$, as was shown in~\cite{abks}. As is evident from the table, these cases have been much harder to resolve and there are only four pairs $(k;s)$ where the solution is known.
We refer the reader to~\cite{stability2} for a more detailed history of the problem and its variants.
This paper is the third in a series (comprising also \cite{psy} with Yilma and~\cite{stability2}) concerning the relationship between the Erd\H{o}s-Rothschild problem and a finite combinatorial optimisation problem, which we now state.

\medskip
\noindent\textbf{Problem~$Q^*$}:
 \it
 Given a sequence $\bm{k} := (k_1,\ldots,k_s) \in \mathbb{N}^s$ of natural numbers, determine
\begin{equation}\label{Qdef}
Q(\bm{k}) := \max_{(r,\phi,\ba) \in \feas^*(\bm{k})} q(\phi,\ba),
\end{equation}
 the maximum value of
 \begin{equation}\label{q}
 q(\phi,\ba) := 2\sum_{1 \leq i < j \leq r} \alpha_i \alpha_j \log_2 |\phi(ij)|
 \end{equation}
 over the set $\feas^*(\bm{k})$ of \emph{feasible solutions}, that is, triples
$(r,\phi,\ba)$ such that $r \in \mathbb{N}$ and
 %satisfying all of the following constraints:
 \begin{itemize}
  \item $\phi \in \Phi_2(r;\bm{k})$, where $\Phi_2(r;\bm{k})$ is the set of all
  functions $\phi : \binom{[r]}{2} \rightarrow 2^{[s]}$
  such that
  \[
  \phi^{-1}(c) := \left\{ ij \in \binom{[r]}{2} : c \in \phi(ij) \right\}
  \]
is $K_{k_c}$-free for every colour $c \in [s]$ and $|\phi(ij)|\ge 2$
for all $ij\in \binom{[r]}{2}$;
  \item $\ba = (\alpha_1,\ldots,\alpha_r) \in \Delta^r$, where $\Delta^r$ is the set of all $\ba \in \mathbb{R}^r$ with $\alpha_i > 0$ for all $i \in [r]$, and $\alpha_1 + \ldots + \alpha_r = 1$.
  %We call $\ba$ a \emph{vertex-weighting}.
  \end{itemize}
\rm

We may assume that $r < R(\bm{k})$, where $R(\bm{k})$ is the \emph{Ramsey number} of $\bm k$ (i.e.\ the minimum $R$ such that $K_R$ admits no $\bm k$-valid
  $s$-edge-colouring).
Note that the maximum in~\eqref{Qdef} is attained, since $q(\phi,\cdot)$ is continuous for each of the finitely many pairs $(r,\phi)$, and $\feas^*(\bm{k})$ with the weaker restrictions $\alpha_i\geq 0$ for every $i$ and $\ba$ is a (non-empty) compact space. 
(If the maximum is obtained at $(r,\phi,\ba)$ with some $\alpha_i=0$ then we can simply remove $i$ to obtain a solution with smaller $r$.)
We call $\phi \in \Phi_2(r;\bm{k})$ a \emph{colour pattern} and $\ba \in \Delta^r$ a \emph{vertex weighting}.
A triple $(r,\phi,\ba)$ is called \emph{optimal} if it attains the maximum, that is, 
$(r,\phi,\ba)\in \feas^*(\bm{k})$ and $q(r,\phi,\ba)=Q(\bm{k})$.
Let $\opt^*(\bm{k})$ be the set of optimal triples $(r,\phi,\ba)$.
Also let
$$
\wt(\bm{k}) := \{(r,\ba): \exists (r,\phi,\ba)\in\textstyle{\opt^*}(\bm{k})\}
$$
and
$$
\pat(r,\ba) := \{\phi \in \Phi_2(r; \bm{k}) : (r,\phi,\ba)\in \textstyle{\opt^*}(\bm{k})\}.
$$
%Given vectors $\bm{a}=(a_1,\ldots,a_s)$ and $\bm{b}=(b_1,\ldots,b_t)$, write $\bm{a}\leq\bm{b}$ if $a_i \leq b_i$ for all $i \in \max\{s,t\}$ where $a_i:=0$ for all $i>s$ and $b_i:=0$ for all $i>t$.
In~\cite{stability2} and~\cite{psy} we considered Problems $Q_0,Q_1,Q_2$ which are as Problem $Q^*$ except we relax $\alpha_i>0$ to $\alpha_i \geq 0$ and Problem $Q_t$ considers feasible $\phi \in \Phi_t(r;\bm{k})$ where
$$
\Phi_t(r;\bm{k}) := \left\{ \phi: \binom{[r]}{2} \to 2^{[s]} : \phi^{-1}(c) \text{ is }K_{k_c}\text{-free } \forall c \in [s], |\phi(ij)| \geq t \ \forall ij \in \binom{[r]}{2}\right\}
$$
(and $q(\phi,\ba)$ only sums over pairs $ij$ with $\phi(ij) \neq \emptyset$). 
It is not hard to show that the optimal value of each problem is the same. 
Clearly Problem $Q^*$ has the smallest feasible set.
However, we will sometimes consider also for $t=0,1$ the set $\opt_t(\bm{k})$ of triples $(r,\phi,\ba)$ with $\ba \in \Delta^r$ and $\phi \in \Phi_t(r;\bm{k})$ which attain the maximum value.
We write $\| \bm{a}-\bm{b}\|_1 := \sum_i|a_i-b_i|$ for the \emph{$\ell^1$-distance} between finite tuples $\bm{a}$ and $\bm{b}$ of real numbers (where in the sum we add trailing $0$'s to make $\bm{a},\bm{b}$ equal length).
In this paper we always take $\log$ to the base $2$; from now on we omit any subscript. 
We define $\mathbb{N} := \{1,2,3,\ldots\}$.

The goal of this series of works has been to verify the following meta-conjecture:
\begin{itemize}
\item[] \emph{To solve the Erd\H{o}s-Rothschild problem, it suffices to solve Problem $Q^*$.}
\end{itemize}

The main result of the present paper is that to determine $F(n;\bm{k})$ \emph{exactly}, it suffices to solve the optimisation problem, when $\bm{k}$ satisfies a certain condition.
This improves the main result of our previous paper~\cite{stability2} which had `approximately' in place of `exactly' -- we use the results of~\cite{stability2} as a crucial tool in the present paper.
In~\cite{psy} we proved with Yilma that for every $n \in \mathbb{N}$, at least one of the $\bm{k}$-extremal graphs of order $n$ is complete multipartite, and that
\begin{equation}\label{FQeq}
F(n;\bm{k}) = 2^{Q(\bm{k})\binom{n}{2}+o(n^2)}, \quad\text{ so }F(\bm{k})=Q(\bm{k}),\text{ where }\quad
F(\bm{k}) := \lim_{n\to\infty}\frac{\log_2 F(n;\bm{k})}{\binom{n}{2}}.
\end{equation}

In~\cite{stability2}, we proved a stability theorem for $\bm{k}$ satisfying a certain condition which we call the \emph{extension property}. This is a property of optimal solutions, that says an infinitesimal part added onto such a solution in an optimal way must look like a clone of an existing part.

\begin{definition}[Clones and extension property]\label{extprop}
\rm
Let $s \in \mathbb{N}$ and $\bm{k} \in \mathbb{N}^s$. 
Given $r \in \mathbb{N}$ and $\phi \in \Phi_0(r;\bm{k})$, say that $i \in [r]$ is
\begin{itemize}
\item a \emph{clone of $j \in [r]\setminus \{ i \}$} \emph{(under $\phi$)} if $\phi(ik) = \phi(jk)$ for all $k \in [r] \setminus \{ i,j\}$ and $|\phi(ij)| \leq 1$;
\item a \emph{strong clone of $j$} if additionally $\phi(ij) = \emptyset$.
\end{itemize}
We say that $\bm{k}$ has
\begin{itemize}
\item the \emph{extension property} if, for every $(r^*,\phi^*,\ba^*) \in \opt^*(\bm{k})$ and $\phi \in \Phi_0(r^*+1;\bm{k})$ such that $\phi|_{\binom{[r^*]}{2}} = \phi^*$ and $\ext(\phi,\ba^*)=Q(\bm{k})$, where
$$
\ext(\phi,\ba^*) := \sum_{\substack{i \in [r^*]\\ \phi(\{i,r^*+1\}) \neq \emptyset}}\alpha^*_i\log|\phi(\{ i,r^*+1 \})|,
$$
%which is the ``normalised contribution'' of the zero-weighted vertex $r+1$ to $q(\phi,(\alpha_1,\ldots,\alpha_r,0))$.
there exists $j \in [r^*]$ such that $r^*+1$ is a clone of $j$ under $\phi$;
\item the \emph{strong extension property} if in fact $r^*+1$ is a strong clone of $j$.
\end{itemize}
\end{definition}

The strong extension property holds in all but one of the cases where the problem has been solved (this was proved in~\cite{stability2} apart from in the case $(3;6)$ which we include in the present paper).

\begin{lemma}\label{recover}
Every $\bm{k}$ in Theorem~\ref{knownthm} apart from $(3;5)$ has the strong extension property. %Thus Theorem~\ref{stabilitysimp} applies to every such $\bm{k}$.
\end{lemma}

The stability theorem says that whenever $\bm{k}$ has the extension property, every almost extremal graph; that is, $G$ on $n$ vertices such that $F(G;\bm{k})=F(n;\bm{k})\cdot 2^{o(n^2)}$, looks like the blow-up of an optimal solution to Problem $Q^*$.
(For a definition of a $(\delta,d)$-regular pair see Section~\ref{sectools}.)

\begin{theorem}[Stability {\cite[Theorem~1.4]{stability2}}]\label{stabilitysimp}
Let $s \in \mathbb{N}$ and suppose that $\bm{k} \in \mathbb{N}^s$ with $k_1\geq\ldots\geq k_s$ has the extension property.
Then for all $\delta > 0$ there exist $n_0 \in \mathbb{N}$ and $\eps > 0$ such that the following holds.
Let $G$ be a graph on $n \geq n_0$ vertices such that
$$
\frac{\log F(G;\textbf{k})}{\binom{n}{2}} \geq Q(\textbf{k}) - \eps.
$$
Then, for at least $(1-2^{-\eps n^2})\cdot F(G;\bm{k})$ colourings $\chi : E(G) \rightarrow [s]$ which are $\bm{k}$-valid, there are $(r^*,\phi^*,\ba^*) \in \opt^*(\bm{k})$ and a partition $Y_1 \cup \ldots \cup Y_{r^*} = V(G)$ such that the following hold.
\begin{enumerate}[label=(\roman*),ref=(\roman*)]
\item\label{stabilityi} For all $i \in [r^*]$, we have that $|\,|Y_i| - \alpha^*_in\,| < 1$.
\item\label{stabilityii} for all $c \in \phi^*(ij)$ and $ij \in \binom{[r^*]}{2}$, we have that $\chi^{-1}(c)[Y_i,Y_j]$ is $(\delta,|\phi^*(ij)|^{-1})$-regular. In particular, $e_G(Y_i,Y_j) \geq (1-s\delta)|Y_i||Y_j|$.
\item\label{stabilityiii} Suppose $\sum_{i \in [r^*]}e(G[Y_i]) > \delta n^2$. Then
$\bm{k}$ does not have the strong extension property, and
all but at most $\delta n^2$ edges in $\bigcup_{i \in [r^*]}G[Y_i]$ are coloured with $1$ under $\chi$.
Moreover if $\bm{\ell} := (\ell_1,\ldots,\ell_{r^*}) \in \mathbb{N}^{r^*}$ is such that at least $\delta n^2$ edges need to be removed from $G[Y_i]$ to make it $K_{\ell_i}$-free, then
$\|\bm{\ell}\|_1 \leq k_{1}-1$.
\end{enumerate}
\end{theorem}

The following corollary is a simple consequence of Theorem~\ref{stabilitysimp} and shows that, when $\bm{k}$ has the strong extension property, every asymptotically extremal graph is close to complete partite.

\begin{corollary}[Stability for the strong extension property {\cite[Corollary~1.5]{stability2}}]\label{uniform}
Let $s \in \mathbb{N}$ and suppose that $\bm{k} \in \mathbb{N}^s$ with $k_1\geq\ldots \geq k_s$ has the strong extension property.
Then for all $\delta > 0$ there exist $n_0 \in \mathbb{N}$ and $\eps > 0$ such that the following holds.
Let $G$ be a graph on $n \geq n_0$ vertices such that
$$
\frac{\log F(G;\bm{k})}{\binom{n}{2}} \geq Q(\bm{k}) - \eps.
$$
Then there are $(r^*,\ba^*) \in \wt(\bm{k})$ and a partition V(G) = $V_1 \cup \ldots \cup V_{r^*}$ with $\left|\,|V_i|-\alpha^*_i n \,\right| < 1$ for all $i \in [r^*]$ such that
the number of adjacencies in $G$ that need to be changed to obtain $K[V_1,\ldots,V_{r^*}]$ is at most $\delta n^2$.
Moreover, for at least $(1-2^{-\eps n^2}) \cdot F(G;\bm{k})$ $\bm{k}$-valid $s$-edge-colourings $\chi$ of $G$, there exists $(r^*,\phi^*,\ba) \in \opt^*(\bm{k})$ such that $\|\ba-\ba^*\|_1 \leq \delta$ and $\chi^{-1}(c)[V_i,V_j]$ is $(\delta,|\phi^*(ij)|^{-1})$-regular for all $ij \in \binom{[r^*]}{2}$ and $c \in \phi^*(ij)$.
\end{corollary}

\subsection{A general exact result}

Our first main result is an exact version of Corollary~\ref{uniform}, that is, an exact result for $\bm{k}$ with the strong extension property. One consequence of this result is that, for such $\bm{k}$, every large $\bm{k}$-extremal graph is complete partite.

Indeed, we prove that every large $\bm{k}$-extremal graph $G^*$ is a complete $r^*$-partite graph, whose classes $X_1,\ldots, X_{r^*}$ are approximately $\ba^*$-weighted, for some $(r^*,\ba^*) \in \wt(\bm{k})$.
Moreover, almost every $\bm{k}$-valid colouring $\chi$ of $G^*$ is \emph{perfect (with respect to $(\phi;X_1,\ldots,X_{r^*})$)}, which means that there is a pattern $\phi \in \pat(r^*,\ba)$ 
for some $\ba$ close to $\ba^*$
so that every edge $e$ between $X_i$ and $X_j$ in $G^*$ satisfies $\chi(e) \in \phi(ij)$.

\begin{theorem}[Exactness for the strong extension property]\label{exact1}
Let $s \in \mathbb{N}$ and $\bm{k} \in \mathbb{N}^s$ have the strong extension property.
Then for all $\eps > 0$, there exist $\delta>0$ and $n_0 \in \mathbb{N}$ such that whenever $G$ is a $\bm{k}$-extremal graph on $n \geq n_0$ vertices, the following hold.
\begin{enumerate}[label=(\roman*),ref=(\roman*)]
\item\label{exact1i} $G$ is a complete multipartite graph: more precisely, there exist $(r^*,\ba^*) \in \wt(\bm{k})$ and $\bb \in \Delta^{r^*}$ with $\| \ba^*-\bb\|_1 < \eps$ such that $G = K(X_1,\ldots,X_{r^*})$, where $|X_i| = \beta_i n$ for all $i \in [r^*]$;
\item\label{exact1ii} for at least $(1-2^{-\delta n})\cdot F(G;\bm{k})$ $\bm{k}$-valid colourings $\chi$ of $G$, there is $(r^*,\ba) \in \wt(\bm{k})$ with $\|\ba^*-\ba\|_1 < \eps$ and $\phi \in \pat(r^*,\ba)$ such $\chi$ is perfect with respect to $(\phi;X_1,\ldots,X_{r^*})$.
\end{enumerate}
\end{theorem}

We show that, by solving a further optimisation problem (which is not too difficult in all known cases), one can determine $F(n;\bm{k})$ up to a multiplicative error of $1+o(1)$, and also determine the extremal graphs -- whenever the solution(s) of Problem $Q^*$ are known. Thus, for $\bm{k}$ with the strong extension property, it is true that `to solve the Erd\H{o}s-Rothschild problem, it suffices to solve Problem $Q^*$'.
Theorem~\ref{exact1} allows us to systematically recover all known exact results for the Erd\H{o}s-Rothschild problem.

\begin{problem}[Perfect colouring problem]\label{perfprob}
Let $(r^*,\ba^*) \in \wt(\bm{k})$ and $n \in \mathbb{N}$.
Maximise
$$
\perf_{r^*,\ba^*}(\bm{m}) := \sum_{\phi \in \pat(r^*,\ba^*)} \prod_{ij \in \binom{[r^*]}{2}}|\phi(ij)|^{m_im_j}
$$
subject to $\bm{m} \in P_{r^*}(n)$, where $P_{r^*}(n)$ is the set of $\bm{m} \in \mathbb{N}^{r^*}$ with $\|\bm{m}\|_1 = n$ and $m_1\geq \ldots \geq m_{r^*}$.
\end{problem}

Notice that $\perf_{r^*,\ba^*}(\bm{m})$ is essentially the number of perfect colourings of a complete $r^*$-partite graph $K_{\bm{m}}$ with vertex classes of size $m_1,\ldots,m_{r^*}$ (with some overcounting due to e.g.~distinct $\phi,\phi'$ with $\phi(ij)\cap \phi'(ij) \neq \emptyset$ for all $ij \in \binom{[r^*]}{2}$).
However, being able to solve the perfect colouring problem may not \emph{a priori} allow one to precisely determine the $\bm{k}$-extremal graphs.
Indeed, since not every colouring of $G$ is perfect, a graph with the maximal number of perfect colourings may not be $\bm{k}$-extremal if there is another graph with only slightly fewer.
For this reason, we make the following definition.

\begin{definition}[Solubility of $\bm{k}$]
\rm
Let $s \in \mathbb{N}$ and $\bm{k} \in \mathbb{N}^s$.
We say that $\bm{k}$ is \emph{soluble} if there exist $c > 0$ and $n_0 \in \mathbb{N}$, 
such that for all $n \geq n_0$ there are $(r^*,\ba^*) \in \wt(\bm{k})$ and $\bm{m}^* \in P_{r^*}(n)$ such that 
for all $(r',\ba') \in \wt(\bm{k})$ we have
$$
\perf_{r^*,\ba^*}(\bm{m}^*) > (1+c)\perf_{r',\ba'}(\bm{m})
$$
for all $\bm{m} \in P_{r'}(n)$ for which $(r^*,\ba^*,\bm{m}^*)$ is distinct from $(r',\ba',\bm{m})$.
We say that $\bm{m}^*$ is the \emph{supersolution} to $\bm{k}$ at $n$ and write $\perf(\bm{m}^*) := \perf_{r^*,\ba^*}(\bm{m}^*)$.
\end{definition}

Theorem~\ref{exact1}(i) implies that every $\bm{k}$-extremal graph $G$ on $n \geq n_0$ vertices is such that there is some $\bm{m} \in \mathbb{N}^{r^*}$ with $\|\bm{m}-\ba^*n\|_1 < \eps n$ so that $G \cong K_{m_1,\ldots,m_{r^*}}$.
If $\bm{k}$ is soluble, then Theorem~\ref{exact1}(ii) implies that for all large $n$ and all $m_1+\ldots+m_{r^*}=n$, we have $F(K_{m_1^*,\ldots,m_{r^*}^*},n) > F(K_{m_1,\ldots,m_{r^*}},n)$, where $\bm{m}^*$ is the supersolution at $n$ and $(m_1,\ldots,m_{r^*}) \neq (m_1^*,\ldots,m_{r^*}^*)$.
This gives the following corollary:

\begin{corollary}[Perfect exactness for the strong extension property]\label{exact2}
Let $s \in \mathbb{N}$ and let $\bm{k} \in \mathbb{N}^s$ be soluble and have the strong extension property.
Then 
$$
F(n;\bm{k}) = (1+o(1)) \cdot \perf(\bm{m}^*_n)\quad\text{as }n \to \infty
$$
and there exists $n_0 \in \mathbb{N}$ such that for all $n \geq n_0$,
the unique $\bm{k}$-extremal graph on $n$ vertices is $K_{\bm{m}_n^*}$, where $\bm{m}^*_n$ is the supersolution to $\bm{k}$ at $n$.
\end{corollary}

\begin{proposition}\label{recoverthm}
Every exact result in Theorem~\ref{knownthm} (that is, except the case $(3;5)$) follows from Corollary~\ref{exact2}.
\end{proposition}

\subsection{The triangle problem}

Our second main result concerns the case $\bm{k}=(3;s)$ (Erd\H{o}s and Rothschild's original problem was $s=2$).
%We use the shorthand $\bm{k}=(3;s)$.
Previously, $F(n;\bm{k})$ had been determined for $s=2,3,4,6$ and asymptotically for $s=5$.
We solve the problem for $s=7$, and use the same method to give a new proof for $s=6$.

\begin{theorem}\label{trianglethm}
There exists $n_0>0$ such that for every $n \geq n_0$ the Tur\'an graph $T_8(n)$ is the unique $(3;6)$-extremal graph, and the unique $(3;7)$-extremal graph.
\end{theorem}

In fact our results imply that
$$
F(n;(3;6)) = (C_j+o(1))\cdot (4^3 3^4)^{t_8(n)/7}\quad\text{and}\quad  F(n;(3;7)) = (C+o(1))\cdot 4^{t_8(n)}
$$
where $C$ is a constant and $C_j$ is a constant depending only on the remainder $j$ when $n$ is divided by $8$. These constants can be explicitly determined if desired.

A key component in our proof is our general exact result, Corollary~\ref{exact2}, which reduces the task to solving Problem $Q^*$.
For other cases of $\bm{k}$, Problem $Q^*$ has been solved by considering a linear relaxation, with variables which are essentially graph densities, and linear constraints which replace combinatorial constraints such as being $K_k$-free by the associated density bound given by Tur\'an's theorem.
However, the solutions to this relaxation generally do not correspond to feasible solutions to Problem $Q^*$, so additional constraints are needed.
The main new ingredient is a lemma about the density of the union of two dense triangle-free graphs (Lemma~\ref{RB}), which may be of independent interest.
This allows us to introduce a new constraint which yields a meaningful solution.

In the known cases of the triangle problem, the perfect colourings of extremal graphs are closely related to Hadamard matrices.
A \emph{Hadamard matrix} $H$ of order $n$ is an $n \times n$ matrix with entries in $\{ -1,+1 \}$ such that $HH^\intercal = nI_n$. 
A Hadamard matrix of order $n$ exists only if $n=1,2$ or $n$ is divisble by $4$; that this is sufficient was conjectured by Hadamard in 1893~\cite{hadamard}.
This conjecture remains open -- at the time of writing the smallest multiple of $4$ for which order there is no known Hadamard matrix is $4\times 167$.
It is however easy to construct arbitrarily large Hadamard matrices by using smaller ones as building blocks.
For example, Sylvester (see e.g.~\cite{sylvester}) observed that if $H$ is a Hadamard matrix, then so is $\left(\begin{smallmatrix}
H & H\\
H & -H
\end{smallmatrix}\right)$.
Hadamard matrices have the largest absolute value of the determinant among all complex square matrices with entries of absolute value at most one (see e.g.~\cite{hadamard}).

The connection to Problem $Q^*(3;s)$ is that there is a decomposition of the multigraph $2t K_{4t}$ into copies of $K_{2t,2t}$ if and only if a Hadamard matrix of order $4t$ exists (and the locations of the copies can be read off the matrix),
see~\cite{biclique} or Lemma~\ref{hadamard} for a proof.
It is plausible that in an optimal $(2r,\phi,\ba)$, every $\phi^{-1}(c) \cong K_{r,r}$, and $\ba$ is uniform 
(and also plausible that the number of vertices in an optimal solution is even).
If also $|\phi(ij)|=p$ for all pairs $ij$, then such a solution exists if and only if $p K_{2r}$ has a decomposition into copies of $K_{r,r}$. Comparing edge counts, the number of copies is $s=p\binom{2r}{2}/r^2=p(2r-1)/r$. So $p = tr$ for some integer $t$, since $r$ and $2r-1$ are coprime for $r \geq 2$. 
If $t=1$ then $r=2\ell$ for some integer $\ell$, since otherwise there is no decomposition due to the non-existence of the corresponding Hadamard matrix.
%$2rK_{2r}$ versus $
%Letting $\ba$ of length $r$ be uniform, we have
%$q(r,\phi,\ba) = \binom{r}{2}\cdot \frac{1}{r^2} \log(tr)$.
So $(4\ell,\phi_{4\ell-1},\bm{u}) \in \feas^*(3;4\ell-1)$, where $|\phi_{4\ell-1}(ij)|=2\ell$ for all $ij$ and $\bm{u}$ is uniform.
For $i=1,2,3$, $(4\ell,\phi_{4\ell-1}^{-i},\bm{u}) \in \feas^*(3;4\ell-1-i)$, where $\phi_{4\ell-1}^{-i}$ is obtained from $\phi_{4\ell-1}$ by removing $i$ colours.
Also, $(4(\ell-1),\phi_{4\ell-5}^{+4-i},\bm{u}) \in \feas^*(3;4\ell-1-i)$, where $\phi_{4\ell-5}^{+4-i}$ is obtained from $\phi_{4\ell-5}$ by duplicating $4-i$ colours.

This suggests the following problem, which would probably be very difficult to resolve.
%We dare not pose it as a conjecture given the dearth of evidence.

\begin{problem}
%Write $(3;s)$ for the length-$s$ tuple consisting of $3$'s.
Are the following true for sufficiently large $n$:
For all $\ell \geq 2$, is the unique $(3;4\ell-1)$-extremal graph $T_{4\ell}(n)$?
For all $\ell \geq 2$, is $T_{4(\ell-1)}(n)$ or $T_{4\ell}(n)$ one of the $(3;4\ell-1-i)$-extremal graphs for all $i=1,2,3$?
\end{problem}

All existing results fit the pattern described (see Subsection~\ref{hadamardsec}).

\subsection{The two colour problem}

Having proved a general exact result for the case when $\bm{k}$ has the strong extension property, we now investigate what happens when $\bm{k}$ has the extension property but not the strong extension property (we will say that such $\bm{k}$ have the \emph{weak extension property}). This appears to be much more difficult.
The reason for this is the possible presence of small parts in a complete multipartite extremal graph (when $\bm{k}$ has the strong extension property, no part has size $o(n)$ (see Lemma~\ref{lemma2}\ref{solutions})).
In a perfect colouring, only the neighbourhood into a large part is controlled and we have \emph{a priori} no information about colourings between small parts.
However, it could be the case that this part of the colouring is forced by the perfect colouring in the rest of the graph, and thus the number of small parts can also be controlled.

In our final theorem, we consider the simplest case when $\bm{k}$ has the weak extension property, namely $\bm{k} = (k+1,k)$ for $k \geq 3$.
For small $k$, we determine the (unique) $(k+1,k)$-extremal graph, which turns out to have a part of size $O(k)$, and the size of this part depends on the value of $n$ modulo $k-1$.
The proof is already fairly involved and relies heavily on a strong stability theorem for complete partite graphs (Theorem~\ref{stabilitycomp}), which is the main ingredient in the proof of Theorem~\ref{exact1}.
Similar arguments could probably be used to determine the $(k+\ell,k)$-extremal graph for (very) small $\ell$ and small $k$, which we discuss further in Section~\ref{concludesec}.

\begin{theorem}\label{2colapproxthm}
For all integers $k \geq 3$, there exists $n_0>0$ such that for all integers $n \geq n_0$, we have
$$
F(n;k+1,k) = O_k(1) \cdot 2^{t_{k-1}(n)},
$$
and every $(k+1,k)$-extremal graph is complete $k$-partite, with one part of size at most $2(k-1)$, and the other part sizes all within $2(k-1)$ of each other. Moreover, the constant $O_k(1)$ is at least $2$.
\end{theorem}

(Recall that the results of~\cite{abks} imply that $F(n;k+1,k) \geq F(n;k,k) = 2^{t_{k-1}(n)}$.)
We conjecture that the $O_k(1)$ multiplicative constant has a special form.

\begin{conjecture}\label{2colconj}
For all integers $k \geq 3$ and $0 \leq j \leq k-1$, there exists $n_0>0$ such that the following holds for all integers $n \geq n_0$. Let
$$
f(k,j,\ell) :=  \begin{cases} 
%2^{\ell(j-\ell)+\binom{j-\ell}{2}-\binom{j}{2}}
2^{-\binom{\ell}{2}}\left(k-1-\frac{j-\ell}{2}\right)^\ell &\mbox{if } 0 \leq \ell \leq j \\
%2^{-\ell(\ell-j)+\binom{\ell-j}{2}-\binom{j}{2}}
2^{j-\binom{\ell+1}{2}}\left(k-1+\ell-j\right)^\ell &\mbox{if } j < \ell \leq k-1
\end{cases}
$$
and let $\ell^*$ be such that $\max_{0 \leq \ell \leq k-1}f(k,j,\ell) = f(k,j,\ell^*)$.
Then $\ell^*$ is unique, and whenever $n \equiv j~(\!\!\!\!\mod k-1)$,
$$
F(n;k+1,k) = (1+o(1))\cdot f(k,j,\ell^*) \cdot 2^{t_{k-1}(n)}.
$$
Moreover the unique $(k+1,k)$-extremal graph on $n$ vertices is complete $k$-partite, with one part of size $\ell^*$ and the other parts as equal as possible.
\end{conjecture}

Given Theorem~\ref{2colapproxthm}, it is easy to prove this conjecture for small $k$, using computer assistance to test each of the possible extremal graphs.

\begin{theorem}\label{2colthm}
Conjecture~\ref{2colconj} holds for $3 \leq k \leq 10$.
\end{theorem}

For example, $F(n;4,3) = (2+o(1)) \cdot 2^{t_2(n)}$ when $n$ is even and $F(n;4,3) = (\frac{9}{4}+o(1))\cdot 2^{t_2(n)}$ when $n$ is odd; in both cases the unique extremal graph has parts of size $2,\lfloor\frac{n-2}{2}\rfloor,\lceil\frac{n-2}{2}\rceil$.

\subsection{Organisation of the paper}
In Section~\ref{exactsec}, we prove Theorem~\ref{exact1}, which follows from Theorem~\ref{stabilitycomp}, a strong stability theorem for complete partite graphs.
Section~\ref{trianglesec} solves Problem $Q^*$ in the case when $s=7$ and $k=3$.
Section~\ref{applicationsec} contains the proofs of Proposition~\ref{recoverthm} and Theorem~\ref{trianglethm}, which follow easily by combining Theorem~\ref{exact1} and solving Problem $Q^*$ in the relevant cases.
In Section~\ref{2colsec} we use the full strength of Theorem~\ref{stabilitycomp} to prove Theorems~\ref{2colapproxthm} and~\ref{2colthm}.
Section~\ref{concludesec} contains some concluding remarks and directions for future research.

There are five python programs used in the paper to check various numerical claims (\texttt{6check.py}, \texttt{6config.py}, \texttt{7ext.py}, \texttt{dcheck.py} and \texttt{smallpart.py}) which are all written in Python 3. They are included in the ancillary folder.

\section{A general exact result}\label{exactsec}

The aim of this section is to use Theorem~\ref{stabilitysimp} to prove Theorem~\ref{exact1}, a strengthening of Corollary~\ref{uniform} for $\bm{k}$ with the strong extension property.

%To obtain a completely \emph{exact} result -- i.e. determining precisely which graphs are $\bm{k}$-extremal, one must solve a further optimisation problem, which we describe.
%In some cases, we are able to solve this problem.

This section is organised as follows.
In Subsection~\ref{sectools} we collect some tools that we will need, including some lemmas on optimal solutions from~\cite{psy} and~\cite{stability2}, as well as standard regularity tools.
In Subsection~\ref{seccols} we define increasingly strict properties of colourings with respect to a partition.
The main result of Subsection~\ref{secstab} is Theorem~\ref{stabilitycomp}, a version of our stability result (Theorem~\ref{stabilitysimp}) for complete multipartite graphs.
The new part of this theorem is the last part which says that if every vertex has large contribution to the number of valid colourings, then there are many `perfect' colourings which follow a colour pattern exactly.
In Subsection~\ref{secexactproof} we prove Theorem~\ref{exact1}.

\subsection{Tools}\label{sectools}

The following are tools concerning Problem $Q^*$ from the previous papers in this series. Given a triple $(r,\phi,\ba)$, we have
$$
q(\phi,\ba) = \sum_{i \in [r]}\alpha_i q_i(\phi,\ba)\quad\text{where}\quad
q_i(\phi,\ba) := \sum_{\substack{j \in [r] \setminus \{ i \}\\ \phi(ij)\neq\emptyset}}\alpha_j\log|\phi(ij)|
$$
is the normalised contribution of vertex $i$ to the sum $q(\phi,\ba)$.

\begin{proposition}\label{continuitylagrange}
Let $s,r \in \mathbb{N}$ and $\bm{k} \in \mathbb{N}^s$.
The following hold.
\begin{enumerate}[label=(\roman*),ref=(\roman*)]
\item\cite[Proposition 11]{psy}\label{continuitylagrangei} Let $\phi \in \Phi_0(r;\bm{k})$ and $\ba,\bb \in \Delta^r$.
Then
$$
|q(\phi,\ba)-q(\phi,\bb)| < 2\left(\log s\right)\|\ba-\bb\|_1.
$$
\item\cite[Proposition~2.1]{stability2}\label{continuitylagrangeii} Let $(r,\phi,\ba) \in \opt_0(\bm{k})$.
For every $i \in [r]$ with $\alpha_i > 0$, we have that
$
q_i(\phi,\ba) = Q(\bm{k})
$.
\item\cite[Lemma~2.5]{stability2}\label{continuitylagrangeiii} Let $(r^*,\phi^*,\ba^*) \in \opt^*(\bm{k})$ and suppose $k_1 \geq \ldots \geq k_s$. Then $r^* \geq k_2-1$ and $\phi^{-1}(c)$ is maximally $K_{k_c}$-free for all $c \in [s]$.\qed
\end{enumerate}
\end{proposition}

\begin{lemma}\label{lemma2}
Let $s \in \mathbb{N}$ and suppose that $\bm{k} \in \mathbb{N}^s$ has the extension property, where $k_1 \geq \ldots \geq k_s$.
Then the following hold.
\begin{enumerate}[label=(\roman*),ref=(\roman*)]
\item\cite[Lemma~2.8]{stability2}\label{solutions} There exists $\mu >  0$ such that $\alpha_i^* > \mu$ for all $(r^*,\phi^*,\ba^*) \in \opt^*(\bm{k})$ and $i \in [r^*]$.
\item\cite[Lemma~2.9]{stability2}\label{clone}
There exists $\eps > 0$ such that the following holds.
Let $(r^*,\phi^*,\ba^*) \in \opt^*(\bm{k})$ and
let $\phi' \in \Phi_0(r^*+1,\bm{k})$ be such that $\phi'|_{\binom{[r^*]}{2}} = \phi^*$ and
$
\ext(\phi',\ba^*) \geq Q(\bm{k}) -\eps
$.
Then there exists $j \in [r^*]$ such that $r^*+1$ is a clone of $j$ under $\phi'$.
If $\bm{k}$ has the strong extension property then $r^*+1$ is a strong clone.
\item\cite[Lemma~2.10(iii)]{stability2}\label{strong}
Let $(r^*,\phi^*,\ba^*) \in \opt^*(\bm{k})$ and $\phi \in \Phi_0(r^*+1;\bm{k})$ be such that $\phi|_{\binom{[r^*]}{2}}=\phi^*$ and $r^*+1$ is a clone of $i \in [r^*]$ under $\phi$. Then $\phi(\{ i,r^*+1 \}) \subseteq \{ 1 \}$.\qed
\end{enumerate}
\end{lemma}

We need various lemmas concerning Szemer\'edi's regularity, starting with the following definitions.

\begin{definition}[Edge density, regularity of pairs and partitions]
\rm
Given a graph $G$ and disjoint non-empty sets $A,B \subseteq V(G)$, we define the \emph{edge density} between $A$ and $B$ to be
$$
d_G(A,B) := \frac{e_G(A,B)}{|A||B|},
$$
where $e_G(A,B)$ is the number of edges in $G$ with one vertex in $A$ and one vertex in $B$.
Given $\eps,d > 0$, the pair $(A,B)$ is called
\begin{itemize}
\item \emph{$\eps$-regular} if for every $X \subseteq A$ and $Y \subseteq B$ with $|X| \geq \eps|A|$ and $|Y| \geq \eps|B|$, we have that $|d_G(X,Y) - d_G(A,B)| \leq \eps$.
\item \emph{$(\eps, d)$-regular} if it is 
$\eps$-regular and $d_G(A,B) = d \pm \eps$, i.e.~$d-\eps\leq d_G(A,B)\leq d+\eps$,
\item \emph{$(\eps,\geq\! d)$-regular} if it is $\eps$-regular and $d_G(A,B) \geq d-\eps$.
\end{itemize}
\end{definition}

\begin{lemma}[Embedding Lemma {\cite[Theorem~2.1]{komsim}}]\label{embed}
For every $\eta >0$ and integer $k \geq 2$ there exist $\eps>0$ and $m_0 \in \mathbb{N}$ such that the following holds.
Suppose that $G$ is a graph with a partition $V(G) = V_1 \cup \ldots \cup V_k$ such that $|V_i| \geq m_0$ for all $i \in [k]$, and every pair $(V_i,V_j)$ for $1 \leq i < j \leq k$ is $(\eps,\geq \eta)$-regular.
Then $G$ contains $K_k$.\qed
\end{lemma}

\begin{proposition}[{\cite[Proposition~4.4]{stability2}}]\label{badred}
Let $A,B$ be disjoint sets of vertices, $s \in \mathbb{N}$ and $\eps > 0$ satisfying $1/|A|,1/|B| \ll \eps \ll 1/s$.
Let $G_1,\ldots,G_s$ be pairwise edge-disjoint subgraphs of $K[A,B]$.
Suppose that not all of $G_1,\ldots,G_s$ are $(\eps,s^{-1})$-regular graphs.
Then there exist $c \in [s]$ and $X \subseteq A$, $Y \subseteq B$ with $|X| = \lceil \eps|A|\rceil$ and $|Y| = \lceil \eps |B| \rceil$ such that 
$$
d_{G_c}(X,Y) \leq \frac{1}{s}\left( 1 - \frac{\eps}{2}\right).\qed
$$
\end{proposition}

\begin{proposition}\label{adjust}
Let $(A,B)$ be an $(\eps, d)$-regular pair and let $(A',B')$ be a pair such that $|A' \bigtriangleup A| \leq \alpha|A'|$ and $|B' \bigtriangleup B| \leq \alpha|B'|$ for some $0 \leq \alpha \leq 1$.
Then $(A',B')$ is an $(\eps+7\sqrt{\alpha}, d)$-regular pair.\qed
\end{proposition}

\begin{proposition}[see e.g.~\cite{komsim}]\label{badrefine}
Let $\eps,\delta$ be such that $0 < 2\delta \leq \eps < 1$.
Suppose that $(X,Y)$ is a $\delta$-regular pair, and let $X' \subseteq X$ and $Y' \subseteq Y$.
If
$$
\min \left\{ \frac{|X'|}{|X|} , \frac{|Y'|}{|Y|} \right\} \geq \frac{\delta}{\eps},
$$
then $(X',Y')$ is $\eps$-regular.\qed
\end{proposition}

The following estimate will be useful when counting colourings; it can be proved by looking at the tail of the binomial distribution.

\begin{proposition}[see e.g.~{\cite[Corollary~4.8]{stability2}}]\label{pain}
Let $n,k \in \mathbb{N}$  and $\delta \in \mathbb{R}$, where $0 < 1/n \ll \delta \ll 1/k$.
Then
$$
\sum_{i=0}^{\lfloor (k^{-1}-\delta) n \rfloor}\binom{n}{i} (k-1)^{n-i} \leq e^{-\delta^2k n/3} \cdot k^n.\qed
$$
\end{proposition}

Finally we need the following simple fact.

\begin{proposition}[see e.g.~{\cite[Claim~2.5.1]{stability2}}]\label{maxfree}
Let $k \in \mathbb{N}$ and let $H$ be maximally $K_k$-free. Then every $x \in V(H)$ lies in a copy of $K_{k-1}$ if and only if $|V(H)| \geq k-1$.\qed
\end{proposition}

\subsection{Hierarchy of colourings}\label{seccols}

We will now define three types of colouring, each stricter than the last, namely \emph{good}, \emph{locally good} and \emph{perfect}.
Each type is defined with respect to a partition of the graph.
Theorem~\ref{stabilitysimp} states that almost every valid colouring of a near-extremal graph $G$ is `good' with respect to a partition weighted like an optimal vertex weighting.
Our aim in Theorem~\ref{exact1} is to prove that almost all of these colourings are in fact `perfect'.
We achieve this via the property of being `locally good': where a colouring is such if, looking at a single colour between a pair in a partition, we see a regular graph of the right density, and furthermore, the coloured neighbourhood of a vertex or pair of vertices in every part has the right size.

Given a partition $X_1,\ldots,X_p$ of a set $S$, 
we say that a partition $Y_1, \ldots, Y_r$ of $S$ is a \emph{coarsening} of $X_1,\ldots, X_p$ if for all $i \in [p]$ there is a $j \in [r]$ such that $X_i \subseteq Y_j$.

\begin{definition}[$\delta$-good, ($\gamma,\delta)$-locally good, $(\gamma,\delta)$-perfect]\label{def:colourings}
\rm
Let $s \in \mathbb{N}$ and $\bm{k} \in \mathbb{N}^s$. Let $G$ be a complete $m$-partite graph with vertex partition $X_1,\ldots,X_m$. Let $\phi \in \Phi_0(r;\bm{k})$, let $Y_0,Y_1,\ldots,Y_r$ be a coarsening of $X_1,\ldots,X_m$ (where we allow $Y_0$ to be empty), let $\mathcal{Y} := (Y_1,\ldots,Y_{r})$, and let $\delta,\gamma > 0$.
We say that a $\bm{k}$-valid colouring $\chi$ of $G$ is:
\begin{itemize} 
\setlength\itemsep{0.75em}
\item \emph{$\delta$-good with respect to $(\phi; \mathcal{Y})$} if the following hold.
\begin{itemize}
\item For all $ij \in \binom{[r]}{2}$ and $c \in \phi(ij)$, we have that $\chi^{-1}(c)[Y_i,Y_j]$ is $(\delta,|\phi(ij)|^{-1})$-regular.
\item $
\sum_{i \in [r]}|\, e_G(Y_i) - |\chi^{-1}(1)[Y_i]| \,| < \delta n^2$.
\end{itemize} 
Write $\mathcal{G}_\delta(G;\phi,\mathcal{Y})$ for the set of these colourings.
\item \emph{$(\gamma,\delta)$-locally good with respect to $(\phi; \mathcal{Y})$} if the following hold.
\begin{itemize}
\item $\chi$ is $\delta$-good with respect to $(\phi;\mathcal{Y})$.
\end{itemize}
For all $x \in V(G)$, there exists $i=i_x \in [r]$ such that
\begin{itemize}
\item For all $j \in [r]\setminus \{ i \}$, parts $X \subseteq Y_j$ (so $X=X_t$ for some $t \in [m]$)
and $c \in \phi(ij)$, we have $|\chi^{-1}(c)[x,X]|=|\phi(ij)|^{-1}|X| \pm \frac{\delta}{\ell_j}|Y_j|$, where $\ell_j$ is the number of parts inside $Y_j$.  In particular $|\chi^{-1}(c)[x,Y_j]| = (|\phi(ij)|^{-1}\pm \delta)|Y_j|$.
\item $d_G(x,Y_i) - |\chi^{-1}(1)[x,Y_i]| < \delta n$. 
\item For all $i' \in [r]\setminus\{i\}$ we have $d_G(x,Y_{i'}) \geq (1-\gamma)|Y_{i'}|$.
\end{itemize}
For all distinct $y,z \in V(G)$, $j \in [r]\setminus\{i_y,i_z\}$ and parts $X \subseteq Y_j$ we have
\begin{itemize}
\item $|N_{\chi^{-1}(c_y)}(y,X) \cap N_{\chi^{-1}(c_z)}(z,X)| = |\phi(i_yj)|^{-1}|\phi(i_zj)|^{-1}|X| \pm \frac{\delta}{\ell_j}|Y_j|$ for all $c_y \in \phi(i_yj)$ and $c_z \in \phi(i_zj)$. 
\end{itemize}
Write $\mathcal{LG}_{\gamma,\delta}(G; \phi,\mathcal{Y})$ for the set of these colourings.
\item \emph{$(\gamma,\delta)$-perfect with respect to $(\phi; \mathcal{Y})$} if 
\begin{itemize}
\item $\chi$ is $(\gamma,\delta)$-locally good with respect to $(\phi;\mathcal{Y})$.
\item For all $x \in V(G)$, there exists $i \in [r]$ such that for all $j \in [r]$ and $y \in N_G(x,Y_j)$, we have $\chi(xy) = 1$ if $j=i$, and $\chi(xy) \in \phi(ij)$ otherwise.
%\item Moreover, if there exists $i' \in [r]$ such that $d_G(a,Y_{i'}) < (1-\gamma)|Y_{i'}|$, then $i = i'$.
\end{itemize}
Write $\mathcal{P}_{\gamma,\delta}(G;\phi,\mathcal{Y})$ for the set of these colourings.
\end{itemize}
\end{definition}

Note that if $\gamma' \geq \gamma$ and $\delta' \geq \delta$, then a $\delta$-good colouring is also $\delta'$-good; a $(\gamma,\delta)$-locally good colouring is also $(\gamma',\delta')$-locally good and a $(\gamma,\delta)$-perfect colouring is also $(\gamma',\delta')$-perfect.

It is a fairly straightforward consequence of Lemma~\ref{embed} (the Embedding Lemma) that, if a colouring is $(\gamma,\delta)$-locally good with respect to $(\phi;\mathcal{Y})$ and no member of $\mathcal{Y}$ is too small, then it is $(\gamma,\delta)$-perfect.

\begin{lemma}\label{perfect}
Let $s,r \in \mathbb{N}$ and $\bm{k} \in \mathbb{N}^s$, where $k_1 \geq \ldots \geq k_s$ and $r \geq k_2-1$, and let $\mu > 0$.
Then there exist $n_0 \in \mathbb{N}$ and $\delta > 0$ such that the following hold.
Let $G$ be a complete partite graph on $n \geq n_0$ vertices with parts $X_1,\ldots,X_m$, and let $Y_0,\ldots,Y_r$ be a partition of $V(G)$ which is a coarsening of $X_1,\ldots,X_m$.
Let $\mathcal{Y} := (Y_1,\ldots,Y_{r})$ and assume that $|Y_i| \geq \mu n$ for all $i \in [r]$.
Let $\phi \in \Phi_0(r;\bm{k})$ be such that $\phi^{-1}(c)$ is maximally $K_{k_c}$-free for all $c \in [s]$.
Then, for all $\gamma > 0$,
$$
\mathcal{LG}_{\gamma,\delta}(G;\phi,\mathcal{Y}) = \mathcal{P}_{\gamma,\delta}(G;\phi,\mathcal{Y}).
$$
\end{lemma}

\bpf
Let $\eps > 0,m_0 \in \mathbb{N}$ be the output of Lemma~\ref{embed} applied with $s^{-1},k_c$ playing the roles of $\eta,k$ for every $c \in [s]$.
Let $n_0 := 2s^2 m_0/\mu$ and $\delta := \eps^2$.
By decreasing $\eps$ and increasing $m_0$ if necessary, we may assume that $\delta \ll \mu$.

Certainly $\mathcal{LG}_{\gamma,\delta}(G;\phi,\mathcal{Y}) \supseteq \mathcal{P}_{\gamma,\delta}(G;\phi,\mathcal{Y})$ by definition.
Suppose that there exists $\chi \in \mathcal{LG}_{\gamma,\delta}(G;\phi,\mathcal{Y})$ such that the $(\gamma,\delta)$-perfect condition fails at some $x \in V(G)$.
Since $\chi$ is $(\gamma,\delta)$-locally good, there exists $i \in [r]$ such that for all $j \in [r]\setminus \{ i \}$, we have that $|\chi^{-1}(c)[x,Y_j]|=(|\phi(ij)|^{-1} \pm \delta)|Y_j|$ and $d_G(x,Y_i) - |\chi^{-1}(1)[x,Y_i]| < \delta n$.
We will show that for all $j \in [r]\setminus \{ i \}$ and $y \in N_G(x,Y_j)$ we have $\chi(xy) \in \phi(ij)$, and for all $y \in N_G(x,Y_i)$ we have $\chi(xy) = 1$.

Suppose first that there exists $j \in [r]\setminus \{ i \}$ and $y \in N_G(x,Y_j)$ such that $\chi(xy) =: c \notin \phi(ij)$.
Let $k := k_{c}$ and $J := \chi^{-1}(c)$.
Now, $\phi^{-1}(c)$ is maximally $K_{k}$-free.
Since $ij \notin \phi^{-1}(c)$, we have that $\phi^{-1}(c) \cup \{ ij \}$ contains a copy of $K_{k}$.
So there exist $i_3,\ldots,i_k \in [r]\setminus \{ i,j \}$ such that $i_1,\ldots,i_k$ span a copy of $K_k$ in $\phi^{-1}(c)$, where $i_1 := i$ and $i_2 := j$.

Since $\chi$ is $(\gamma,\delta)$-locally good, for all $\ell=3,\ldots,k$, we have, by the pairs condition, taking the union over all parts $X \subseteq Y_{i_\ell}$, that
$$
|N_J(x,Y_{i_\ell}) \cap N_J(y,Y_{i_\ell})| \geq (|\phi(ii_\ell)|^{-1}|\phi(ji_\ell)|^{-1} - \delta)|Y_{i_\ell}| \geq \frac{|Y_{i_\ell}|}{2s^2} \geq \frac{\mu n}{2s^2} \geq m_0.
$$
Let $U_\ell := N_J(x,Y_{i_\ell}) \cap N_J(y,Y_{i_\ell})$.
Proposition~\ref{badred} implies that $J[U_\ell,U_{\ell'}]$ is a $(\sqrt{\delta},\geq\! s^{-1})$-regular pair for all distinct $\ell,\ell' \in \{ 3,\ldots, k \}$.
Lemma~\ref{embed} implies that $J$ contains a copy of $K_{k-2}$.
Together with $x,y$, we obtain a copy of $K_k$ in $\chi^{-1}(c)$, a contradiction.

Suppose instead that there is some $y \in N_G(x,Y_i)$ such that $\chi(xy) =: c \neq 1$.
Let $k := k_c$ and $J := \chi^{-1}(c)$.
Since $r \geq k_2-1 \geq k-1$, Proposition~\ref{maxfree} implies that $i$ lies in a copy of $K_{k-1}$ in the graph $\phi^{-1}(c)$.
Let $i_1 := i$ and let $i_2,\ldots,i_{k-1}$ be the other vertices in this copy.
%Then each of $d_J(x,Y_{i_\ell}),d_J(y,Y_{i_\ell})=(|\phi(i_1i_\ell)|^{-1}\pm\delta)|Y_{i_\ell}|$ for all $\ell \in [2,k-1]$.
As before, setting $U_\ell := N_J(x,Y_{i_\ell}) \cap N_J(y,Y_{i_\ell})$ for all $\ell \in [2,k-1]$, we have that $J[U_2,\ldots,U_{k-1}]$ contains a copy of $K_{k-2}$.
Together with $x,y$, this gives a copy of $K_k$ in $J$, a contradiction.
%Now, it remains to prove that, if there is $i' \in [r]$ such that $x \in Y_{i'}$ and $d_G(x,Y_{i'}) < (1-\gamma)|Y_{i'}|$, then $i=i'$.
%But this follows from the fact that $\chi$ is $(\gamma,\delta)$-locally good.
\epf

\subsection{Stability for complete multipartite graphs}\label{secstab}

%Let %$s \in \mathbb{N}$, 
Given a graph $G$, a subgraph $H$ of $G$, and an $s$-edge-colouring $\chi$ of $H$,
we say that $\overline{\chi}$ is an \emph{extension} of $\chi$ if $\overline{\chi}$ is an $s$-edge-colouring of $G$ such that $\overline{\chi}|_{H} = \chi$.

\begin{theorem}[Stability for complete multipartite graphs]\label{stabilitycomp}
Let $s \in \mathbb{N}$ and suppose that $\bm{k} \in \mathbb{N}^s$ with $k_1 \geq \ldots \geq k_s$ has the extension property.
Then for all $\delta > 0$ there exist $n_0 \in \mathbb{N}$ and $\eps \in \mathbb{R}$ with $0 < \eps < \delta$ such that the following holds.
Let $m \in \mathbb{N}$ and $G = K(X_1,\ldots,X_m)$ be a complete $m$-partite graph on $n \geq n_0$ vertices such that
$$
\frac{\log F(G;\textbf{k})}{\binom{n}{2}} \geq Q(\textbf{k}) - \eps.
$$
Then for at least $(1-2^{-\eps n^2}) \cdot F(G;\bm{k})$ $s$-edge-colourings $\chi$ of $G$ which are $\bm{k}$-valid, there are $(r^*,\phi^*,\ba^*) \in \opt^*(\bm{k})$ with $r^* \leq m$ and a coarsening $Z_0,\ldots,Z_{r^*}$ of $X_1,\ldots,X_m$ such that the following hold.
\begin{enumerate}[label=(\roman*),ref=(\roman*)]
\item\label{stabilityci}
$
\sum_{i \in [r^*]}|\,|Z_i|-\alpha_i^*n\,| < \delta n$,
and $\chi$ is $\delta$-good with respect to $(\phi^*; Z_1,\ldots,Z_{r^*})$.
%Furthermore, if $|X_i| \leq \delta^2n$, then $X_i \subseteq Z_0$. 
\item\label{stabilitycii} If, for all $i \in [r^*]$, $\ell_i$ is the number of $X_j$ in $Z_i$, then $\ell_1 = \ldots = \ell_{r^*} = 1$; or $\bm{k}$ does not have the strong extension property and $\ell_1 + \ldots + \ell_{r^*} \leq k_1-1$.
\end{enumerate}
Furthermore, if we have
\begin{equation}
\begin{aligned}\label{def}
&\log F(G;\bm{k}) - \log F(G-x;\bm{k}) \geq (Q(\bm{k}) - 2\eps)n \quad\forall \ x \in V(G), \ \text{ and} \\
&\log F(G;\bm{k}) - \log F(G-y-z;\bm{k}) \geq (Q(\bm{k}) - 2\eps)(n+n-1) \quad\forall\text{ distinct }y,z \in V(G),
\end{aligned}
\end{equation}
then there are at least $(1-2^{-\eps n})\cdot F(G;\bm{k})$ $s$-edge-colourings $\chi$ of $G$ for which there are $(r^*,\phi^*,\ba^*) \in \opt^*(\bm{k})$ and $Z_0,\ldots,Z_{r^*}$ as above such that $\chi$ is $(0,\delta)$-perfect with respect to $(\phi^*; Z_1,\ldots,Z_{r^*})$.
\end{theorem}

%\bpf
{\it Proof.}
For brevity, write $Q := Q(\bm{k})$, $R := R(\bm{k})$, $F(G) := F(G;\bm{k})$, and abbreviate similarly elsewhere.
Lemma~\ref{lemma2}\ref{solutions} implies that there exists $\mu > 0$ such that for all $(r^*,\ba^*) \in \wt(\bm{k})$, we have that $\alpha_i^* > \mu$ for all $i \in [r^*]$.
Applying Lemma~\ref{lemma2}\ref{clone} gives the constant $\eps_0$.
We may assume that $\delta \ll \eps_0,\mu,1/R$ since a smaller $\delta$ yields a stronger conclusion.
Apply Lemma~\ref{embed} with $(2s)^{-1},k_c$ playing the roles of $\eta,k$ for all $c \in [s]$ to obtain $m_0 \in \mathbb{N}$ and $\eps' > 0$ such that its conclusions hold.
Choose constants $\delta_0,\ldots,\delta_5 \in \mathbb{R}$ such that $0 < \delta_0 \ll \ldots \ll \delta_5 \ll \delta$.
We may assume that $\eps' \geq 3\sqrt{\delta}_1$.
Apply Theorem~\ref{stabilitysimp} with $\delta_0$ playing the role of $\delta$ to obtain $n_0 \in \mathbb{N}$ and $\eps > 0$.
Without loss of generality, we may assume that $0 < 1/n_0 \ll \eps \ll \delta_0$ and
\begin{equation}\label{m0}
\sqrt{\delta_1}\mu n_0 \geq 2m_0.
\end{equation}
We may further assume that the conclusion of Lemma~\ref{perfect} holds with $s,\bm{k},\mu/2$ playing the roles of inputs $s,\bm{k},\mu$; and with outputs $n_0 \in \mathbb{N}, \delta_5 > 0$.
Thus our constants form the hierarchy
\begin{equation}\label{hier2}
0 < 1/n_0 \ll \eps \ll \delta_0 \ll \ldots \ll \delta_5 \ll \delta \ll \eps_0,\mu,1/R.
\end{equation}
Let $G = K(X_1,\ldots,X_m)$ be a complete $m$-partite graph on $n \geq n_0$ vertices such that $\log F(G) \geq (Q-\eps)\binom{n}{2}$.
Note that $m < R$ (or $F(G)=0$).
Let also
\begin{equation}\label{IZ0}
Y_0 := \bigcup_{i \in [m] : |X_i| \leq \delta_1^4 n}X_i.
\end{equation}
Given $(r^*,\ba^*) \in \wt(\bm{k})$, let $\coars(\ba^*)$ be the set of partitions $\mathcal{Y} := (Y_1,\ldots,Y_{r^*})$ of $V(G)\setminus Y_0$ such that $Y_0,\ldots,Y_{r^*}$ is a coarsening of $X_1,\ldots,X_m$, and
$$
\sum_{i \in [r^*]}|\,|Y_i|-\alpha_i^*n\,| < \delta_1^2n.
$$
(We do not yet know that $\coars(\ba^*)$ is non-empty.)

Given $\eta>0$, let $\mathcal{G}_{\eta}(G)$ be the set of $\bm{k}$-valid colourings $\chi$ of $G$ for which there exist $(r^*,\phi^*,\ba^*) \in \opt^*(\bm{k})$ and $\mathcal{Y} \in \coars(\ba^*)$ such that $\chi \in \mathcal{G}_{\eta}(G;\phi^*,\mathcal{Y})$.
Call the elements of $\mathcal{G}_{\eta}(G)$ \emph{$\eta$-good}.
Given $\eta_1,\eta_2$, let $\mathcal{LG}_{\eta_1,\eta_2}(G)$ and $\mathcal{P}_{\eta_1,\eta_2}(G)$ be as $\mathcal{G}_{\eta}(G)$ but with $\mathcal{LG}_{\eta_1,\eta_2}(G;\phi^*,\mathcal{Y})$ and $\mathcal{P}_{\eta_1,\eta_2}(G;\phi^*,\mathcal{Y})$ replacing $\mathcal{G}_{\eta}(G;\phi^*,\mathcal{Y})$ respectively.

\begin{claim}\label{delta1}
We have the following:
\begin{enumerate}[label=(\roman*),ref=(\roman*)]
\item\label{delta1i} $|\mathcal{G}_{\delta_1}(G)| \geq (1-2^{-\eps n^2})\cdot F(G)$.
\item\label{delta1ii} For all $(r^*,\phi^*,\ba^*) \in \opt^*(\bm{k})$, $\mathcal{Y} \in \coars(\ba^*)$ and $\chi \in \mathcal{G}_{\delta_1}(G;\phi^*,\mathcal{Y})$, let $t_i$ be the number of $X_j$ in $Y_i$ for all $i \in [r^*]$.
Then either $t_1 = \ldots = t_{r^*} = 1$; or $\bm{k}$ does not have the strong extension property and $t_1+\ldots+t_{r^*} \leq k_1-1$.
\end{enumerate}
\end{claim}

\bcpf
By Theorem~\ref{stabilitysimp} applied to $G$ with parameter $\delta_0$, we have that there are at least $(1-2^{-\eps n^2}) \cdot F(G)$ colourings $\chi : E(G) \rightarrow [s]$ which are $\bm{k}$-valid and for which there is some $(r^*,\phi^*,\ba^*) \in \opt^*(\bm{k})$ and a partition $V_1 \cup \ldots \cup V_{r^*} = V(G)$ such that:
 \begin{itemize}
\item For all $i \in [r^*]$ we have $|\,|V_i|-\alpha^*_i n\,| \leq 1$.
\item $\chi$ is $\delta_0$-good with respect to $(\phi^*;V_1,\ldots,V_{r^*})$.
\item Suppose $\sum_{i \in [r^*]}e(G[V_i]) > \delta_0 n^2$. Then
$\bm{k}$ does not have the strong extension property, and
all but at most $\delta_0 n^2$ edges in $\bigcup_{i \in [r^*]}G[V_i]$ are coloured with $1$ under $\chi$.
Moreover if $\bm{\ell} := (\ell_1,\ldots,\ell_{r^*}) \in \mathbb{N}^{r^*}$ is such that at least $\delta_0|V_i|^2$ edges need to be removed from $G[V_i]$ to make it $K_{\ell_i}$-free, then
$\|\bm{\ell}\|_1 \leq k_{1}-1$.
\end{itemize}
We will show that every such $\chi$ lies in $\mathcal{G}_{\delta_1}(G)$, which implies the first part of the claim.
Fix $\chi$ and its associated $(r^*,\phi^*,\ba^*)$ and $V_1,\ldots, V_{r^*}$ (recall that both $(r^*,\phi^*,\ba^*)$ and the partition $V_1\cup\ldots\cup V_{r^*}$ may be different for different $\chi$).
For all $ij \in \binom{[r^*]}{2}$, the $\delta_0$-goodness of $\chi$ implies that
$$
e_G(V_i,V_j) \geq \sum_{c \in \phi^*(ij)}|\chi^{-1}(c)[V_i,V_j]| \geq \sum_{c \in \phi^*(ij)}(|\phi^*(ij)|^{-1} - \delta_0)|V_i||V_j| \geq (1-s\delta_0)|V_i||V_j|,
$$
so $e_{\overline{G}}(V_i,V_j) \leq s\delta_0n^2$.
Suppose that $X_k$, some $k \in [m]$, is such that $|X_k \cap V_i|> \sqrt{s\delta_0}n$.
Then, for every $j \in [m]\setminus\{k\}$,
$$
|X_k \cap V_j| = \frac{e_{\overline{G}[X_k]}(V_i,V_j)}{|X_k \cap V_i|} \leq \frac{e_{\overline{G}}(V_i,V_j)}{|X_k \cap V_i|} < \sqrt{s\delta_0}n. 
$$
So for all $k \in [m]$ with $|X_k| > \delta_1^4n > R\sqrt{s\delta_0}n$, there is a unique $i_k \in [r^*]$ such that $|X_k \cap V_{i_k}| > \sqrt{s\delta_0}n$.
For all $i \in [r^*]$, let
$$
Y_i := \bigcup_{k \in [m] : i_k = i}X_k.
$$
Recall that we already defined $Y_0$ in~(\ref{IZ0}).
Then $Y_0,\ldots,Y_{r^*}$ is a coarsening of $X_1,\ldots,X_m$, and $|X_j| \leq \delta_1^4 n$ if and only if $X_j \subseteq Y_0$.
For all $i \in [r^*]$ we have that
\begin{align}\label{YiVi}
|Y_i \bigtriangleup V_i| &\leq |Y_0| + \sum_{k \in [m]: i_k \neq i}|X_k \cap V_i| + \sum_{k \in [m]: i_k = i}|X_k \setminus V_i| \leq m\delta_1^4n + 2Rm\sqrt{s\delta_0}n \leq \delta_1^3 n
\end{align}
and
\begin{equation}\label{Yisize}
\sum_{i \in [r^*]}|\,|Y_i| - \alpha^*_in \,| \leq \sum_{i \in [r^*]}\left(|\, |V_i| - |Y_i|\, | + |\,|V_i|-\alpha^*_i n\,| \right) \leq R\delta_1^3 n + R < \delta_1^{2} n,
\end{equation}
so $(Y_1,\ldots,Y_{r^*}) \in \coars(\bm{\alpha^*})$.
Moreover, for all $i \in [r^*]$ we have
\begin{equation}\label{Yiagain}
|Y_i| \geq (\alpha^*_i - \delta^2_1)n-1 \geq \frac{\mu n}{2}.
\end{equation}
It remains to prove that $\chi$ is $\delta_1$-good with respect to $(\phi^*;Y_1,\ldots,Y_{r^*})$.
Proposition~\ref{adjust} and~(\ref{YiVi}) imply that for all $ij \in \binom{[r^*]}{2}$ and all $c \in \phi^*(ij)$, we have that $\chi^{-1}(c)[Y_i,Y_j]$ is a $(\delta_1,|\phi^*(ij)|^{-1})$-regular pair.
Moreover,
\begin{eqnarray*}
&\phantom{=}& \sum_{i \in [r^*]}|\, e_G(Y_i) - |\chi^{-1}(1)[Y_i]| \,| \leq \sum_{i \in [r^*]}|\, e_G(V_i) - |\chi^{-1}(1)[V_i]| \,| + \sum_{i \in [r^*]}|Y_i \triangle V_i|^2\\
&\stackrel{(\ref{YiVi})}{\leq}& \delta_0 n^2 + R\delta_1^6 n^2 < \delta_1^5 n^2.
\end{eqnarray*}
Thus $\chi \in \mathcal{G}_{\delta_1}(G)$.
This completes the proof of the first part of the claim.

For the second part, 
fix $(r^*,\phi^*,\ba^*) \in \opt^*(\bm{k})$, $\mathcal{Y} \in \coars(\ba^*)$ and $\chi \in \mathcal{G}_{\delta_1}(G;\phi^*,\mathcal{Y})$. 
For each $i \in [r^*]$, let $t_i$ be the number of parts $X_j$ which lie in $Y_i$.
Then $t_i$ is at most the number of parts $X_j$ which have intersection at least $\sqrt{s\delta_0}n$ with $V_i$.
Suppose first that $\sum_{i \in [r^*]}e_G(V_i) \leq \delta_0 n^2$.
If $t_i \geq 2$ for some $i \in [r^*]$, then $e_G(V_i) \geq s\delta_0 n^2$, a contradiction.
So $t_1 = \ldots = t_{r^*} = 1$.

Suppose now that $\sum_{i \in [r^*]}e_G(V_i)> \delta_0 n^2$.
Then by Theorem~\ref{stabilitysimp}\ref{stabilityiii}, $\bm{k}$ does not have the strong extension property.
Since $G[V_i]$ is a complete multipartite graph containing at least $t_i$ parts of size at least $\sqrt{s\delta_0}n$, we have that at least $s\delta_0 n^2 > \delta_0|V_i|^2$ edges need to be removed from $G[V_i]$ to make it $K_{t_i}$-free.
Thus $t_1+\ldots+t_{r^*} \leq k_1-1$, proving the second part of the claim.
\ecpf

This proves parts~\ref{stabilityci} and~\ref{stabilitycii} of Theorem~\ref{stabilitycomp}.
Namely, take any $\chi \in \mathcal{G}_{\delta_1}(G) \neq \emptyset$.
Let $(r^*,\phi^*,\ba^*) \in \opt^*(\bm{k})$ and $\mathcal{Y} := (Z_0,\ldots,Z_{r^*}) \in \coars(\ba^*)$ witness $\chi \in \mathcal{G}_{\delta_1}(G;\phi^*,\mathcal{Y})$.
Then they satisfy items~\ref{stabilityci} and~\ref{stabilitycii} of Theorem~\ref{stabilitycomp}.

Suppose now that~(\ref{def}) holds, but
\begin{equation}\label{LG}
|\mathcal{G}_{\delta_1}(G) \setminus \mathcal{LG}_{\delta_4,\delta_5}(G)| > 2^{-\delta_2 n} \cdot F(G).
\end{equation}
For most of the next part of the proof, we will establish a contradiction to this assumption.
(Recall that, by Lemma~\ref{perfect}, a direct contradiction to the statement of the theorem would replace $\mathcal{LG}_{\delta_4,\delta_5}(G)$ by $\mathcal{LG}_{0,\delta_5}(G)$.)

For every $(r^*,\phi^*,\ba^*) \in \opt^*(\bm{k})$, Lemma~\ref{continuitylagrange}\ref{continuitylagrangeiii} implies that $(\phi^*)^{-1}(c)$ is maximally $K_{k_c}$-free for all $c \in [s]$, and that $r^* \geq k_2-1$.
Lemma~\ref{perfect} and~(\ref{Yiagain}) imply that for every $\mathcal{Y} \in \coars(\ba^*)$, we have that $\mathcal{P}_{\delta_4,\delta_5}(G;\phi^*,\mathcal{Y}) = \mathcal{LG}_{\delta_4,\delta_5}(G;\phi^*,\mathcal{Y})$.
So
\begin{equation}\label{perfectcol}
\mathcal{LG}_{\delta_4,\delta_5}(G) = \mathcal{P}_{\delta_4,\delta_5}(G).
\end{equation}
Given $\mathcal{Y} = (Y_1,\ldots,Y_{r^*}) \in \coars(\ba^*)$ and $x \in V(G)$, write $\mathcal{Y}-x$ to denote the partition $(Y_1\setminus\{x\},\ldots,Y_{r^*}\setminus\{x\})$ and define $\mathcal{Y}-y-z$ similarly.

\begin{claim}\label{goodcol}
At least one of the following hold.

There exist $x \in V(G)$ and a $\bm{k}$-valid $s$-edge-colouring $\chi$ of $G-x$ such that the following two statements hold.
\begin{enumerate}[label=(\roman*),ref=(\roman*)]
\item\label{goodcoli} $\chi \in \mathcal{G}_{2\delta_1}(G-x;\phi^*,\mathcal{Y}-x)$, for some $(r^*,\phi^*,\ba^*) \in \opt^*(\bm{k})$ and $\mathcal{Y} := (Y_1,\ldots,Y_{r^*}) \in \coars(\ba^*)$.
\item\label{goodcolii} There is a set $\mathrm{Ext}(\chi)$ of at least $2^{(Q-\delta_3)n}$ $\bm{k}$-valid extensions $\overline{\chi}$ of $\chi$ to $G$ such that $\overline{\chi} \in \mathcal{G}_{\delta_1}(G;\phi^*,\mathcal{Y})$ but $\overline{\chi}$ is not $(\delta_4,\delta_5)$-locally good with respect to $(\phi^*;Y_1,\ldots,Y_{r^*})$ at $x$.
\end{enumerate}

There exist $zy \in \binom{V(G)}{2}$ and a $\bm{k}$-valid $s$-edge-colouring $\xi$ of $G-z-y$ such that the following two statements hold.
\begin{enumerate}[label=(\roman*),ref=(\roman*)]
\setcounter{enumi}{2}
\item\label{goodcoliii} $\xi \in \mathcal{G}_{2\delta_1}(G-z-y;\phi^*,\mathcal{Y}-z-y)$, for some $(r^*,\phi^*,\ba^*) \in \opt^*(\bm{k})$ and $\mathcal{Y} := (Y_1,\ldots,Y_{r^*}) \in \coars(\ba^*)$.
\item\label{goodcoliv} There is a set $\mathrm{Ext}(\xi)$ of at least $2^{(Q-\delta_3)(n+n-1)}$ $\bm{k}$-valid extensions $\overline{\xi}$ of $\xi$ to $G$ such that $\overline{\xi} \in \mathcal{G}_{\delta_1}(G;\phi^*,\mathcal{Y})$ but $\overline{\xi}$ is not $(\delta_4,\delta_5)$-locally good with respect to $(\phi^*;Y_1,\ldots,Y_{r^*})$ at $z,y$.
\end{enumerate}
\end{claim}

\bcpf
Suppose first that for at least half of the colourings in $\mathcal{G}_{\delta_1}(G)\setminus \mathcal{LG}_{\delta_4,\delta_5}(G)$, the locally good condition fails at some vertex (rather than only at pairs).
We will show that there exists an $x$ satisfying~\ref{goodcoli} and~\ref{goodcolii}.
By~(\ref{LG}) and~(\ref{perfectcol}), there is some 
%$\phi \in \textsc{Pat}(r^*, \ba^*)$ and 
$x \in V(G)$ such that there are at least $\frac{1}{2}\cdot\frac{1}{n} \cdot 2^{-\delta_2 n} \cdot F(G) \geq 2^{-2\delta_2 n} \cdot F(G)$
valid colourings $\chi$ which are $\delta_1$-good (with respect to some optimal solution and partition) but for which the $(\delta_4,\delta_5)$-locally good condition fails at $x$ (for all optimal solutions and partitions).
Call this set of colourings $\mathcal{L}_x(G)$.

%Note that $s^{-n}\cdot F(G) \leq F(G-x) \leq F(G)$. 
We have $\log F(G-x) \geq \log (s^{-n}\cdot F(G)) \geq (Q-\eps)\binom{n}{2} - n\log s \geq (Q-2\eps)\binom{n-1}{2}$,
so a version of Claim~\ref{delta1}\ref{delta1i} applied to $G-x$ implies that
$$
|\mathcal{G}_{2\delta_1}(G-x)| \geq (1-2^{-2\eps (n-1)^2})\cdot F(G-x).
$$

%Note that every colouring of $G-x$ which is not $2\delta_1$-good is such that no extension to $G$ is $\delta_1$-good.
%Therefore Claim~\ref{delta1}\ref{delta1i} implies that
%$$
%F(G-x) - |\mathcal{G}_{2\delta_1}(G-x)| \leq F(G) - |\mathcal{G}_{\delta_1}(G)| \leq 2^{-\eps n^2} F(G).
%$$
Suppose that $x$ does not satisfy Claim~\ref{goodcol}.
Then for each $\bm{k}$-valid colouring $\chi$ of $G-x$, either it does not lie in $\mathcal{G}_{2\delta_1}(G-x)$, or it does lie in $\mathcal{G}_{2\delta_1}(G-x)$ but has at most $2^{(Q-\delta_3)n}$ $\bm{k}$-valid extensions which lie in $\mathcal{L}_x(G)$.
We have that
\begin{eqnarray*}
2^{-2\delta_2n}\cdot F(G) &\leq& |\mathcal{L}_x(G)| \leq |\mathcal{G}_{2\delta_1}(G-x)| \cdot 2^{(Q-\delta_3)n} + (F(G-x)-|\mathcal{G}_{2\delta_1}(G-x)|) \cdot s^n\\
&\leq& F(G-x)\cdot 2^{(Q-\delta_3)n} + s^{n} \cdot 2^{-2\eps (n-1)^2}F(G-x)\\
&\leq& %F(G-x)(2^{(Q-\delta_3)n}+s^{2n}\cdot 2^{-\eps n^2}) \leq 
F(G-x)\cdot 2^{(Q-\delta_3/2)n},
\end{eqnarray*}
and so
$$
\log F(G) - \log F(G-x) \leq (Q - \delta_3/2 + 2\delta_2)n \leq (Q-\delta_3/3)n < (Q-\eps)n,
$$
a contradiction to~(\ref{def}).

Thus we may assume that at least half the good but not locally good colourings fail due to the pair condition.
Again there is a pair $y,z$ of distinct vertices that appear in at least $2^{-2\delta_2 n}\cdot F(G)$ such colourings,
and an identical argument gives the required $\xi$, satisfying~\ref{goodcoliii} and~\ref{goodcoliv}.
\ecpf

\medskip
\noindent
Suppose there exist $x,\chi$ as in Claim~\ref{goodcol}\ref{goodcoli} and~\ref{goodcolii}.
So there are $(r^*,\phi^*,\ba^*) \in \opt^*(\bm{k})$ and $\mathcal{Y} = (Y_1,\ldots,Y_{r^*}) \in \coars(\ba^*)$ and $\mathrm{Ext}(\chi)$ as in Claim~\ref{goodcol}.
Define $\bb := (\beta_1,\ldots,\beta_{r^*})$ by setting $\beta_i := |Y_i|/n$ for all $i \in [r^*]$.
Then~(\ref{Yisize}) implies that
\begin{equation}\label{alphabeta}
\|\bb - \ba^*\|_1 < \delta_1^2.
\end{equation}
%For each $\chi \in \mathcal{G}(G)$, let $(r^*,\phi^*,\ba^*)$ be \emph{the associated optimal solution of $\chi$} and let $Z_1 \cup \ldots \cup Z_{r^*}$ be the \emph{associated partition of $\chi$}.
Note that $\mathrm{Ext}(\chi) \subseteq \mathcal{G}_{3\delta_1}(G;\phi^*,\mathcal{Y})$.
For every $\overline{\chi} \in \mathrm{Ext}(\chi)$,
define $\phi = \phi(\overline{\chi}) : \binom{[r^*+1]}{2} \rightarrow 2^{[s]}$ by setting
$$
\phi(ij) := \begin{cases} 
\phi^*(ij) &\mbox{if }  ij \in \binom{[r^*]}{2}; \\
\{ c \in [s] : |\overline{\chi}^{-1}(c)[x,Y_j]| \geq \sqrt{\delta_1} |Y_{j}| \} &\mbox{if } i = r^*+1, j \in [r^*]. \\
\end{cases}
$$
%(So if $Y_{i_0}$ is a single part of $G$, then $\phi^*(\{ i_0,r^*+1 \}) = \emptyset$.)
Fix the pattern $\phi$ that appears for the largest number of
extensions $\overline{\chi} \in \mathrm{Ext}(\chi)$, and let $\mathrm{Ext}_{\phi}(\chi)$ be the set of these $\overline{\chi}$.
By Claim~\ref{goodcol}\ref{goodcolii},
\begin{equation}\label{ext'}
|\mathrm{Ext}_{\phi}(\chi)| \geq 2^{-sr^*} \cdot |\mathrm{Ext}(\chi)| \geq 2^{(Q-2\delta_3)n}.
\end{equation}

\begin{claim}\label{claimphi1}
 $\phi \in \Phi_0(r^*+1,\bm{k})$.\end{claim}

\bcpf
Suppose not.
Then there is some $c \in [s]$ such that $\phi^{-1}(c)$ contains a copy of $K_{k_c}$.
%Since $\phi$ is an extension of $\phi^* \in \Phi(r^*,\bm{k})$, the vertex $r^*+1$ lies in this copy.
Since $\phi^* \in \Phi_2(r^*;\bm{k})$ is the restriction of $\phi$ to $\binom{[r^*]}{2}$, the vertex $r^*+1$ must lie in this copy.
Let $z_1,\ldots,z_{k_c-1}$ be the other vertices.
Let $\overline{\chi} \in \mathrm{Ext}_{\phi}(\chi)$ be arbitrary.
For each $k \in [k_c-1]$, let $Z_k := \{ y \in Y_{z_k} : \overline{\chi}(xy) = c \}$.
Then, by the definition of $\phi$, we have
\begin{equation}\label{Zk}
|Z_k| > \sqrt{\delta_1}|Y_{z_k}| \stackrel{(\ref{m0})}{\geq} m_0.
\end{equation}
Let $kk' \in \binom{[k_c-1]}{2}$.
Now, $\overline{\chi}^{-1}(c)[Y_{z_k},Y_{z_{k'}}]$ is $(3\delta_1,\geq\!\! s^{-1})$-regular since $\overline{\chi} \in \mathrm{Ext}(\chi) \subseteq \mathcal{G}_{3\delta_1}(G; \phi^*,\mathcal{Y})$. Therefore
$$
d(\overline{\chi}^{-1}(c)[Z_k,Z_{k'}]) \geq d(\chi^{-1}(c)[Y_{z_k},Y_{z_{k'}}]) - 3\delta_1 \geq |\phi^*(z_kz_{k'})|^{-1} - 6\delta_1 \geq 1/(2s). 
$$
Now Proposition~\ref{badrefine} implies that $\overline{\chi}^{-1}(c)[Z_k,Z_{k'}]$ is $3\sqrt{\delta_1}$-regular.
Therefore $\overline{\chi}^{-1}(c)[Z_k,Z_{k'}]$ is $(3\sqrt{\delta_1},\geq\!\!1/(2s))$-regular.
Now~(\ref{Zk}) and Lemma~\ref{embed} imply that $\chi^{-1}(c)[Z_1,\ldots,Z_{k_c-1}]$ contains a copy of $K_{k_c-1}$.
Together with $x$, this gives a $c$-coloured copy of $K_{k_c}$ in $\overline{\chi}$, contradicting the fact that $\overline{\chi}$ is $\bm{k}$-valid. 
\ecpf

\medskip
\noindent
Next we will show that there is some $i \in [r^*]$ such that $r^*+1$ is a clone of $i$ under $\phi$.
Suppose for a contradiction that this is not the case.
Thus by our choice of $\eps_0$ we have $\ext(\phi,\ba^*) \leq Q-\eps_0$.
Then, using (a version of) Proposition~\ref{continuitylagrange}\ref{continuitylagrangei} and Lemma~\ref{lemma2}\ref{clone},
\begin{eqnarray}
\label{extcount} 
\frac{1}{n}\log\prod_{j \in [r^*]}|\phi(\{j,r^*+1\})|^{|Y_j|} &=& \ext(\phi,\bb) \leq \ext(\phi,\ba^*) + \|\ba^*-\bb\|_1 \cdot \log s\\
\nonumber &\stackrel{(\ref{alphabeta})}{<}& Q - \eps_0 + \delta_1^2\log s \leq Q - \eps_0/2.
\end{eqnarray}

The number of possible patterns $\phi$ on $r^*+1$ vertices which are extensions of $\phi^*$ is at most $2^{sr^*}$, and, by definition, the number of edges with endpoint $x$ which are not coloured according to $\phi$ is at most $s\sqrt{\delta_1}n$.
Therefore the \emph{total} number of $\bm{k}$-valid extensions $\overline{\chi} \in \mathrm{Ext}(\chi)$ of $\chi$ to $G$ is at most
\begin{equation}\label{total}
2^{sr^*} \cdot \binom{n}{\leq s\sqrt{\delta_1}n} \cdot s^{s\sqrt{\delta_1}n} \cdot s^{|Y_0|} \cdot \prod_{j \in [r^*]} |\phi(\{ j,r^*+1 \})|^{|Y_j|} \stackrel{(\ref{IZ0}),(\ref{extcount})}{\leq} 2^{(Q-\eps_0/3)n} < 2^{(Q-2\delta_3)n}.
\end{equation}
So certainly there are at most this number of extensions which lie in $\mathrm{Ext}_{\phi}(\chi)$, contradicting~(\ref{ext'}).

Therefore we may assume that there is some $i^* \in [r^*]$ such that $r^*+1$ is a clone of $i^*$ under $\phi$.
Thus we have
\begin{align}\label{phi'}
\phi(\{ j,r^*+1 \}) = \phi^*(i^*j) \ \text{ for all } \ j \in [r^*]\setminus \{ i^*\} \ \text{ and } \ \phi(\{ i^*,r^*+1 \}) \subseteq \{ 1 \},
\end{align}
where the last inclusion follows from Lemma~\ref{lemma2}\ref{strong}.

\begin{claim}\label{badcol}
We have the following for all $\overline{\chi} \in \mathrm{Ext}_\phi(\chi)$:
\begin{enumerate}[label=(\roman*),ref=(\roman*)]
\item\label{badcoli} There exist $j^* \in [r^*]\setminus \{ i^* \}$, a part $X \subseteq Y_{j^*}$ and $c \in \phi^*(i^*j^*)$ such that $|\overline{\chi}^{-1}(c)[x,X]| \neq |\phi(i^*j^*)|^{-1}|X| \pm \frac{\delta_5}{\ell_{j^*}}|Y_{j^*}|$.
\item\label{badcolii} If there exists $h \in [r^*]$ such that the part $X$ of $G$ containing $x$ lies in $Y_{h}$, and $|X| > \delta_4|Y_{h}|$, then $i^* = h$.
\end{enumerate}
\end{claim}

\bcpf
Define $\ell \in \{ 0,\ldots, r^* \}$ as follows.
If there is $i' \in [r^*]$ such that $d_G(x,Y_{i'}) < (1-\delta_4)|Y_{i'}|$ (noting that there can be at most one such $i'$), let $\ell := i'$.
If there is no such $i'$, let $\ell :=0$.
Let $X$ be the part of $G$ which contains $x$.
%If $\ell \neq 0$, then $d_G(x,Y_\ell) < (1-\delta_4)|Y_\ell|$.
%So $|X| \geq \delta_4 |Y_\ell| \geq \delta_4^2 n$.

\medskip
\noindent
\textbf{Case 1:}
\emph{$\ell = 0$.}

\medskip
\noindent
%Then, for all $j \in [r^*]$, we have that $d_G(x,Y_j) > (1-\delta_4)|Y_j| > s\sqrt{\delta_1}|Y_j|$ and so $\phi(\{ j,r^*+1\}) \neq \emptyset$.
%Then~(\ref{phi'}) implies that
%$\phi(\{ i^*,r+1 \}) = \{ 1 \}$.
Suppose that~(i) is false.
Since $\overline{\chi}$ is not $(\delta_4,\delta_5)$-locally good at $x$, we must have that
$d_G(x,Y_{i^*}) - |\overline{\chi}^{-1}(1)[x,Y_{i^*}]| > \delta_5 |Y_{i^*}|$.
By~(\ref{phi'}) and the definition of $\phi$, we have that 
\begin{align}\label{Yi*} 
d_G(x,Y_{i^*}) - |\overline{\chi}^{-1}(1)[x,Y_{i^*}]| &= \sum_{c \in \{ 2,\ldots, s \}}|\overline{\chi}^{-1}(c)[x,Y_{i^*}]| < (s-1)\sqrt{\delta_1}|Y_{i^*}| < \delta_5 |Y_{i^*}|,
\end{align}
a contradiction.
For~(ii), let $X \subseteq Y_h$ be the part of $G$ containing $x$ and suppose $h \in [r^*]$. Then $|Y_h|-|X|=d_G(x,Y_h) \geq (1-\delta_4)|Y_h|$. So in this case, $|X| \leq \delta_4|Y_h|$ and~(ii) is vacuous.

\medskip
\noindent
\textbf{Case 2:}
\emph{$\ell \in [r^*]$.}

\medskip
\noindent
We will first show that (ii) holds.
Note that $h=\ell$ since $X \subseteq Y_h$ and $x\in X$ is not adjacent to the whole of $Y_\ell$.
If $Y_{h} = X$, then $N_G(x,Y_{h}) = \emptyset$ so $\phi(\{ r^*+1,h \}) = \emptyset$.
So~(\ref{phi'}) implies that $h = i^*$, as required.

Otherwise, $Y_{h} \neq X$.
Suppose that $i^* \neq h$.
Then $|\phi^*(i^*h)| \geq 2$, and so
\begin{align*}
\frac{1}{n}\left( \log \prod_{j \in [r^*]} |\phi(\{ j,r^*+1 \})|^{d_G(x,Y_j)}\right) &\stackrel{(\ref{phi'})}{=} \frac{1}{n}\sum_{j \in [r^*]\setminus \{ i^* \}} |Y_j\setminus X| \log|\phi^*(i^*j)|\\
&= \sum_{j \in [r^*]\setminus \{ i^*,h \}}\beta_j \log|\phi^*(i^*j)| + \left(\beta_h -\frac{|X|}{n}\right) \log|\phi^*(i^*h)|\\
%&\leq \sum_{j \in [r^*]\setminus \{ i^* \}}\alpha^*_j\log|\phi^*(i^*j)| 
&\leq q_{i^*}(\phi^*,\ba^*) + \left( \|\ba^*-\bb\|_1 - \frac{|X|}{n}\right)\log s\\
&\stackrel{(\ref{alphabeta})}{\leq} Q + (\delta_1^2 - \delta_4^2)\log s < Q - 4\delta_3,
\end{align*}
where the penultimate inequality follows from Proposition~\ref{continuitylagrange}\ref{continuitylagrangeii}.
Combining this with a very similar calculation to~(\ref{total}) implies that the total number of $\bm{k}$-valid extensions $\overline{\chi} \in \mathrm{Ext}(\chi)$ of $\chi$ to $G$ is less than $2^{(Q-2\delta_3)n}$, contradicting~(\ref{ext'}).
Therefore $i^*=h$, proving (ii).
The proof of (i) now proceeds exactly as in Case~1.
\ecpf

\medskip
\noindent
Claim~\ref{badcol}\ref{badcoli} implies that for every $\overline{\chi} \in \mathrm{Ext}_{\phi}(\chi)$, there exist $j^* \in [r^*]\setminus\{i^*\}$, a part $X \subseteq Y_{j^*}$ and $c^* \in \phi^*(i^*j^*)$
such that, by averaging,
\begin{equation}\label{need}
|\overline{\chi}^{-1}(c^*)[x,X]| \leq |\phi^*(i^*j^*)|^{-1}|X|-\frac{\delta_5}{s\ell_{j^*}}|Y_{j^*}|.
\end{equation}
In particular, since $\ell_{j^*} \leq k_1 < R$, we have
$|X| \geq \delta_5|Y_{j^*}|/(s \cdot R) \geq (\delta_5)^2|Y_{j^*}|$.
The number of ways of adding edges coloured with $\phi^*$ between $x$ and $X$ with this reduced density 
(choosing the set of colour-$c^*$ neighbours, and then colouring every other edge with any available colour other than $c^*$)
is at most
\begin{align*}
&\phantom{=} \sum_{k=0}^{\lfloor |\phi^*(i^*j^*)|^{-1}|X|-\frac{\delta_5}{s\ell_{j^*}}|Y_{j^*}|\rfloor} \binom{|X|}{k}(|\phi^*(i^*j^*)|-1)^{|X|-k}\\
&\leq  e^{-\delta_5^2|X|/3(s^2R^2)} \cdot |\phi^*(i^*j^*)|^{|X|} \stackrel{(\ref{Yiagain})}{\leq} e^{-\delta_5^5 n} \cdot |\phi^*(i^*j^*)|^{|X|},
\end{align*}
where the first inequality is a consequence of Lemma~\ref{pain}.
Using Proposition~\ref{continuitylagrange}\ref{continuitylagrangeii}, we have that
\begin{eqnarray}\label{qextend}
%\log \left( \prod_{j \in [r^*]\setminus \{ i^* \}} |\phi^*(\{ j,r^*+1 \})|^{|Y_j|} \right)
q_{i^*}(\phi^*,\bb) %n\sum_{j \in [r^*]\setminus \{ i^* \}}\beta_j\log|\phi^*(i^*j)| 
\leq Q + \|\ba^*-\bb\|_1 (2\log s)  \stackrel{(\ref{alphabeta})}{\leq} Q+2\delta_1^2\log s.
\end{eqnarray}
Now, we can generate every $\overline{\chi} \in \mathrm{Ext}_{\phi}(\chi)$ from $\chi$ by doing the following.
For each $j \in [r^*]$, choose at most $\sqrt{\delta_1}|Y_j|$ vertices $y \in N_G(x,Y_j)$ and colour $xy$ arbitrarily.
Arbitrarily colour all edges $xy_0$ where $y_0 \in Y_0$.
Then choose $j^*,X,c^*$ as above and colour every uncoloured edge (with endpoint $x$) according to $\phi$ so that~(\ref{need}) holds (with a small adjustment to account for the $\sqrt{\delta_1}|Y_j|$ edges that have already been coloured).
That this will indeed generate every $\overline{\chi} \in \mathrm{Ext}_{\phi}(\chi)$ is a consequence of the definition of $\phi$ and Claim~\ref{badcol}.
Therefore
\begin{eqnarray*}
|\mathrm{Ext}_{\phi}(\chi)| &\leq& \binom{n}{\leq \sqrt{\delta_1}n} \cdot s^{\sqrt{\delta_1}n} \cdot s^{|Y_0|} \cdot e^{-\delta_5^5 n} \cdot \prod_{j \in [r^*]\setminus \{ i^*\}} |\phi^*(i^*j)|^{|Y_j|}\\
&\stackrel{(\ref{IZ0}),(\ref{qextend})}{\leq}& 2^{-\delta_5^6 n} \cdot 2^{(Q+2\delta_1^2 \log s)n} < 2^{(Q-2\delta_3)n},
\end{eqnarray*}
contradicting~(\ref{ext'}).
So our assumption that there are $x,\chi$ as in Claim~\ref{goodcol}\ref{goodcoli} and~\ref{goodcolii} are false.

Thus there must be $z,y,\xi$ such that Claim~\ref{goodcol}\ref{goodcoliii} and~\ref{goodcoliv} hold.
So there are $(r^*,\phi^*,\ba^*) \in \opt^*(\bm{k})$ and $\mathcal{Y} = (Y_1,\ldots,Y_{r^*}) \in \coars(\ba^*)$ and $\mathrm{Ext}(\xi)$ as in Claim~\ref{goodcol}.
Again we define colour patterns $\phi_y, \phi_z$ in analogy with $\phi$ to be the pair of patterns that appear together for the greatest number of extensions $\overline{\xi} \in {\rm Ext}(\xi)$, and write ${\rm Ext}_{\phi_y,\phi_z}(\xi)$ for the set of extensions $\overline{\xi} \in {\rm Ext}(\xi)$ with this pair of patterns. Again, by Claim~\ref{goodcol}\ref{goodcoliii},
\begin{equation}\label{ext'2}
|\mathrm{Ext}_{\phi_y,\phi_z}(\xi)| \geq 2^{(Q-2\delta_3)(n+n-1)}.
\end{equation} 
Similarly to Claim~\ref{claimphi1}, we have $\phi_y,\phi_z \in \Phi_0(r^*+1;\bm{k})$.
Applying Lemma~\ref{lemma2}\ref{clone}, we have that there are $i_y, i_z \in [r^*]$ such that in $\phi_y,\phi_z$ respectively, the vertex $r^*+1$ is a clone of $i_y,i_z$. And, in particular, $\phi_y(\{i_y,r^*+1\}) \subseteq \{1\}$ and $\phi_z(\{i_z,r^*+1\}) \subseteq \{1\}$.
The analogue of Claim~\ref{badcol}\ref{badcoli} for pairs holds in the same way, that is,
there exist $j^* \in [r^*]\setminus\{i_y,i_z\}$, a part $X \subseteq Y_{j^*}$, $c_y \in \phi_y(i_yj^*)$ and $c_z \in \phi_z(i_zj^*)$ such that
$$
|N_{\xi^{-1}(c_y)}(y,X) \cap N_{\xi^{-1}(c_z)}(z,X)| \neq |\phi^*(i_yj^*)|^{-1}|\phi^*(i_zj^*)|^{-1}|X| \pm \frac{\delta_5}{\ell_{j^*}}|Y_{j^*}|.
$$
%\item[(ii)] If there exists $h_y \in [r^*]$ such that the part $X_y$ of $G$ containing $y$ lies in $Y_{h_y}$, and $|X_y| > \delta_4|Y_{h_y}|$, then $i_y = h_y$. Similarly for $z$.
It remains to show that this implies there are few extensions of $\chi$ to $G$.
By averaging, there are $c_y^* \in \phi^*(i_yj^*)$ and $c_z^* \in \phi^*(i_zj^*)$ such that
$$
|N_{\xi^{-1}(c_y^*)}(y,X) \cap N_{\xi^{-1}(c_z^*)}(z,X)| < |\phi^*(i_yj^*)|^{-1}|\phi^*(i_zj^*)|^{-1}|X| - \frac{\delta_5}{s^2\ell_{j^*}}|Y_{j^*}|.
$$
The number of ways of adding edges coloured with $\phi^*$ between $y,z$ and $X$ with this reduced density (choosing the set of $v$ such that $yv$ is coloured $c_y^*$ and $zv$ is coloured $c_z^*$, and then colouring every other $yu,zu$ with any pair other than $c_y^*,c_z^*$) is at most
\begin{align*}
&\phantom{=} \sum_{k=0}^{\lfloor |\phi^*(i_yj^*)|^{-1}|\phi^*(i_zj^*)|^{-1}|X|-\frac{\delta_5}{s^2\ell_{j^*}}|Y_{j^*}|\rfloor} \binom{|X|}{k}(|\phi^*(i_yj^*)||\phi^*(i_zj^*)|-1)^{|X|-k} \cdot s\\
%\leq  e^{-\delta_5^2|Y_{j^*}|/2s^2} \cdot (|\phi^*(i_yj^*)||\phi^*(i_zj^*)|)^{|Y_{j^*}|} 
&\hspace{4cm}\stackrel{(\ref{Yiagain})}{\leq} e^{-\delta_5^5 n} \cdot (|\phi^*(i_yj^*)||\phi^*(i_zj^*)|)^{|X|}.
\end{align*}
Repeating the calculations in the single vertex case, we see that $|{\rm Ext}_{\phi_y,\phi_z}(\xi)| < 2^{(Q-2\delta_3)(n+n-1)}$, contradicting~{\ref{ext'2}).

Thus our assumption~(\ref{LG}) is false.
Therefore, combining its negation with Claim~\ref{delta1}, we see that
\begin{eqnarray*}
|\mathcal{P}_{\delta_4,\delta_5}(G)| \geq |\mathcal{G}_{\delta_1}(G)| - 2^{-\delta_2 n} \cdot F(G) \geq (1-2^{-\eps n}) \cdot F(G).
\end{eqnarray*}

Let $\chi \in \mathcal{P}_{\delta_4,\delta_5}(G)$.
Then there exists $(r^*,\phi^*,\ba^*) \in \opt^*(\bm{k})$ and $\mathcal{Y} = (Y_1,\ldots,Y_{r^*}) \in \coars(\ba^*)$ such that $\chi$ is $(\delta_4,\delta_5)$-perfect with respect to $(\phi^*;Y_1,\ldots,Y_{r^*})$.
So for all $x \in V(G) \setminus Y_0$, there exists $i(x) \in [r^*]$ such that for all $j \in [r^*]$ and $b \in N_G(x,Y_j)$, we have $\chi(xb) = 1$ if $j = i(x)$; and $\chi(xb) \in \phi^*(i(x)j)$ otherwise.
Moreover, (the proof of) the second part of Claim~\ref{badcol} implies that, if there exists $h(x) \in [r^*]$ such that the part $X$ of $G$ containing $x$ lies in $Y_{h(x)}$, and $|X| > \delta_4|Y_{h(x)}|$, then $i(x) = h(x)$.

Now, define a new partition $Z_0,\ldots,Z_{r^*}$ of $V(G)$ by setting, for all $i \in [r^*]$,
$$
Z_i := \{ x \in V(G) : i(x) = i \text{ and }x \in X \subseteq Y_{h(x)} \text{ with }|X| \geq \delta_4|Y_{h(x)}| \}
\quad\text{and}\quad Z_0 := V(G) \setminus \textstyle\bigcup_{i \in [r^*]}Z_i.
$$
Then $Z_0,\ldots,Z_{r^*}$ is a coarsening of $X_1,\ldots,X_m$, and $\chi$ is $(0,\delta_5)$-perfect with respect to $(\phi^*;Z_1,\ldots,Z_{r^*})$ by definition.
Moreover, 
for all $i \in [r^*]$, $Z_i \subseteq Y_i$, and $|Y_i \setminus Z_i| \leq R\delta_4|Y_i|$.
So
$$
\sum_{i \in [r^*]}|\, |Z_i| - \alpha^*_in\, | \leq \sum_{i \in [r^*]}(|Y_i| - |Z_i|) + \sum_{i \in [r^*]}|\, |Y_i| - \alpha^*_in \,| \leq R\delta_4 n + \delta_1^2 n < \delta n, 
$$
as required.
Finally, let $\ell_i$ be the number of $X_j$ in $Z_i$.
Then $\ell_i \leq t_i$ for all $i \in [r^*]$.
The second part of Claim~\ref{delta1} yields the desired conclusion.
%\epf
\hfill $\qed$
\medskip

Note that the statement of Theorem~\ref{stabilitycomp} can be made much simpler in the case when $\bm{k}$ has the strong extension property as in this case we have that $(r^*,\ba^*) \in \wt(\bm{k})$ and the partitions $Z_0,\ldots,Z_{r^*}$ are identical for all of the at least $(1-2^{-\eps n^2})F(G;\bm{k})$ colourings specified in the theorem.
Indeed, suppose that $(r^*,\ba^*) \in \wt(\bm{k})$ and its associated partition $Z_0,\ldots,Z_{r^*}$ are outputs of Theorem~\ref{stabilitycomp} for some specified colouring.
Then each $Z_1,\ldots,Z_{r^*}$ is a part of $G$ by Theorem~\ref{stabilitycomp}(ii).
Part (i) implies that for each $i \in [r^*]$ we have that $|Z_i| \geq \alpha^*_i n - \delta n \geq \mu n/2$.
Furthermore, $|Z_0| \leq \delta n \ll \mu n/2$.
So, provided $\delta$ is chosen to be smaller than $\mu/2$, the structure of $G$ itself determines $(r^*,\ba^*)$ and $Z_0,\ldots,Z_{r^*}$.
Thus we have the following corollary, which will be used to prove Theorem~\ref{exact1}.

\begin{corollary}\label{stabilitycomp2}
Let $s \in \mathbb{N}$ and suppose that $\bm{k} \in \mathbb{N}^s$ has the strong extension property.
Then for all $\delta > 0$ there exist $n_0 \in \mathbb{N}$ and $\eps \in \mathbb{R}$ with $0 < \eps < \delta$ such that the following holds.
Let $m \in \mathbb{N}$ and $G = K(X_1,\ldots,X_m)$ be a complete $m$-partite graph on $n \geq n_0$ vertices such that $|X_1| \geq \ldots \geq |X_m|$ and
$$
\frac{\log F(G;\textbf{k})}{\binom{n}{2}} \geq Q(\textbf{k}) - \eps.
$$
Then the following hold.
\begin{enumerate}[label=(\roman*),ref=(\roman*)]
\item\label{sci} There is $(r^*,\ba^*) \in \wt(\bm{k})$ with $r^* \leq m$ such that
$
\sum_{i \in [r^*]}|\,|X_i|-\alpha_i^*n\,| < \delta n
$.
%Moreover, $|\mathcal{G}_{\delta_1}(G;\phi^*, X_1,\ldots,X_{r^*})| \geq (1-2^{-\eps n^2})\cdot F(G;\bm{k})$.
\item\label{scii} For at least
$(1-2^{-\eps n^2}) \cdot F(G;\bm{k})$ $s$-edge-colourings $\chi$ of $G$ which are $\bm{k}$-valid there is $\phi^* \in \pat(\ba^*;\bm{k})$ such that
$\chi$ is $\delta$-good with respect to $(\phi^*; X_1,\ldots,X_{r^*})$.

\item\label{sciii} Furthermore, if we have
\begin{equation}
\begin{aligned}\label{def2}
&\log F(G;\bm{k}) - \log F(G-x;\bm{k}) \geq (Q(\bm{k}) - 2\eps)n \quad\forall \ x \in V(G), \ \text{ and} \\
&\log F(G;\bm{k}) - \log F(G-y-z;\bm{k}) \geq (Q(\bm{k}) - 2\eps)(n+n-1) \quad\forall\text{ distinct }y,z \in V(G),
\end{aligned}
\end{equation}
%then $|\mathcal{P}_{0,\delta}(G;\phi^*,X_1,\ldots,X_{r^*})| \geq (1-2^{-\eps n})\cdot F(G;\bm{k})$. 
then for at least $(1-2^{-\eps n})\cdot F(G;\bm{k})$ valid $s$-edge-colourings $\chi$ of $G$ there is $(r^*,\phi^*,\ba) \in \opt^*(\bm{k})$ with $\|\ba-\ba^*\|_1 \leq \delta$ such that $\chi$ is $(0,\delta)$-perfect with respect to $(\phi^*; X_1,\ldots,X_{r^*})$.\qed
\end{enumerate}
\end{corollary}

\subsection{The proof of Theorem~\ref{exact1}}\label{secexactproof}

The next observation is a simple but key ingredient of our proof, which allows us to only consider complete multipartite graphs.
If there were a $\bm{k}$-extremal graph $H$ which is not complete multipartite, one can use symmetrisation to obtain from $H$ a new graph $H'$ which is complete multipartite (Theorem~1 in~\cite{psy}).
By \emph{symmetrisation}, we mean replacing a vertex $u$ with a copy, or \emph{twin}, of $v \notin N_H(u)$.
Crucially, we can do this in such a way that we end up with a part containing a single vertex
(which is connected to every other vertex).

\begin{lemma}\label{small}
Let $s \in \mathbb{N}$ and $\bm{k} \in \mathbb{N}^s$.
Let $G$ be a $\bm{k}$-extremal graph which is not complete multipartite.
Then there exists a $\bm{k}$-extremal graph $G'$ which is complete multipartite and has a part of size one.
\end{lemma}

\bpf
Since $G$ is not complete multipartite, there exist distinct non-adjacent vertices $u,v \in V(G)$ such that $N_G(u) \neq N_G(v)$.
For any graph $H$, let $\chi(H)$ be the set of $\bm{k}$-valid colourings of $H$.
For each $\chi \in \chi(G-u-v)$, let $\chi_u, \chi_v$ denote the number of valid extensions of $\chi$ to $G-v$ and $G-u$ respectively. 
Since $u$ and $v$ are non-adjacent,
$$
F(G;\bm{k}) = \sum_{\chi \in \chi(G-u-v)}\chi_u\chi_v.
$$
Let $G_u$ denote the graph obtained from $G$ by replacing $v$ by a twin of $u$.
Define $G_v$ similarly.
The operation of passing from $G$ to $G_u$ or $G_v$ is a \emph{symmetrisation}.
We have
$$
F(G_u;\bm{k}) = \sum_{\chi \in \chi(G-u-v)}\chi_u^2 \ \text{ and } \ F(G_v;\bm{k}) = \sum_{\chi \in \chi(G-u-v)}\chi_v^2.
$$
Then
$$
0 \geq F(G_u;\bm{k}) + F(G_v;\bm{k}) - 2F(G;\bm{k}) = \sum_{\chi \in \chi(G-u-v)} (\chi_u - \chi_v)^2 \geq 0. 
$$
Therefore $G_u$ and $G_v$ are both $\bm{k}$-extremal.

Let $\mathcal{G}$ be the directed graph whose vertex set contains all $n$-vertex graphs (up to isomorphism) that can be obtained from $G$ by a sequence of symmetrisations, and add a directed edge from $H$ to $H' \neq H$ if $H'\cong H_u$ for some vertex $u$.
Note that $H$ has outdegree equal to $0$ in $\mathcal{G}$ if and only if $H$ is complete multipartite.
By \cite[Theorem~3]{psy} there is at least one sequence of symmetrisations which leads to a ($\bm{k}$-extremal) complete partite graph $H$; among all choices pick one such that the number $m$ of parts in $H$ is as small as possible.
We may assume that every part has size at least two, or we are done.
Let $H^-$ be an inneighbour of $H$ in $\mathcal{G}$.
Observe that $H^-$ is not complete multipartite since it does not have $0$ outdegree.
Then there exists $x \in V(H^-)$ such that $H^--x$ is a complete $m$-partite graph with parts $V_1,\ldots,V_m$.

\begin{claim}
$x$ has at least one neighbour in each of $V_1,\ldots, V_m$ in $H^-$.
\end{claim}

\bcpf
Suppose for a contradiction that $x$ does not have a neighbour in $V_1$, say.
Since $H^-$ is not complete multipartite, without loss of generality, there is $y \in V_2$ such that $xy \notin E(H^-)$.
Replace every $u \in V_2 \setminus \{ y \}$ with a twin of $y$ to obtain a graph $J$ with vertex partition $\{ x \}, V_1,\ldots, V_m$, such that $J-x \cong K(V_1,\ldots,V_m)$, and $xz \notin E(J)$ for all $z \in V_1 \cup V_2$.
This is a sequence of symmetrisations, so there is an oriented path from $H^-$ to $J$ in $\mathcal{G}$ (and in particular $J \in V(\mathcal{G})$).
Now, given $v_1 \in V_1$ and $v_2 \in V_2$, we have $xv_1,xv_2 \notin E(J)$ but $v_1v_2 \in E(J)$.
Therefore we can replace $V_1 \cup V_2$ with a set $X$ of $|V_1|+|V_2|$ twins of $x$ to obtain a new graph $J'$, which has vertex partition $X, V_3,\ldots,V_m$, and $J'-X \cong K(V_3,\ldots,V_m)$.
Again $J' \in V(\mathcal{G})$.

Suppose that $x \in V(J)$ is such that $xw \notin E(J)$ for some $w \in V_3 \cup \ldots \cup V_m$.
Then $xw \notin E(J')$ for all $x \in X$.
Replace every $x \in X$ with a twin of $w$ to obtain a complete $(m-2)$-partite graph which is a vertex in $\mathcal{G}$.
Otherwise, $x \in V(J)$ is adjacent in $J$ to all of $V_3 \cup \ldots \cup V_m$, and so $J'$ is a complete $(m-1)$-partite graph which is a vertex in $\mathcal{G}$.
In both cases, we obtain a contradiction to the choice of $m$.
This proves the claim.
\ecpf

\medskip
\noindent
Therefore $x$ has a neighbour in each of $V_1,\ldots,V_m$.
Then, for each $i \in [m]$, $V_i$ has partition $A_i,B_i$, where $xa \in E(H^-)$ for all $a \in A_i$, and $xb \notin E(H^-)$ for all $b \in B_i$; and $A_i \neq \emptyset$.
Observe that every $A_i$ is a set of twins, and $B_i$ is a set of twins.
For each $i \in [m]$, replace $B_i$ by a set of $|B_i|$ twins of $a \in A_i$.
Thus obtain a graph $H' \cong K(\{ x \}, V_1,\ldots,V_m)$ which is a vertex of $\mathcal{G}$, as required.
\epf

\begin{lemma}\label{defective}
Let $s \in \mathbb{N}$ and $\bm{k} \in \mathbb{N}^s$.
Then, for all $\eps,\gamma > 0$, there exist $\eta > 0$ and $n_0 \in \mathbb{N}$ such that the following holds for all graphs $G$ on $n \geq n_0$ vertices with $\log F(G;\bm{k}) \geq (Q(\bm{k})-\eta)\binom{n}{2}$.
For every $x \in V(G)$ such that $G$ contains at least $\gamma n$ twins of $x$, we have that
\begin{equation}\label{def2}\log F(G; \bm{k}) \geq (Q(\bm{k})-\eps)n + \log F(G-x; \bm{k}).
\end{equation}
\end{lemma}

\bpf
As before, we omit $\bm{k}$ from our notation where it is clear from the context, so e.g.~$F(H) := F(H;\bm{k})$.
%For all graphs $H$ and integers $n$, write $F(H) := F(H;\bm{k})$ and $F(n) := F(n;\bm{k})$.
Choose constants $\delta,\eta$ such that $0 < \eta \ll \delta \ll \gamma,\eps$.
By~(\ref{FQeq}), we can choose $n_0$ to be such that whenever $N \geq n_0/2$, we have that
\begin{equation}\label{etabound}
(Q - \eta)\binom{N}{2} \leq \log F(N) \leq (Q + \eta)\binom{N}{2}. 
\end{equation}
Without loss of generality, we may suppose that $1/n_0 \ll \eta$.
Let $G$ be a graph on $n \geq n_0$ vertices.
Suppose that there is some $x \in V(G)$ which does not satisfy (\ref{def2}).
Let $T \subseteq V(G)$ be the set of twins of $x$ (including $x$).
Let $\chi(G-T)$ be the set of $\bm{k}$-valid colourings of $G-T$ and let $c(\chi,x)$ be the number of extensions of $\chi$ to $G - (T \setminus \{ x \})$.
For each $\chi \in \chi(G-T)$, since every pair of vertices of $X$ are twins, we have that $c(\chi,x) \equiv c_\chi$ for all $x \in T$.
Therefore, if we list the vertices of $T$ as $x_1,\ldots,x_{|T|}$, and let $t \in [\,|T|\,]$, we have that
\begin{equation}\label{power}
F(G-x_1-\ldots-x_t) = \sum_{\chi \in \chi(G-T)}c_\chi^{|T|-t}.
\end{equation}
%Suppose that~(\ref{def2}) does not hold.
Choose $t$ with $1 \leq t \leq |T|$ to be maximal such that
\begin{equation}\label{contra}
\log F(G-x_1-\ldots-x_{i-1}) < (Q-\eps)n + \log F(G-x_1 - \ldots - x_i)
\end{equation}
for all $i \in [t]$.
(Since (\ref{def2}) does not hold, $t:=1$ satisfies the inequality).
Suppose $t < |T|$.
By the Cauchy-Schwarz inequality,
\begin{eqnarray*}
F(G-x_1-\ldots - x_t)^2 &=& \left(\sum_{\chi \in \chi(G-T)}c_\chi^{|T|-t} \right)^2 \leq \left(\sum_{\chi \in \chi(G-T)} c_\chi^{|T|-t-1} \right) \left( \sum_{\chi \in \chi(G-T)} c_\chi^{|T|-t+1} \right)\\
&\stackrel{(\ref{power})}{=}& F(G-x_1-\ldots-x_{t+1}) F(G-x_1-\ldots-x_{t-1})\\
&\stackrel{(\ref{contra})}{<}& 2^{(Q-\eps)n} F(G-x_1-\ldots-x_{t+1}) F(G-x_1-\ldots-x_t).
\end{eqnarray*}
So
$F(G-x_1-\ldots-x_t) < 2^{(Q-\eps)n} \cdot F(G-x_1-\ldots-x_{t+1})$.
Thus (\ref{contra}) holds for $t+1$ which is a contradiction to the maximality of $t$. 
So $t=|T|$.
Therefore, inductively, for all non-empty $T' \subseteq T$,
$$
\log F(G) \leq (Q-\eps)|T'|n + \log F(G-T').
$$

Choose $T' \subseteq T$ with $|T'| = \lceil \delta n\rceil$.
Then
\begin{align*}
\log F(G) &\leq (Q-\eps)|T'|n + \log F(G-T') \stackrel{(\ref{etabound})}{\leq} Q|T'|n -\eps \delta n^2 + (Q+\eta)\binom{n-|T'|}{2}\\
&\leq Q\binom{n}{2}  + (-\eps \delta + \eta +Q\delta^2)n^2 < Q\binom{n}{2} + (\eta + \delta^2 \log s - \eps \delta )n^2\\
&< (Q-\eta)\binom{n}{2},
\end{align*}
a contradiction to~(\ref{etabound}).
\epf

We will now prove Theorem~\ref{exact1}. The idea is that if there is an extremal graph $G$ not satisfying the conclusion of the theorem, there must be a complete partite extremal graph $G'$ with an induced subgraph satisfying the conclusion of the theorem which is almost the whole of $G'$. We can either set $G'=G$ or symmetrise to obtain $G'$ (Lemma~\ref{small}), and the required structure of $G'$ follows from Theorem~\ref{stabilitycomp}.
In a typical colouring, of which there are many by Theorem~\ref{stabilitycomp}, every exceptional vertex (not in the good induced subgraph) has deficient contribution to $F(G')$, which is a consequence of the strong extension property.
We can form a new graph by replacing each exceptional vertex by a twin of a good vertex to obtain a graph with more valid colourings than $G'$, which is the required contradiction.

\bpf[Proof of Theorem~\ref{exact1}]
Again we omit $\bm{k}$ from our notation where possible.
Let $\eps > 0$.
Let $\mu > 0$ be the constant obtained from Lemma~\ref{lemma2}\ref{solutions} and let $\eps_0$ be the constant obtained from Lemma~\ref{lemma2}\ref{clone}.
Without loss of generality, we may assume that $0 < \eps_0 \ll \mu,\eps,1/R,1/s$.
Choose $\delta_4,\delta_3,\delta_2,\delta_1>0$ in this order such that 
$\delta_3$ is at most the output of Corollary~\ref{stabilitycomp2} with input $\delta_4$,
and similarly $\delta_2$ from $\delta_3$, and 
$\delta_1$ from $\delta_2^2$.
%$\delta_3 \ll \eps_0$ and let $\delta_2,n_3$ be the outputs of Corollary~\ref{stabilitycomp2} applied with parameter $\delta_3$.
%Let $\delta_1,n_2$ be the outputs of Corollary~\ref{stabilitycomp2} applied with parameter $\delta_2^2$.
Further, let $n_0$ be at least the integer output of all of these applications.
By increasing $n_0$ and decreasing $\delta_1,\ldots,\delta_4$, we may assume that $0<1/n_0 \ll \delta_1 \ll \ldots \ll \delta_4 \ll \eps_0$ and for all $n \geq n_0$, by~(\ref{FQeq}),
\begin{equation}\label{logFn}
\left(Q - \delta_1 \right)\binom{n}{2} \leq \log F(n) \leq \left(Q + \delta_1 \right)\binom{n}{2}
\end{equation}
and that Lemma~\ref{defective} applied with $\delta_3,\mu/2$ playing the roles of $\eps,\gamma$ respectively has output $n_1 \leq n_0$ and $\eta > 2\delta_1$.
%\begin{equation}\label{logFG}
%\log F(H;\bm{k}) \geq (Q(\bm{k})-\delta_2)n + \log F(H-x;\bm{k}) \quad\text{for all }\bm{k}\text{-extremal }H \text{ and }x \in V(G) \text{ with at least }\mu n/2 \text{ twins in }H.
%\end{equation}
Altogether
we have the hierarchy
\begin{equation}\label{hierexact}
0 < 1/n_0 \ll \delta_1 \ll \delta_2 \ll \delta_3 \ll \delta_4 \ll \eps_0 \ll \eps,\mu,1/R,1/s.
\end{equation}

We first prove part~\ref{exact1i}.
Towards a contradiction, let $G$ be a $\bm{k}$-extremal graph on $n \geq 2n_0$ vertices such that either
\begin{itemize}
\item[(a)] $G$ is not complete multipartite; or
\item[(b)] there is some $m \in \mathbb{N}$ and $\ba \in \Delta^m$ such that $G \cong K_{\alpha_1n,\ldots,\alpha_mn}$ but there is no $(m,\ba^*) \in \wt(\bm{k})$ with $\|\ba-\ba^*\|_1 < \delta_3 < \eps$.
\end{itemize}
If (a) holds, apply Lemma~\ref{small} to obtain a $\bm{k}$-extremal graph $G'$ on $n$ vertices which is complete $m$-partite with vertex partition $X_1,\ldots,X_{m}$, where $|X_{1}| \geq \ldots \geq |X_m| = 1$.
Observe that $m < R$.
If (b) holds, set $G' := G$.
In both cases, apply Corollary~\ref{stabilitycomp2} with parameter $\delta_2^2$ to $G'$.
Part~\ref{sci} implies that there exists $(r^*,\ba^*) \in \wt(\bm{k})$ with $r^* \leq m$ such that, defining $Y_i := X_i$ for $i \in [r^*]$ and $Y_0$ to be the union of the remaining parts of $G$, we have
\begin{equation}\label{Xbds}
\sum_{i \in [r^*]}|\,|Y_i|-\alpha_i^*n\,| < \delta_2^2n < \delta_2 n.
\end{equation}
%Observe that, in both cases (a) and (b), we have $0 < |Y_0| < \delta_2 n$.
%Suppose, without loss of generality, that $|Y_1| \geq \ldots \geq |Y_m|$.

Together with Proposition~\ref{continuitylagrange}\ref{continuitylagrangeii} we have that $|Y_i| \geq \mu n/2$ for all $i \in [r^*]$.
So if either (a)  or (b) holds, we have $m > r^*$.
%For each $i \in [r^*]$, let $\ell_i$ be the set of $k \in [m]$ such that $X_k \subseteq Y_i$.
Corollary~\ref{stabilitycomp2} implies that each of $Y_1,\ldots,Y_{r^*}$ is a part of $G'$ (but $Y_0$ may contain several parts).
We have that $Y_0,\ldots,Y_{r^*}$ is a partition of $V(G)$, and $0 < |Y_0| < R\delta_2^2 n < \delta_2 n$.
From now on, we make no distinction between cases (a) and (b), and only use this fact about the size of $Y_0$.
Let
$$
H := G'[Y_1 \cup \ldots \cup Y_{r^*}] \cong K[Y_1,\ldots,Y_{r^*}]
$$
be the \emph{core} of $G'$, and let $N := |V(H)|$. So
\begin{equation}\label{n-N}
0 < n - N = |Y_0| < R\delta_2^2 n < \delta_2 n.
\end{equation}
%Apply Corollary~\ref{stabilitycomp2} with parameter $3\sqrt{\delta_2}\log s$ to $H$ to see that 
%there are at least $(1-2^{-3\sqrt{\delta_2}\log s n^2})\cdot F(H)$ colourings which are

%By definition, every $\delta$-good colouring of $G'$ restricts to a $\delta$-good colouring of $H$,
%and each such colouring of $H$ is the restriction of at most $s^{(n-N)n}$ colourings of $G'$.
%Thus
%\begin{eqnarray}
%|\mathcal{G}_{\delta_2}(H;Y_1,\ldots,Y_{r^*})| \geq s^{-(n-N)n}\cdot |\mathcal{G}_{\delta_2}(H;Y_1,\ldots,Y_{r^*})| \geq s^{-(n-N)n}| \geq s^{-(n-N)n}(1-2^{-\delta_1 n^2}) 2^{(Q(\bm{k})-\delta_1)\binom{n}{2}}.
%\end{eqnarray}

Every $\bm{k}$-valid colouring of $G'$ can be obtained by extending a $\bm{k}$-valid colouring of $H$, so $F(H) \geq F(G') \cdot s^{-|Y_0|n}$.
Therefore
\begin{eqnarray}
\nonumber \log F(H) &\stackrel{(\ref{Xbds})}{\geq}& \log F(G') - R\delta_2^2 \log s \cdot n^2 \stackrel{(\ref{logFn})}{\geq}  \left(Q - \delta_1 \right) \binom{n}{2} - R\delta_2^2 \log s \cdot n^2\\
\label{FH} &\stackrel{(\ref{n-N})}{\geq}& (Q-\delta_2)\binom{N}{2}.
\end{eqnarray}
Apply Corollary~\ref{stabilitycomp2} to $H$ to see that $|\mathcal{G}(H)| \geq (1-2^{-\delta_3 N^2})\cdot F(H)$, where
(recalling Definition~\ref{def:colourings}) 
$$\mathcal{G}(H) := \bigcup_{\substack{\phi^* \in \pat(r^*,\ba):\\ (r^*,\ba) \in \wt(\bm{k}), \|\ba-\ba^*\|_1 \leq \delta_4}}\mathcal{G}_{\delta_4}(H;\phi^*,Y_1,\ldots,Y_{r^*}).$$
Define $\bb \in \Delta^{r^*}$ by setting $\beta_i := |Y_i|/N$ for all $i \in [r^*]$.
So~(\ref{Xbds}) implies that for all $\ba$ with $\|\ba-\ba^*\|_1 \leq \delta_4$,
\begin{equation}\label{normdiff}
\|\bb - \ba\|_1 < 2\delta_2^2 + \delta_4 \leq 2\delta_4.
\end{equation}
%Corollary~\ref{stabilitycomp2} implies that $|\mathcal{G}_{\delta_2}(G';Y_1,\ldots,Y_{r^*})| \geq (1-2^{-\delta_1 n^2}) \cdot F(G)$, so
%$$
%|\mathcal{G}_{\delta_2}(H;Y_1,\ldots,Y_{r^*})| \geq (1-2^{-\delta_1 n^2}) \cdot s^{-\sqrt{\delta_2}n} \cdot F(G) \geq (1-2^{-\delta_1 n^2/2}) \cdot F(H).
%$$
Let $\mathcal{B}(H)$ be the set of $\bm{k}$-valid colourings of $H$ not in $\mathcal{G}(H)$.
For each $\chi \in \mathcal{G}(H)$ and $v \in Y_0$, do the following.
Let $\phi^* \in \Phi_2(r^*;\bm{k})$ be the pattern of $\chi$ (that is, there is $\ba$, which depends on $\chi$, with $(r^*,\ba) \in \wt(\bm{k})$ and $\|\ba-\ba^*\|_1 \leq \delta_4$ and $\phi^* \in \pat(r^*,\ba)$)
such that $\chi \in \mathcal{G}_{\delta_4}(H;\phi^*,Y_1,\ldots,Y_{r^*})$.
For each valid extension $\overline{\chi}$ of $\chi$ to $G'[V(H) \cup \{ v \}]$ and every
$i\in [s]$, define $\phi = \phi(\overline{\chi},v) : \binom{[r^*+1]}{2} \rightarrow 2^{[s]}$ by setting
$$
\phi(ij) := \begin{cases} 
\phi^*(ij) &\mbox{if } ij \in \binom{[r^*]}{2} \\
\{ c \in [s] : |\overline{\chi}^{-1}(c)[v,Y_i]| \geq \sqrt{\delta_4} |Y_i| \} &\mbox{if } i \in [r^*], j = r^*+1. \\
\end{cases}
$$
Fix the $\phi$ that appears for the largest number of
extensions $\overline{\chi}$ over all $\chi \in \mathcal{G}(H)$ and all $v \in Y_0$.

\begin{claim} $\phi \in \Phi_0(r^*+1,\bm{k})$.\end{claim}

\bpf
This is almost identical to Claim~\ref{claimphi1} so we omit the proof.
%Suppose the statement does not hold.
%Then there is some $c \in [s]$ such that $\phi^{-1}(c)$ contains a copy of $K_{k_c}$.
%Since $\phi$ is an extension of $\phi^* \in \Phi(r^*,\bm{k})$, the vertex $r^*+1$ lies in this copy.
%Let $z_1,\ldots,z_{k_c-1}$ be the other vertices.
%For each $i \in [k_c-1]$, let $Z_i := \{ x \in Y_{z_i} : \chi'(vx) = c \}$.
%Then $|Z_i|/|Y_i| \geq \sqrt{\delta_3}$ by the definition of $\phi$.
%Let $ij \in \binom{[k_c-1]}{2}$.
%Then
%$$
%d(\chi^{-1}(c)[Z_i,Z_j]) \geq d(\chi^{-1}(c)[Y_1,Y_j]) - \delta_3 \geq |\phi(ij)|^{-1} - \delta_3 \geq 1/2s. 
%$$
%Now Proposition~\ref{badrefine} implies that $\chi^{-1}(c)[Z_i,Z_j]$ is $\sqrt{\delta_3}$-regular.
%Therefore $\chi^{-1}(c)[Z_i,Z_j]$ is $(\sqrt{\delta_3},\geq 1/2s)$-regular.
%Now Lemma~\ref{embed} implies that $\chi^{-1}(c)[Z_1,\ldots,Z_{k_c-1}]$ contains a copy of $K_{k_c-1}$.
%Together with $x$, this gives a $c$-coloured copy of $K_{k_c}$ in $\chi$, contradicting the fact that $\chi$ is $\bm{k}$-valid. 
\epf

Proposition~\ref{continuitylagrange}\ref{continuitylagrangei} implies that
\begin{equation}\label{qb}
q(\phi^*,\bb) \geq q(\phi^*,\ba) - 2\log s \|\bb-\ba\|_1 \stackrel{(\ref{normdiff})}{\geq} Q - 4(\log s)\delta_4.
\end{equation}
Since every $v \in Y_0$ is incident to the whole of $Y_1 \cup \ldots \cup Y_{r^*}$ in $G'$ and $s\sqrt{\delta_4}<1$, we have that $\phi(\{ i,r^*+1 \}) \neq \emptyset$ for all $i \in [r^*]$.
Thus $r^*+1$ is not a strong clone of any vertex in $[r^*]$ under $\phi$.
Therefore, applying Lemma~\ref{lemma2}\ref{clone} for the second inequality,
\begin{eqnarray}
\label{betasum}\ext(\phi,\bb) &\leq&  \ext(\phi,\ba) + \log s \cdot \|\bb-\ba\|_1 \stackrel{(\ref{normdiff})}{<} Q - \eps_0 + 4(\log s)\delta_4 < Q - \frac{\eps_0}{2}.
\end{eqnarray}
Similarly to the derivation of~(\ref{total}), we can now bound the number of extensions of $\chi$.
Indeed, by our choice of $\phi$, for any $\chi \in \mathcal{G}(H)$ and any $v \in Y_0$, the number of ways to extend $\chi$ to a valid colouring $\overline{\chi}$ of $G'[V(H) \cup \{ v \}]$ is at most
\begin{align*}
&\phantom{\leq} 2^{sr^*} \cdot \prod_{i \in [r^*]}|\phi(\{ i,r^*+1 \})|^{\beta_i N} \cdot \binom{N}{\leq s\sqrt{\delta_4} N} \cdot s^{sr^*\sqrt{\delta_4} n} \stackrel{(\ref{betasum})}{\leq}
%2^{\delta_4 N} \cdot 2^{(Q(\bm{k}-\eps_0/2)N} \leq 
2^{(Q-\eps_0/3)N}.
\end{align*}
Here, the first term is the number of possible patterns $\phi$; the second term is the maximum number of extensions given the $\phi$ that appears most often; and the third and fourth terms are an upper bound on the number of ways to choose and colour those edges with uncommon colours.
Now we can give an upper bound for the number of valid colourings of $G'$ as follows:
Any valid colouring of $G'$ is either an extension of $\chi \in \mathcal{G}(H)$ (where we must additionally colour the edges between the $n-N$ vertices of $Y_0$ and $H$, and the edges induced by $Y_0$); or an extension of $\chi \in \mathcal{B}(H)$.
We have 
\begin{eqnarray}
\label{G'} F(G') &\leq& \sum_{\chi \in \mathcal{G}(H)} \left( s^{\binom{n-N}{2}} \cdot 2^{(Q-\eps_0/3)N(n-N)} \right) + \sum_{\chi \in \mathcal{B}(H)}s^{(n-N)n}\\
\nonumber &\leq& s^{\binom{n-N}{2}} \cdot 2^{(Q-\eps_0/3)N(n-N)} \cdot F(H) + 2^{-\delta_3 N^2} F(H) \cdot s^{(n-N)n}\\
\nonumber &\stackrel{(\ref{n-N})}{\leq}& 2^{(Q - \eps_0/3 + \delta_2 \log s) N(n-N)} \cdot F(H) + 2^{-\delta_3 N^2/2} \cdot F(H)\\
\nonumber &\leq& 2^{(Q-\eps_0/4)N(n-N)} \cdot F(H).
\end{eqnarray}

We will now form a new graph on $n$ vertices which shares the same core $H$, but $Y_0$ is replaced by $n-N$ clones of some vertex in another part.

By (\ref{qb}), for any $\phi^* \in \pat(r^*,\ba)$,
$$
Q - 4(\log s)\delta_4 \leq q(\phi^*,\bb) = \sum_{i \in [r^*]}\beta_i q_i(\phi^*,\bb).
$$
So, by averaging, there exists $i^* \in [r^*]$ such that
\begin{equation}\label{average}
q_{i^*}(\phi^*,\bb) \geq Q - 4(\log s)\delta_4.
\end{equation}
Let $H' := K[W_1,\ldots,W_{r^*}]$, where $W_j = Y_j$ for $j \in [r^*]\setminus \{ i^* \}$, and $W_{i^*} := Y_{i^*} \cup Y_0$.
%For each $\chi \in \mathcal{G}$ and $v \in Y_{i^*}\setminus X_{i^*}$, the number of ways of extending $\chi$ to $H+v$ is at least
%$$
%\prod_{j \in [r^*]\setminus \{ i^* \}}|\phi^*(i^*j)|^{\beta_j N} \geq 2^{(Q(\bm{k})-4\log s \delta)N}.
%$$
Every colouring that follows $\phi^*$ is valid. Thus
\begin{eqnarray*}
\log F(H') &\geq& \log\left( \prod_{ij \in \binom{[r^*]}{2}}|\phi^*(ij)|^{|Y_i||Y_j|} \cdot \prod_{k \in [r^*]\setminus\{i^*\}}|\phi^*(i^*k)|^{|Y_0||Y_k|} \right)\\
&=& N^2 \cdot q(\phi^*,\bb)/2 + (n-N)N \cdot q_{i^*}(\phi^*,\bb)/2 \geq nN(Q-4(\log s)\delta_4)/2\\
&\geq& (Q-\sqrt{\delta_4}/2)\binom{n}{2} \stackrel{(\ref{logFn})}{\geq} F(H) + (Q-\sqrt{\delta_4})N(n-N),
\end{eqnarray*}
which, together with $n-N>0$ and~(\ref{G'}) implies that $\log F(H') > \log F(G')$, a contradiction to the $\bm{k}$-extremality of $G'$, and hence $G$.
%\begin{eqnarray*}
%F(H') &\geq& \sum_{\chi \in \mathcal{G}(H)} \left( \prod_{j \in [r^*]\setminus \{ i^* \}}|\phi^*(i^*j)|^{\beta_jn} \right)^{n-N} \stackrel{(\ref{average})}{\geq} \sum_{\chi \in \mathcal{G}(H)}2^{(Q(\bm{k})-2\delta_2\log s)N(n-N)}\\
%&\geq& (1-2^{-\eps N^2}) \cdot F(H) \cdot 2^{(Q(\bm{k})-2\delta_2\log s)N(n-N)} \geq 2^{(Q(\bm{k})-3\delta_2\log s)N(n-N)} \cdot F(H)\\
%&\stackrel{(\ref{hierexact})}{>}& 2^{(Q(\bm{k})-\eps_0/4)N(n-N)} \cdot F(H) \cdot 2^{\eps_0n/5} \stackrel{(\ref{G'})}{\geq} 2^{\eps_0 n/5} \cdot F(G'),
%\end{eqnarray*}
%a contradiction to the $\bm{k}$-extremality of $G'$, and hence of $G$.
This completes the proof of part (i).

\medskip
\noindent
We have proved that $G \cong K[Y_1,\ldots,Y_{r^*}]$ such that $\|\bb-\ba^*\|_1 < \delta_2$.
Now~Lemma~\ref{defective} and our choice of parameters implies that, for all $x \in V(G)$,
$$
\log F(G) \geq (Q-\delta_3)n + \log F(G-x).
$$
Note that $F(G) \leq s^n \cdot F(G-x)$, so, by~(\ref{logFn}), $\log F(G-x) \geq (Q-2\delta_1)\binom{n-1}{2}$.
The hypotheses of Lemma~\ref{defective} still hold for $G-x$, so, for all $y \in V(G)\setminus\{x\}$,
$$
\log F(G) \geq (Q-\delta_3)n + \log F(G-x) \geq (Q-\delta_3)(n+n-1) + \log F(G-x-y).
$$

By Corollary~\ref{stabilitycomp2}, there is a set $\mathcal{G}$ of at least $(1-2^{-\delta_3 n})\cdot F(n)$ $\bm{k}$-valid colourings $\chi : E(G) \rightarrow [s]$ for which there exists $(r^*,\ba) \in \wt(\bm{k})$ with $\|\ba^*-\ba\|_1 \leq \delta_4$ and $\phi \in \pat(r^*,\ba)$ such that $\chi$ is $(0,\delta_4)$-perfect with respect to $(\phi,Y_1,\ldots,Y_{r^*})$.
\epf

\bpf[Proof of Corollary~\ref{exact2}]
Suppose that $\bm{k}$ is soluble, let $c$ be the constant implied by solubility and let $0 < \eps \ll c$.
Let $\delta,n_0$ be the output of Theorem~\ref{exact1} applied with parameter $\eps$.
We may assume that $1/n_0 \ll \delta \ll \eps$.
Let $n \geq n_0$ and let $G$ be an $n$-vertex $\bm{k}$-extremal graph.
Theorem~\ref{exact1} implies that $G$ is a complete $r^*$-partite graph for some $r^* \in \mathbb{N}$, where,
writing $\bb \in \Delta^{r^*}$ for the vector of the part ratios of $G$, we have $\|\bb-\ba^*\|_1 \leq \eps$ for some $(r^*,\ba^*) \in \wt(\bm{k})$.
Let $\bm{d}^*$ be the supersolution of $n$ with $\perf(\bm{d}^*) = \perf_{r,\ba}(\bm{d}^*)$ for some $(r,\ba) \in \wt(\bm{k})$.
Suppose for a contradiction that $\bb n \neq \bm{d}^*$.

The number of perfect colourings of $G$ is at least $(1-2^{-\delta n})F(G;\bm{k})$.
Then
$$
F(K_{\bm{d}^*}) \geq \perf(\bm{d}^*) \geq (1+c) \cdot \perf_{r^*,\ba^*}(\bb n) \geq (1+c)\cdot (1-2^{-\delta n})F(G) > F(G),
$$
contradicting the $\bm{k}$-extremality of $G$.
\epf

\section{Forbidden triangles in seven colours}\label{trianglesec}

In this section, we solve Problem $Q^*$ in the case $\bm{k}=(3;7)$.

\begin{theorem}\label{37thm}
Let $\bm{k}:=(3;7)$.
Then every $(r,\phi,\ba)\in\opt^*(\bm{k})$ has $r=8$, $\ba$ uniform, $|\phi(ij)|=4$ for all $ij \in \binom{[8]}{2}$ and $\phi^{-1}(c) \cong K_{4,4}$ for all $c \in [7]$. 
\end{theorem}

We also solve it in the case $\bm{k}=(3;6)$. Recall that $F(n;\bm{k})$ was already determined exactly in this case in~\cite{skokan}.

\begin{theorem}\label{36thm}
Let $\bm{k}:=(3;6)$.
Then every $(r,\phi,\ba)\in\opt^*(\bm{k})$ has $r=8$, $\ba$ uniform, $|\phi(ij)|\in \{3,4\}$ for all $ij \in \binom{[8]}{2}$, $\{ij: |\phi(ij)|=3\} \cong K_{4,4}$ and $\phi^{-1}(c) \cong K_{4,4}$ for all $c \in [6]$. 
\end{theorem}

\subsection{The union of dense triangle-free graphs}

The key new idea is the following lemma, a $2$-coloured version of Mantel's theorem, that allows us to add a new constraint to the linear relaxation of the optimisation problem.
The constraint ensures that the union of any two colour graphs $R := \phi^{-1}(c), B := \phi^{-1}(c')$ has density at most $\frac{3}{4}$, whenever the individual graphs have large density.
This is attained by two complete balanced bipartite graphs whose overlap is minimal.
The trivial bound for the density of $R \cup B$ is $\frac{4}{5}$. Indeed, $R \cup B$ is $K_6$-free, otherwise it would contain a monochromatic triangle, so this claim follows from Tur\'an's theorem. However, if $R \cup B$ has density $\frac{4}{5}$, then each of $R,B$ has density $\frac{2}{5}$, coming from the unique red-blue colouring of $K_5$ without monochromatic triangles.
The lemma states that if the sum of densities of $R,B$ is larger, closer to the maximum of $\frac{1}{2}+\frac{1}{2}$, then $R \cup B$ has a density significantly smaller than $\frac{4}{5}$.

\begin{lemma}\label{RB}
Let $(a,b)$ be $(\frac{19}{25},\frac{89}{100})$ or $(\frac{3}{4},\frac{19}{20})$.
Let $R,B$ be two triangle-free graphs on the same vertex set of size $n$ with $|R|+|B| \geq b n^2/2$.
Then $|R \cup B| \leq a n^2/2 + o(n^2)$.
\end{lemma}

\bpf
Let $n$ be a sufficiently large integer. Let $R$ and $B$ be two triangle-free graphs with vertex set $[n]$. We will sometimes denote their edge sets as $R$ and $B$ respectively too.
Write $R \cup B$ for the set of edges that lie in at least one of $R,B$ (a simple graph), and $R+B$ for the multiset union of $R,B$ (a multigraph).
%So $d(R+B)=d(R)+d(B)$ etc.
Suppose that $\ed(R+B) > b$ and $\ed(R \cup B)>a$, where $\ed(E):=|E|/\binom{n}{2}$ is the edge density of a set $E$ of edges.
Let $Y \subseteq [n]$ be a maximal set with the property that $Y$ has a partition $Y_1\cup\ldots\cup Y_t$ into sets of size $5$ where $(R \cup B)[Y_i]\cong K_5$ for all $i \in [t]$. Since there is a unique $2$-edge-colouring of $K_5$ which avoids monochromatic triangles, we can label the vertices in each $Y_i$ as $y^i_{1},\ldots,y^i_{5}$ where $y^i_{j}y^i_{j+1} \in B$ for all $j \in [5]$, where $y^i_{6}:=y^i_{1}$, and every other pair is in $R$ (so there are no double edges).
Let $X \subseteq \overline{Y} := [n]\setminus Y$ be a maximal set with the property that $X$ has a partition $X_1 \cup \ldots \cup X_s$ into pairs where $(R \cap B)[X_i]\cong K_2$ for all $i \in [s]$.
Write $|Y|=yn$, $|X|=xn$ and $Z:=[n]\setminus (X \cup Y)$.
%Suppose for a contradiction that $y>\eps$.

We claim that
\begin{align}\label{R+B}
\nonumber &0 \leq \ed(R+B) - b \leq 2q(x,y) + o(1)\quad\text{where}\\
&q(x,y) := \frac{x^2}{2}+\frac{3}{4}\cdot\frac{(1-x-y)^2}{2}+x(1-x-y)+\frac{4}{5}\cdot\frac{y^2}{2}+\frac{4}{5}\cdot y(1-y) - \frac{b}{2}.
\end{align}
This is a consequence of the following observations on the densities of $R$ and $B$ in various subsets of the vertex set.
For a graph $G$ with disjoint $U,U' \subseteq V(G)$ we write $\ed_G(U) := |E(G[U])|/\binom{|U|}{2}$ and $\ed_G(U,U') := |E(G[U,U'])|/(|U||U'|)$.
\begin{enumerate}[label=(\roman*),ref=(\roman*)]
\item\label{dX} $\ed_R(X) \leq \frac{1}{2}+o(1)$ by Tur\'an's theorem (or rather Mantel's theorem), since $R$ is triangle-free. Similarly for $\ed_B(X)$.
\item\label{dZ} $\ed_{R+B}(Z) \leq \frac{3}{4}+o(1)$ by Tur\'an's theorem, since $R,B$ are disjoint on $Z$ and so $R+B=R \cup B$ is $K_5$-free, so $\ed_{R+B}(Z)=\ed_{R \cup B}(Z) \leq \frac{3}{4}+o(1)$.
\item\label{dXZ} $\ed_{R+B}(X,Z) \leq 1$.
If not, by averaging, there is $j \in [s]$ with $e_{R+B}(X_j,Z)>2|Z|$, so without loss of generality there is $v \in Z$ such that $v$ sends a red edge to both vertices in $X_j$, a contradiction.
\end{enumerate}
Thus~(\ref{R+B}) holds for $y=o(1)$.
So for the rest of the derivation we may suppose $y=\Omega(1)$.
\begin{enumerate}[label=(\roman*),ref=(\roman*)]
\setcounter{enumi}{3}
\item\label{dY} $\ed_{R+B}(Y) \leq \frac{4}{5}+o(1)$. If not, by averaging, there are distinct $i,j \in [t]$ such that
$$
e_{R+B}(Y_i,Y_j) \geq \frac{(\frac{4}{5}+\Omega(1))\binom{yn}{2} - 10\cdot yn/5}{\binom{yn/5}{2}} >  20
$$
(if $\ed_{R+B}(Y) > \frac{4}{5}+c$ then this holds already for $yn>\frac{1}{c}-1$).
Then, without loss of generality, $d_{R+B}(y^j_{1},Y_i) \geq 5$.
There is at least one double edge, say $y^j_{1}y^i_{1}$, for otherwise $(R\cup B)[Y_i \cup \{y^j_{1}\}]\cong K_6$. 
Then $y^j_{1}$ sends no blue edges to $\{y^i_{2},y^i_{5}\}$, and at most one red edge.
Similarly, $y^j_{1}$ sends no red edges to $\{y^i_{3},y^i_{4}\}$, and at most one blue edge.
 Thus in fact $d_{R+B}(y^j_{1},Y_i) \leq 4$, a contradiction.
\item\label{dYYc} For all $W \subseteq X \cup Z$, $\ed_R(Y,W) \leq \frac{2}{5}$.
If not, by averaging, there is $i \in [t]$ and $v\in W$ such that $d_R(v,Y_i) \geq 3$ which yields a red triangle containing $v$ and two vertices in $Y_i$.
Similarly for $\ed_B(Y,W)$.
\end{enumerate}
This proves that~(\ref{R+B}) holds, using in order~\ref{dX}--\ref{dYYc} for each term, bounding the density of $R+B$ in, respectively,
$X,Z,(X,Z),Y,(Y,X \cup Z)$.

We obtain a similar polynomial upper bound for the density of $R\cup B$.
We claim that
\begin{align}
\nonumber &0 \leq \ed(R \cup B) - a \leq 2p(x,y) + o(1)\quad\text{where}
\end{align}
\begin{align}\label{RcupB}
p(x,y) &:= \frac{4}{5}\cdot\frac{y^2}{2}+ \frac{4}{5}(1-x-y)y + \frac{3}{4}\cdot\frac{(1-y)^2}{2} + \frac{7}{10} xy - \frac{a}{2}\\
\nonumber &= -\frac{y}{40}(4x+y-2) - \frac{1}{2}\left(a-\frac{3}{4}\right).
\end{align}
This follows from some more observations.
\begin{enumerate}[label=(\roman*),ref=(\roman*)]
\setcounter{enumi}{5}
\item\label{dYc} $\ed_{R \cup B}(\overline{Y}) \leq \frac{3}{4}+o(1)$ by Tur\'an's theorem, since $(R \cup B)[\overline{Y}]$ is $K_5$-free.
\item\label{dYX} $\ed_{R \cup B}(Y,X) \leq \frac{7}{10}$.
If not, then by averaging there is $i \in [t]$ and $j \in [s]$ such that $e(R \cup B)[Y_i,X_j] > \frac{7}{10}\cdot 5 \cdot 2 = 7$.
%We check that this is impossible by a simple case analysis.
For ease of notation, write $1,\dots,5$ for the vertices of $Y_i$, with blue edges forming the cycle $12,23,34,45,51$ 
and red edges forming the cycle $13,35,52,24,41$, and write $X_j=\{x,y\}$, where $xy$ is both red and blue.
Each of $x,y$ has at most two neighbours in $[5]$ of any one colour, and therefore exactly two or we are done.
For blue these neighbours are a subset of one of
$\{1,3\},\{2,4\},\{3,5\},\{4,1\},\{5,2\}$, and for red a subset of
$\{1,5\},\{3,2\},\{5,4\},\{2,1\},\{4,3\}$.
We have that $N_B(x),N_B(y)$ are disjoint, as are $N_R(x),N_R(y)$.
Also, each of $N_B(x),N_R(x)$ are disjoint, as are $N_B(y),N_R(y)$, since otherwise there are at most $7$ $R \cup B$-edges between $\{x,y\}$ and $\{1,\ldots,5\}$, a contradiction.
Without loss of generality, suppose $N_B(x)=\{1,3\}$.
Then $N_R(x) = \{5,4\}$ and $N_B(y)\in\{\{2,4\},\{5,2\}\}$.
For either possibility, there  is no choice of $N_R(y)$ which satisfies the disjointness conditions. 
\end{enumerate}
Again~(\ref{RcupB}) holds using~\ref{dY}--\ref{dYX} in order for each term, bounding the density of $R \cup B$ in, respectively, $Y,(Z,Y),X \cup Z, (X,Y)$.
This proves the claim.

\begin{claim}
The regions $P:= \{(x,y) \in [0,1]^2: p(x,y) \geq 0\}$ and $Q := \{(x,y) \in [0,1]^2: q(x,y) \geq 0\}$ intersect only when $y=0$.
\end{claim}

\bcpf
Write
\begin{align*}
p_x(y) &:=40p(x,y) = -y^2 - 2y(2x-1) - 5(4a - 3)\quad\text{and}\\
q_x(y) &:=40q(x,y)= -y^2 - 2y(5x-1) +5(-x^2+2x+3-4b)
\end{align*}
which are both quadratic functions of $y$ with negative $y^2$ coefficient.
Note that their discriminants are
\begin{align*}
{\rm disc}(p_x) &= 4(2x-1)^2-20(4a-3)\quad\text{and}\\
{\rm disc}(q_x) &= 4\left(5x-1\right)^2+20\left(-x^2+2x-3-4b\right)= 80\left(x^2+\frac{4}{5}-b\right)\\
&= 80\left(\left(x-\frac{3}{10}\right)\left(x+\frac{3}{10}\right)+\frac{89}{100}-b\right).
\end{align*}
%So for $x \in [0,0.3)$, this is negative.
%That is $, q(x,y)<0$ for all $x \in [0,0.3)$.
Thus, if ${\rm disc}(q_x) \geq 0$, the largest root $y_x$ of $q_x$ is
$$
-5x+1+2\sqrt{5}\cdot \sqrt{x^2+\frac{4}{5}-b}.
$$

Suppose first that $(a,b)=(\frac{19}{25},\frac{89}{100})$.
We claim that $x < \frac{3}{10}$ for all $(x,y) \in P$.
For this, note that ${\rm disc}(p_{\frac{3}{10}})<0$, so $p_{\frac{3}{10}}$ has no real roots, and therefore $p_{\frac{3}{10}}(y)<0$ for all $y \in [0,1]$.
Since $p(x,y)$ is decreasing in $x \geq 0$, we have that $p(x,y)=\frac{1}{40}p_x(y)<0$ for all $x \in [\frac{3}{10},1]$ and $y \in [0,1]$, as required.
On the other hand, we claim that $x\geq \frac{3}{10}$ for all $(x,y) \in Q$.
Indeed, if $x \in [0,\frac{3}{10})$, we have ${\rm disc}(q_x)<0$ and so again $q(x,y)=q_x(y)<0$ for all $y \in [0,1]$.
Thus $P$ and $Q$ are disjoint in this case.

Suppose secondly that $(a,b)=(\frac{3}{4},\frac{19}{20})$.
We claim that $x \leq \frac{1}{2}$ for all $(x,y) \in P$ with $y \neq 0$.
Indeed, $p_x(y)=-y(y+4x-2)$, so we have $y+4x-2 \leq 0$ for $y \in (0,1]$, so $x \leq \frac{1}{2}$.
On the other hand, we claim that $x>\frac{1}{2}$ for all $(x,y) \in Q$.
Indeed, let $x \in [0,\frac{1}{2}]$.
If ${\rm disc}(q_x)<0$ then, since the coefficient of $y^2$ is negative, $(x,y)<0$ for all $y \in [0,1]$.
Otherwise, the largest root $y_x$ of $q_x$ exists, and it is at least $0$ only if $5x^2-10x+4\leq 0$, which is false for all $x \in [0,\frac{1}{2}]$.
So $y_x$ is negative and we see that $q_x(y)<q_x(y_x) = 0$ for all $y \in (y_x,1]$.
Thus $x > \frac{1}{2}$ for all $(x,y) \in Q$, as required.
%But $y_x <0$ for all $x \in [0,\frac{1}{2}]$, and so $q(x,y)<0$ for all $x \in [0,\frac{1}{2}]$, $y \in [0,1]$.
This completes the proof of the claim.
\ecpf
The statement of the lemma then follows easily from the claim, noting that for any $(x,y) \in P \cap Q$ we have $y=0$ and $p(x,0)=\frac{1}{2}(\frac{3}{4}-a) \leq 0$.
\epf

\subsection{Solving Problem $Q^*$ via a linear relaxation}

We define a further optimisation problem.

\medskip
\noindent
\textbf{Problem $L$:}
\it
Given a sequence $\bm{k}:=(k_1,\ldots,k_s)\in \mathbb{N}^s$ of natural numbers,
determine
$
\ell^{\rm max}(\bm{k}) := \max_{\bm{d} \in D(\bm{k})}\ell(\bm{d})
$,
the maximum value of
$$
\ell(\bm{d}) := \sum_{2 \leq t \leq s}\log t \cdot d_t
$$
over the set $D(\bm{k})$ of $(s-1)$-tuples $\bm{d}=(d_2,\ldots,d_s)$
such that $0 \leq d_t \leq 1$ for all $2 \leq t \leq s$, and $\sum_{2 \leq t \leq s}td_t \leq \sum_{c \in [s]}\left(1-\frac{1}{k_c-1}\right)$.
\rm

\medskip
We say that $\bm{d}$ which is feasible for Problem $L$ is \emph{realisable} if there is some $(r,\phi,\ba)\in\feas^*(\bm{k})$ with
\begin{equation}\label{realisable}
d_t=2\sum_{ij \in \binom{[r]}{2}:|\phi(ij)|=t}\alpha_i \alpha_j\quad \text{for all }2 \leq t \leq s
\end{equation}
and call such a feasible triple a \emph{realisation (of $\bm{d}$)}.

We have the following, which is proved by applying Tur\'an's theorem to blow-ups (defined below).
It implies that, in certain special cases, to solve Problem $Q^*$ it suffices to solve Problem $L$.

\begin{lemma}[{\cite[Lemma~5.1]{stability2}}]\label{transfer}
Let $s \in \mathbb{N}$ and $\bm{k} \in \mathbb{N}$.
Then $Q(\bm{k}) \leq \max_{\bm{d} \in D(\bm{k})}\ell(\bm{d})$.
Moreover, the following is true.
Suppose that at least one optimal solution $\bm{d}$ to Problem $L$ is realisable. Then $\max_{\bm{d}\in D(\bm{k})}\ell(\bm{d})=Q(\bm{k})$
and $\opt^*(\bm{k})$ is the set of all $(r,\phi,\ba) \in \feas^*(\bm{k})$ which are realisations of some optimal $\bm{d}$.
\end{lemma}
We wish to add more constraints to Problem $L$.
Indeed, without additional constraints, Problem $L$ only yields realisable solutions in some very special cases, for example $\bm{k}=(k,k)$ or $\bm{k}=(k,k,k)$.
A constraint is \emph{valid} if every $\bm{d}$ which has a realisation $(r,\phi,\ba)\in\opt^*(\bm{k})$ must satisfy the constraint. 
We use $I$ for a set of constraints, each of the type $\sum_{2 \leq f \leq s}a_f d_f \leq b$ for some $a_2,\ldots,a_s,b \in \mathbb{R}$.
%We call this constraint a \emph{$(\mathcal{A},\lambda)$-constraint}.
%A special case is when $\mathcal{A}$ contains some $A \in 2^{[s]}$ if and only if it also contains every $A' \in 2^{[s]}$ with $|A'|=|A|$.
%In this case we can write it as $\sum_{t \in T}d'_t \leq \lambda$ for some $T \subseteq [s]$.
%We only considered this special type of constraint in~\cite{stability2}.
Let Problem $(L,I)$ be Problem $L$ with the constraints in $I$ added to it, and 
let $\ell^{\rm max}_I(\bm{k})$ be the optimal solution of Problem $(L,I)$.
We will still discuss \emph{realisable} solutions $\bm{d}$ and \emph{realisations} of $\bm{d}$ for Problem $(L,I)$ without referring to $I$ when it is clear from the context.

The two types of constraints that we consider are as follows.

\medskip
\noindent
{\it \underline{Universal constraints}.}
Let $\mathcal{A}$ be a set of subsets of $\binom{[s]}{\geq 2}$.
Let $i_2,\ldots,i_s \geq 0$ be such that for each $2 \leq f \leq s$ and $S \in \binom{[s]}{f}$, the number of $A \in \mathcal{A}$ for which $S \in A$ is at least $i_f$. 

Next, given $(r,\phi,\ba)\in\feas^*(\bm{k})$ and $A \subseteq \binom{[s]}{\geq 2}$,
let $H^n_A(\phi,\ba)$ be the `blow-up' graph on $n$ vertices with vertex classes $X_1,\ldots,X_r$ where $|\,|X_i|-\alpha_i n\,| \leq 1$ for all $i \in [r]$ and $xy$ is an edge for $x \in X_i$, $y \in X_j$ if and only if $\phi(ij) \in A$.
Suppose that $d(H^n_A(\phi,\ba)) \leq c_A+o_n(1)$ for all $A \in \mathcal{A}$ and all $(r,\phi,\ba)\in\feas^*(\bm{k})$. Then
$$
i_2 d_{2}+\ldots + i_s d_{s} \leq \sum_{A \in \mathcal{A}}c_A
$$
is a valid constraint.
Indeed, $\sum_{A \in \mathcal{A}}e(H^{n}_A(\phi,\ba)) \geq (i_2d_{2}+\ldots+i_sd_{s})\frac{n^2}{2}+O(n)$.

For example, the basic constraint $\sum_{2 \leq t \leq s}td_t \leq \sum_{c \in [s]}\left(1-\frac{1}{k_c-1}\right)$ is a universal constraint, coming from $\mathcal{A} := \{A_1,\ldots,A_s\}$ and $c_{A_t} := 1-\frac{1}{k_t-1}$, where $A_t := \{S \in \binom{[s]}{\geq 2}: t \in S\}$,
noting that each $S \in \binom{[s]}{f}$ lies in $A_g$ if and only if $g \in S$, so $i_f :=f$.

%Suppose $\bm{d}$ is a feasible solution of Problem $L$ with additional constraints $I$, and $\bm{d}$ has realisation $(r,\phi,\ba)$.
A special case of universal constraint arises when $H_{\mathcal{A}}(\phi) := \{ij \in \binom{[r]}{2}: \phi(ij) \in A \text{ for some }A \in \mathcal{A}\}$ is $K_k$-free for every $(r,\phi,\ba) \in \feas^*(\bm{k})$.
Then $H_{\mathcal{A}}^n(\phi,\ba) := \bigcup_{A \in \mathcal{A}}H_A^n(\phi,\ba)$ is always $K_k$-free, so Tur\'an's theorem implies that $d(H_{\mathcal{A}}^n(\phi,\ba)) \leq 1-\frac{1}{k-1}$.
If $\mathcal{A} := \{\binom{[s]}{t}: t \in T\}$ for some $T \subseteq \{2,\ldots,s\}$, then $i_f=1$ when $f \in T$ and $0$ otherwise, so we have
$$
\sum_{t \in T}d_t \leq 1-\frac{1}{k-1}.
$$
We have the following observations from~\cite{stability2}.
\begin{enumerate}[label=(\Alph*),ref=(\Alph*)]
\item\label{obs1} Suppose there is equality in this constraint. Then there is a partition of $[r]$ into parts $A_1,\ldots,A_{k-1}$ such that $\sum_{i \in A_{i'}}\alpha_{i}=\frac{1}{k-1}$ for all $i' \in [k-1]$, and $ij \in H_{\mathcal{A}}(\phi)$ if and only if $i,j$ lie in different parts $A_{i'},A_{j'}$. 
\item\label{obs2} If $S \subseteq [r]$ has $|S| \leq k$, then $2\sum_{ij \in \binom{S}{2}}\alpha_i\alpha_j \leq \left(1-\frac{1}{k-1}\right)\sum_{i \in S}\alpha_i$
with equality if and only if $\alpha_i=\alpha_j$ for all $ij \in \binom{S}{2}$.
\end{enumerate}

\noindent
{\it \underline{Existential constraints}.}
Let $I$ be a set of valid constraints.
Suppose that there is some $(r,\phi,\ba)\in\feas^*(\bm{k})$ which is the realisation of some feasible $\bm{d}^*$,
and a constraint $a_2d_2+\ldots+a_sd_s \leq b$, and let $I'$ be obtained from adding this constraint to $I$. Suppose that Problem $L$ with constraints $I'$ has optimal value $\ell^{\rm max}_{I'}(\bm{k})< \ell(\bm{d}^*)$.
Then
$$
a_2d_2 + \ldots + a_sd_s \geq b\quad\text{i.e.}\quad -a_2d_2-\ldots-a_sd_s \leq -b
$$ 
is a valid constraint.
Indeed, no optimal solution of Problem $Q^*$ is a realisation of $\bm{d}$ that satisfies the new constraint, so it cannot hold.

\begin{lemma}\label{transfer}
Let $s \in \mathbb{N}$ and $\bm{k} \in \mathbb{N}$.
Let $I$ be a set of valid constraints from Problem $L$ and let $\ell^{\rm max}_I(\bm{k})$ be the optimal value.
Then $Q(\bm{k}) \leq \ell^{\rm max}_I(\bm{k})$.
Moreover, the following is true.
Suppose that at least one optimal solution $\bm{d}$ to Problem $L$ with constraints $I$ is realisable. Then $\ell^{\rm max}_I(\bm{k})=Q(\bm{k})$
and $\opt^*(\bm{k})$ is the set of all $(r,\phi,\ba) \in \feas^*(\bm{k})$ which are realisations of some optimal $\bm{d}$.
\qed\end{lemma}

\begin{figure}\label{hadamardfig}
\begin{center}
\includegraphics{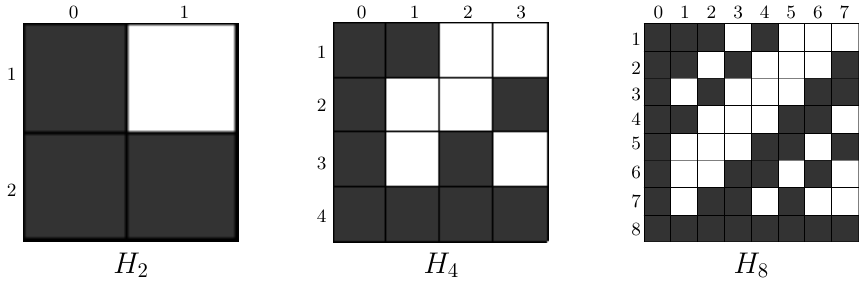}
\end{center}
\caption{Hadamard matrices of order $2,4,8$, which are unique up to equivalence. Here, black represents $1$ and white represents $-1$}
\end{figure}

\subsection{The proof of Theorem~\ref{37thm}.}

First we see why a particular density vector $\bm{d}$ for the $(3;4t-1)$-problem must give rise to a realisation corresponding to the Hadamard matrix $H_{4t}$.

\begin{lemma}\label{hadamard}
Let $t \geq 1$ be an integer and let $\bm{k}=(3;4t-1)$.
Suppose that $\bm{d} = (d_2,\ldots,d_{4t-1})$ with $d_{2t}=1-\frac{1}{4t}$ and all other entries zero is realisable and let $(r,\phi,\ba)\in\feas^*(\bm{k})$ be a realisation.
Then $r=4t$, $\ba$ is uniform, and every $\phi^{-1}(c) \cong K_{2t,2t}$. Moreover,
there is a Hadamard matrix $H_{4t}$ of order $4t$
whose columns are labelled $0,1,\ldots,4t-1$ and rows are labelled $1,2,\ldots,4t$,
normalised so that Column $0$ consists of $1$-entries, and  such that the vertex classes of $\phi^{-1}(c)$ are given respectively by the set of row indices of the $1$-entries, and of the $-1$-entries, in Column $c$ in $H_{4t}$.
\end{lemma}

\bpf
By Tur\'an's theorem, $|\phi^{-1}(c)| \leq \lfloor r^2/4\rfloor$. Also $\bm{d}$ is such that only $d_{2t}$ is non-zero, so $|\phi(ij)|=2t$ for all $ij \in \binom{[r]}{2}$. So $(4t-1)\lfloor r^2/4\rfloor \geq \sum_{c \in [4t-1]}|\phi^{-1}(c)|=2t\binom{r}{2}$.
Solving this yields $r\leq 4t$.
On the other hand, $1-\frac{1}{4t} = d_{2t} = 2\sum_{ij \in \binom{[r]}{2}}\alpha_i\alpha_j$ so Observation~\ref{obs2} implies that $r \geq 4t$.
Thus $r=4t$ and $|\phi^{-1}(c)|=(4t)^2/4$ for all $c \in [4t-1]$ and $\ba$ is uniform.
So the graphs $\phi^{-1}(1),\ldots,\phi^{-1}(4t-1)$, each one isomorphic to $K_{2t,2t}$, decompose $2tK_{4t}$, the complete multigraph on $4t$ vertices with every multiplicity equal to $2t$. That such a decomposition exists and its connection to Hadamard matrices follows from an observation of de Caen, Gregory and Pritikin~\cite{biclique}, as follows.
Let $G_1,\ldots,G_{4t-1}$ be copies of $K_{2t,2t}$ which decompose $2tK_{4t}$.
Construct a $4t\times 4t$ matrix $H$ with the leftmost column consisting of all $1$-s and the $j$-th next column with $1$-s in the rows corresponding to one part of $G_j$ and $-1$-s in the remaining rows. 
For distinct $i,i'$, the number of $j \in [4t-1]$ such that $H_{ij}H_{i'j}=-1$ is exactly the number of $\phi^{-1}(j)$ containing $ii'$, which is $2t$. 
Clearly $H_{i0}H_{i'0}=1$.
Thus for distinct rows $\bm{h},\bm{h}'$ of $H$ we have that the scalar product $\bm{h}.\bm{h}'=2t(-1)+2t(1)=0$,
so the rows are pairwise orthogonal. Since the sum of the squares of entries in each row $\bm{h}$ is $4t$, we see that the collection of $\frac{1}{\sqrt{4t}}\bm{h}$ over rows $\bm{h}$ is an orthonormal basis of $\mathbb{R}^{4t}$ and hence $HH^\intercal = 4tI_{4t}$, so $H$ is a Hadamard matrix.
\epf

\begin{remark}\label{hadamardremark}
\rm
For completeness, we give the other direction of the proof in~\cite{biclique}, that a decomposition can be read off a Hadamard matrix.
Let $H$ be a Hadamard matrix of order $4t$ with columns labelled $0,\ldots,4t-1$ and rows labelled $1,\ldots,4t$.
First note that we can flip the signs of all entries in a given row to get a new Hadamard matrix. Thus we can obtain from $H$ a Hadamard matrix whose first column (of index $0$) is $(1,\ldots,1)^\intercal$.
Since $\frac{1}{\sqrt{4t}}H$ is an orthogonal matrix, every pair of columns are orthogonal (and similarly for rows).
So comparing every column with Column $0$, we see that every other column contains exactly $2t$ $1$-entries. Take any two rows $\bm{h}, \bm{h}'$ with row indices $1 \leq i < i' \leq 4t$ respectively and let $\ell_{ij}$ be the number of columns $j$ where $H_{ij} \neq H_{i'j}$. Since rows are pairwise orthogonal, we have $0 = \bm{h}.\bm{h}' = -\ell_{ij}+(4t-\ell_{ij})=4t-2\ell_{ij}$, so $\ell_{ij}=2t$.
Now for each column $1 \leq j \leq 4t-1$, let $G_j$ be the copy of $K_{2t,2t}$ on vertex set $[4t]$ with vertex classes $A_j,B_j$ where $A_j$ consists of those rows $1 \leq i \leq 4t$ where $H_{ij}=1$.
Given any pair $gh \in\binom{[4t]}{2}$, the number of $G_j$ where $gh$ is an edge is $\ell_{gh}=2t$.
Thus $G_1,\ldots,G_{4t-1}$ decompose $2tK_{4t}$.
\end{remark}

In the next lemma we solve Problem~$L$ for $(3;6)$ and $(3;7)$, using our new tool Lemma~\ref{RB}.

\begin{lemma}\label{37lemma}
Let $\bm{k}=(3;s)$.
\begin{itemize}
\item[(i)] Let $s=6$. Then $Q(\bm{k}) \leq \frac{1}{2}\log 3 + \frac{3}{4}$.
Moreover, if $\bm{d}$ is realisable and satisfies $\ell(\bm{d})= \frac{1}{2}\log 3 + \frac{3}{4}$, then $\bm{d}=(0,\frac{1}{2},\frac{3}{8},0,0)$.
\item[(ii)] Let $s=7$. Then $Q(\bm{k})\leq \frac{7}{4}$.
Moreover, if $\bm{d}$ is realisable and satisfies $\ell(\bm{d})=\frac{7}{4}$, then $\bm{d}=(0,0,\frac{7}{8},0,0,0)$.
\end{itemize}
\end{lemma}

\bpf
By Lemma~\ref{transfer}, it suffices to find a set $I$ of valid constraints for Problem $L$ such that $\ell_I^{\rm max}(\bm{k})$ is at most the required value.
We first prove~(ii). 
Note that
\begin{equation}\label{constraint0}
d_2+\ldots + d_7 \leq 1
\end{equation}
is a valid constraint.
(In fact we could replace $1$ by $1-\frac{1}{R(3;7)-1}$ which is a special case of the universal constraints discussed earlier.)
For each colour $i \in [7]$, let $A_i := \{B \in \binom{[7]}{\geq 2}: i \in B\}$ and 
for each pair $ij \in \binom{[7]}{2}$ of colours, let $A_{ij} := A_i \cup A_j$,
and let $\mathcal{A} := \{A_{ij}:ij \in \binom{[7]}{2}\}$. Then
\begin{equation}\label{Adef}
|\mathcal{A}|=\binom{7}{2} \text{ and } (i_2,\ldots,i_7)=(11,15,18,20,21,21)
\end{equation}
since $i_f=\binom{7}{2}-\binom{7-f}{2}$ for each $2 \leq f \leq 7$.
%So $(i_2,\ldots,i_7)=(11,15,18,20,21,21)$.
Let $(r,\phi,\ba)\in\feas^*(\bm{k})$ be arbitrary and let $n \in \mathbb{N}$ be large.
For brevity write $H^n_i(\phi,\ba)$ for $H^n_{A_i}(\phi,\ba)$ and $H^n_{ij}(\phi,\ba)$ for $H^n_{A_{ij}}(\phi,\ba)$.
Then $H^n_i(\phi,\ba)$ and $H^n_j(\phi,\ba)$ are triangle-free since $\phi^{-1}(i)$ and $\phi^{-1}(j)$ are. Clearly $H^n_{ij}(\phi,\ba)=H^n_i(\phi,\ba) \cup H^n_j(\phi,\ba)$.
%Lemma~\ref{RB} implies that if 

Let $(8,\phi^*,\bm{u})$ be a feasible solution from Remark~\ref{hadamardremark} for $t=2$, so $\bm{u}$ is uniform, $|\phi^*|\equiv 8$ and $q(\phi^*,\bm{u}) = \frac{7}{4}$.
First consider the constraint set $I_1$ consisting of~(\ref{constraint0}) and the single constraint $2d_2+\ldots+7d_7 \leq \frac{339}{100}$.
Multiply~(\ref{constraint0}) by $\log\frac{81}{64}$ and the new constraint by $\log\frac{4}{3}$, and add these together.
(These multipliers come from the (unique) solution to the dual linear program.)
Then each $d_f$ on the left-hand side has coefficient at least $\log f$ (its coefficient in the objective function) and thus $\ell^{\rm max}_{I_1}(\bm{k})$ is at most the right-hand side, which is $\frac{39}{50}+\frac{61}{100}\log 3 < \frac{7}{4}=q(\phi^*,\ba^*)$.
Thus
\begin{equation}\label{constraint1}
2d_2+\ldots+7d_7 \geq \frac{339}{100}
\end{equation}
is a valid (existential) constraint.
Suppose there are $ij \in \binom{[7]}{2}$ and $(r,\phi,\ba)\in\feas^*(\bm{k})$ with $d(H^n_i(\phi,\ba))+d(H^n_j(\phi,\ba)) \leq \frac{89}{100}-\Omega_n(1)$.
Let $\bm{d}$ be such that $(r,\phi,\ba)$ is a realisation of $\bm{d}$.
Since every $d(H^n_{i'}(\phi,\ba)) \leq \frac{1}{2}+o_n(1)$, we have
\begin{equation}\nonumber
2d_2+\ldots + 7d_7 + o_n(1) = \binom{n}{2}^{-1}\sum_{i' \in [7]}d(H^n_{i'}(\phi,\ba)) \leq \frac{5}{2}+\frac{89}{100}-\Omega_n(1) = \frac{339}{100}-\Omega_n(1),
\end{equation}
a contradiction.
Thus Lemma~\ref{RB} implies that $e(H^n_{ij}(\phi,\ba)) \leq \frac{19}{25}n^2/2$ for every $ij$ and $(r,\phi,\ba)$.
With $\mathcal{A}$ defined before~(\ref{Adef}), we have by the above that
\begin{equation}\label{constraint2}
11d_2+15d_3+18d_4+20d_5+21d_6+21d_7 \leq \binom{7}{2}\frac{19}{25}
\end{equation}
is a valid (universal) constraint.
Consider the constraint set $I_2$ consisting of the valid constraint~(\ref{constraint2}) and the new constraint $2d_2+\ldots + 7d_7 \leq \frac{69}{20}$.
Multiply~(\ref{constraint2}) by $\log\frac{32}{27}$ and the new constraint by $\frac{1}{3}\log\frac{9}{8}$, and add these together.
Then each $d_f$ on the left-hand side has coefficient at least $\log f$ and thus $\ell^{\rm max}_{I_2}(\bm{k})$ is at most the right-hand side, which is $\frac{129}{100}+\frac{29}{100}\log 3 < \frac{7}{4}$.
So
\begin{equation}\label{constraint3}
2d_2+\ldots + 7d_7 \geq \frac{69}{20}
\end{equation}
is a valid constraint.
As above, this implies that for every $ij \in \binom{[7]}{2}$ and $(r,\phi,\ba)\in\feas^*(\bm{k})$, which is a realisation of some $\bm{d}$, we have 
$d(H^n_i(\phi,\ba))+d(H^n_j(\phi,\ba)) \geq \frac{69}{20}-\frac{5}{2} + o_n(1) = \frac{19}{20}+o_n(1)$.
Thus Lemma~\ref{RB} implies that $d(H^n_{ij}(\phi,\ba)) \leq \frac{3}{4}+o_n(1)$ for every $ij$.
So
\begin{equation}\label{constraint4}
11d_2+15d_3+18d_4+20d_5+21d_6+21d_7 \leq \binom{7}{2}\frac{3}{4}
\end{equation}
is a valid constraint.
Finally, consider the constraint set $I_3$ consisting of the two valid constraints:~(\ref{constraint4}) and $2d_2+\ldots+7d_7 \leq \frac{7}{2}$, the original universal constraint in Problem $L$.
Multiply the original constraint by $\frac{1}{7}\log\frac{343}{128}$ and~(\ref{constraint4}) by $\frac{1}{21}\log\frac{128}{49}$, and add these together.
Then each $d_f$ on the left-hand side has coefficient at least $\log f$ and thus $\ell^{\rm max}_{I_3}(\bm{k})$ is at most the right-hand side, which is $\frac{7}{4}$.
Moreover, the coefficient of $d_f$ is strictly greater than $\log f$ unless $f=4$, so every optimal solution has $d_f=0$ for all $f \neq 4$.
Thus the unique optimiser has $4d_4 = \frac{7}{2}$, and all other entries equal to $0$.
This completes the proof of~(ii).

The same argument works for $s=6$ to prove (i).
Here $(i_2,\ldots,i_6)=(9,12,14,15,15)$. 
Again, $d_1+\ldots+d_6 \leq 1$ is a valid constraint. We see that $2d_2+\ldots +6d_6 \geq \frac{289}{100}$ is a valid constraint, 
since its negation alone implies that $\ell(\bm{d}) = \sum_{2 \leq f \leq s} fd_f \cdot \frac{\log f}{f} \leq \frac{289}{100}\cdot\frac{\log 3}{3} \leq \frac{1}{2}\log 3 +\frac{3}{4}$.
From this,
the sum of densities of any pair of colour graphs is at most $\frac{289}{100}-4\cdot\frac{1}{2}=\frac{89}{100}$ and so
Lemma~\ref{RB} implies that their union has density at most $\frac{19}{25}$. Thus
\begin{equation}\label{6constraintX}
9d_2+12d_3+14d_4+15d_5+15d_6 \leq \binom{6}{2}\frac{19}{25}
\end{equation}
is a valid constraint. This implies that
\begin{equation}\label{6constraintY}
2d_2+\ldots+6d_6 \geq \frac{59}{20}
\end{equation}
is a valid constraint. 
Indeed, if we add $\frac{1}{3}\log\frac{9}{8}$ times~(\ref{6constraintX}) plus $4-\frac{7}{3}\log 3$ times the negation of~(\ref{6constraintY}), we see $\ell(\bm{d}) \leq \frac{2}{5}+\frac{43}{60}\log 3 < \frac{1}{2}\log 3 + \frac{3}{4}$.
Now the sum of densities of any pair of colour graphs is at most $\frac{59}{20}-4\cdot\frac{1}{2}=\frac{19}{20}$
and so Lemma~\ref{RB} implies that their union has density at most $\frac{3}{4}$. Thus
\begin{equation}\label{6constraintZ}
9d_2+12d_3+14d_4+15d_5+15d_6 \leq\binom{6}{2}\frac{3}{4}
\end{equation}
is a valid constraint. 
Finally, we consider Problem $L$ with the set of two valid constraints:~(\ref{6constraintZ}) and $2d_2+\ldots +6d_6 \leq 3$, the original universal constraint in Problem $L$.
Multiply~(\ref{6constraintZ}) by $\frac{1}{3}\log\frac{9}{8}$ and the original constraint by $4-\frac{7}{3}$ and add these together to see that each $d_f$ on the left-hand side has coefficient at least $\log f$, and is strictly greater than $\log f$ only for $f=2,5,6$, and the right-hand side is $\frac{1}{2}\log 3+\frac{3}{4}$.
Thus the unique optimal solution is $\ell^{\rm max}_I(\bm{k})=\frac{1}{2}\log 3 + \frac{3}{4}$ and the unique optimiser has $d_2=d_5=d_6=0$ and is therefore $\bm{d}=(0,\frac{1}{2},\frac{3}{8},0,0)$.
\epf

\bpf[Proof of Theorem~\ref{37thm}]
This follows immediately from Lemmas~\ref{hadamard} and~\ref{37lemma}. 
\epf

\bpf[Proof of Theorem~\ref{36thm}]
We need to show that every $(r,\phi,\ba) \in \opt^*(\bm{k})$ is as described.
Lemma~\ref{37lemma}(i) implies that $Q(\bm{k}) \leq \frac{1}{2}\log 3 + \frac{3}{4}$, and that, if there is equality, any optimal $(r,\phi,\ba)$ is the realisation of $\bm{d}=(0,\frac{1}{2},\frac{3}{8},0,0)$.
Since $\|\bm{d}\|_1 = \frac{7}{8}$, Observation~\ref{obs2} implies that $r \geq 8$.
Suppose that $|\phi(1i)|=4$ for $i=2,3,4,5$.
Then 
for all distinct $i,i' \in \{2,3,4,5\}$, we have
$\phi(1i) \neq \phi(1i')$, otherwise $|\phi(ii')|\leq 2$ (or there would be a monochromatic triangle).
Let $\psi(ii'):=\{j \in [6]: j \notin \phi(1i) \cap \phi(1i')\}$.
Then $\phi(ii') \subseteq \psi(ii')$.
We use a computer to obtain all $\binom{15}{4}$ possible $\phi(12),\ldots,\phi(15)$ and corresponding $\{\psi(ii')\}$.
Then check each triple $ii'i''$: suppose $\psi(ii'),\psi(i'i''),\psi(i''i')$ share an element $j$.
If they each have size three we have a contradiction since we must delete $j$ from at least one of these sets so we end up with a $\phi$-set of size two.
If one has size $4$ and two have size $3$ then we must delete $j$ from the set of size $4$ (that is, the corresponding $\phi$-sets all have size at most $3$.
After checking all triples, the resulting sets are still supersets of $\phi(ii')$ for distinct $i,i'=2,3,4,5$. In every case, the sets are of the form
$
\{abc, def, xyuv, xyuv',xy'uv,xy'uv'\}$ where $\{a,b,c,d,e,f\}=[6], \{x,y,y'\}=\{a,b,c\}, \{u,v,v'\}=\{d,e,f\}$.
So without loss of generality, the resulting sets are $\{123\},\{456\},\{1245\},\{1246\},\{1345\},\{1346\}$.
Sets $1,3,4$ of sizes $3,4,4$ respectively all contain $1,2$, so without loss of generality we can reduce to
$$
\{123\},\{456\},\{145\},\{246\},\{1345\},\{1346\}.
$$
But now sets $2,4,6$ of sizes $3,3,4$ respectively all contain $4,6$ which is a contradiction because deleting one copy of each gives a $\phi$-set of size at most two.
Thus we have eliminated all possible cases.
We implemented the above in python (\texttt{6check.py}).
Thus when $r \geq 8$, for each $i \in [r]$ there is a set $X_i \subseteq [r]\setminus\{i\}$ of size at most $3$ such that $|\phi(ij)|=4$ if and only if $j \in X_i$.
Therefore $2\cdot 6\lfloor\frac{r^2}{4}\rfloor \geq 2\sum_{c \in [6]}|\phi^{-1}(c)| \geq r\cdot (12+3(r-4)) = 3r^2$. Thus $r$ is even, every $|X_i|=3$, and $\phi^{-1}(c) \cong K_{\frac{r}{2},\frac{r}{2}}$ for all $c \in [6]$.

For $c \in [6]$, let $A_c,B_c$ be the vertex classes of $\phi^{-1}(c)$.
Let $A=(a_{ij})$ be an $r \times 6$ matrix with $\pm 1$ entries, where the $c$-th column represents $\phi^{-1}(c)$, with $a_{jc}=1$ if $j \in A_c$ and $a_{jc}=-1$ if $j \in B_c$.
By relabelling classes, without loss of generality, the $r$-th row of $A$ is $(1,1,1,1,1,1)$.
Notice that $|\phi(ij)|$ is the number of columns $c \in [6]$ where $a_{jc}$ and $a_{ic}$ differ.
So $|\phi(rj)|$ is the number of entries equal to $-1$ in the $j$-th row. But $|\phi(rj)|=3$ for all but $3$ rows $j$.

If $r \geq 10$, then $|\phi(rj)|=3$ for at least $r-1-3 \geq 6$ other rows $j$. Each of these rows, without loss of generality $j=1,2,3,4,5,6$, therefore contains exactly $3$ entries equal to $-1$. By parity, no pair of these rows can differ in exactly three places.
Thus $|\phi(ij)| \neq 3$ and hence $|\phi(ij)|=4$ for all $ij\in\binom{[6]}{2}$, so~e.g.~$|X_1| \geq 5$, a contradiction. 

Recalling that $r \geq 8$ and $r$ is even, we must have $r=8$.
For $\ell =3,4$, let $G_\ell := \{ij \in \binom{[8]}{2}: |\phi(ij)|=\ell\}$.
Since also $\|\bm{d}\|_1=\frac{7}{8}$, $\ba$ is uniform.
One can check via computer (\texttt{6config.py}) that whenever $\phi^{-1}(c) \cong K_{4,4}$ for all $c \in [6]$ are such that every edge multiplicity is $3$ or $4$, then $G_3 \cong K_{4,4}$ and $G_4 \cong K_4 \cup K_4$.
(To reduce computations, first we assume all $\phi^{-1}(c)$ are distinct. Next, if they are not, one needs at least $3=\log 8$ copies of $K_{4,4}$ to cover every edge of $K_8$ at least once, so at least five such graphs are required to cover every edge at least three times. So it remains to check the cases where $\phi^{-1}(1),\phi^{-1}(2)$ are identical and the other colour graphs are distinct.)
Given this structure, $3G_3 \cup 4G_4$ has a decomposition into $6$ copies of $K_{4,4}$ if and only if $4K_8$ has a decomposition into $7$ copies of $K_{4,4}$.
Thus the columns of $A$ are among the $7$ rightmost columns of a normalised order-$8$ Hadamard matrix $H$ (up to permutation of rows and/or negation of rows).
\epf

\subsection{Hadamard matrices and the triangle problem}\label{hadamardsec}

\begin{figure}[h]
\includegraphics[scale=0.88]{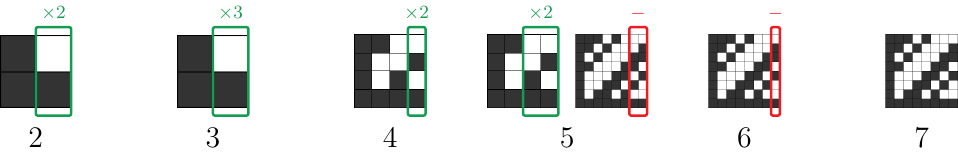}
\caption{Extremal graphs for $(3;s)$ for $2 \leq s \leq 7$}
\label{fighadamard}
\end{figure}

Figure~\ref{fighadamard} shows how the known extremal graphs for the triangle problem relate to Hadamard matrices.
For $s \in [7]\setminus\{5\}$, these are all unique, and for $s=5$ there is an infinite family of asymptotically extremal graphs, of which two are drawn here. Index the columns of the above matrices by $0,1,\ldots$. For each $c \in [s]$, $\phi^{-1}(c)$ is a complete balanced bipartite graph whose vertex classes are given by the set of row indices of the white squares, and of the black squares, in Column $c$. Some columns are repeated, denoted by $\times \ell$, and some are removed, denoted by $-$. Taking $r$ to be the order of the matrix, and $\bm{u} \in \Delta^r$ to be uniform, we have the optimal solution $(r,\phi,\bm{u})$.

\section{Applications of Theorem~\ref{exact1}}\label{applicationsec}

%We have already seen that the strong extension property holds. 
Here we discuss the perfect colouring problem (Problem~\ref{perfprob}), where we have part sizes that maximise the number of perfect colourings. This was implicitly done in previous works (and is trivial in some cases). 
Let $r \in \mathbb{N}$ and let $H \subseteq K_r$ be regular.
Suppose that every $(r^*,\phi^*,\ba^*)\in \opt^*(\bm{k})$ has $r^*=r$, $\ba^* = \bm{u}_r$ uniform and furthermore, there is $M_{\phi^*} \cong H$ such that $|\phi^*(ij)|=t_1$ for all $ij \in M_{\phi^*}$ and $|\phi^*(ij)|=t_2$ for all $ij \in \binom{[r]}{2}\setminus M_{\phi^*}$, where $0 \leq t_1-t_2 \leq 1$. 
Then we denote $\perf(\bm{m}) := \perf_{r,\bm{u}_r}(\bm{m})$.

By symmetry, for some explicit integer $C=C(H)$ we have that
$$
\perf(\bm{m}) = C\cdot \sum_{J \cong H}f(\bm{m},J),\quad\text{where } f(\bm{m},J) := \prod_{ij \in J}t_1^{m_im_j}\prod_{i'j'\notin J}t_2^{m_{i'}m_{j'}},
$$
where the sum is over all copies $J$ of $H$ in $K_r$. Let $\bm{m}':=(m_1-1,m_2,\ldots,m_{r-1},m_r+1)$.
Write $m_{J(i)} := \sum_{j \in N_J(i)}m_j$. Then
\begin{align*}
C^{-1}\perf(\bm{m}') &= \sum_{\substack{J: \{1,r\} \notin J\\ J \cong H}}g_{\rm out}(\bm{m},J)f(\bm{m},J)+ \sum_{\substack{J: \{1,r\} \in J\\ J \cong H}}g_{\rm in}(\bm{m},J) f(\bm{m},J),
\end{align*}
%\begin{align*}
%\perf(\bm{m}') &= \sum_{J: \{1,r\} \notin J}t_1^{-m_{J(1)} + m_{J(r)}}t_2^{m_1-m_r-1+m_{J(1)}-m_{J(r)}}f(\bm{m},J)\\
%&\hspace{2cm}+ \sum_{J: \{1,r\} \in J}t_1^{-m_{J(1)} + m_{J(r)}-1}t_2^{m_1-m_r+m_{J(1)}-m_{J(r)}} f(\bm{m},J).
%\end{align*}
%Using $m_1-m_r \geq \max\{2,m_i-m_j\}$ for any $ij \in \binom{[r]}{2}$, calculations show that the two additional multiplicative terms for $J$ are at least the minimum of
%$$
%t_2, \  \frac{1}{t_1}\left(\frac{t_2^{p-1}}{t_1^{p-2}}\right)^{m_1-m_2}, \ \frac{1}{t_2}\left(\frac{t_2^{2r-4-2p}}{t_1^{2r-4-2p-1}}\right)^{m_1-m_2}
%$$
where
\begin{align}\label{constanteq}
g_{\rm out}(\bm{m},J) &:= \left(\frac{t_2}{t_1}\right)^{m_{J(1)}-m_{J(r)}} t_2^{m_1-m_r-1}\quad\text{and}\quad 
g_{\rm in}(\bm{m},J) := \left(\frac{t_2}{t_1}\right)^{m_{J(1)}-m_{J(r)}} t_2^{m_1-m_r}t_1^{-1}.
\end{align}
%Suppose first that $m_{J(1)}\geq m_{J(r)}$.
%Then for $(i,j)=(2,1)$, this is at least $t_2^{m_r-m_2-1} \geq 2$.
%For $(i,j)=(1,2)$, writing $t_2=\ell-1$ and suppose $J$ is $p$-regular, it is at least
%$$
%\left(\frac{\ell-1}{\ell}\right)^{p(m_1-m_r)}\cdot (\ell-1)^{m_1-m_r-1}.
%$$
%If $m_{J(1)} < m_{J(r)}$ then....
Indeed, if we let $n:= m_1+\ldots+m_r$ and $\bm{m}'':=(m_1-1,m_2,\ldots,m_r)$, then $\perf(\bm{m})/\perf(\bm{m}'') = t_1^{m_{J(1)}}\cdot t_2^{n-m_{J(1)}-m_1}$.
Now combine this with the analogous formula for $\perf(\bm{m}')/\perf(\bm{m}'')$, where e.g.~the power of $t_1$ is $\sum_{j \in N_J(r)}m_j'$ which is $m_{J(r)}$ if $\{1,r\} \notin J$ and $m_{J(r)}-1$ otherwise.

\bpf[Proof of Proposition~\ref{recoverthm}]
Recall that we consider only non-increasing sequences $\bm{m}$.
In each case, letting $r$ be the (unique) number of vertices in every solution in $\opt^*(\bm{k})$, we will show that whenever $m_1-m_r \geq 2$, we have 
$g_{\rm out}(\bm{m},J) \geq 1.01$, say, for all copies $J$ of $H$ with $\{1,r\} \notin J$
and $g_{\rm in}(\bm{m},J) \geq 1.01$ for all $J$ with $\{1,r\} \in J$.
Thus 
we have $\frac{1}{1.01}\cdot \perf(\bm{m}') \geq C\cdot\sum_{J \cong H}f(\bm{m},J) = \perf(\bm{m})$.
This implies that 
the unique supersolution to the perfect colouring problem for $n$ is $\bm{m}^* = (m_1^*,\ldots,m_r^*)$ where $\|\bm{m}^*\|_1=n$ and $m_1^*-m_r^* \leq 1$, and the perfect colouring problem is soluble with constant $c:= 0.01$. If $\bm{k}$ has the strong extension property, Corollary~\ref{exact2} implies that for all sufficiently large $n$, the unique $\bm{k}$-extremal graph on $n$ vertices is $K_{m_1^*,\ldots,m_r^*} \cong T_r(n)$.

All $\bm{k}$ from Theorem~\ref{knownthm} apart from $(3;5)$ have the strong extension property by Lemma~\ref{recover}.
Indeed, we showed in \cite[Theorem~1.7]{stability2} that each $\bm{k}$ among $(k;2), (k;3), (3;4), (4;4)$ has the strong extension property for integers $k \geq 3$, and described the optimal solutions~\cite[Table~2]{stability2}.
Write $g_{\rm in}$ and $g_{\rm out}$ when $\bm{m},J$ are fixed or clear from the context.
Suppose $\bm{k} := (k;s)$ is one of $(k;2),(k;3),(4;4)$. Then every $(r,\phi,\ba)\in\opt^*(\bm{k})$ has $\ba$ uniform and $\phi$ only takes one value $t$. Then for all $\bm{m}$ with $m_1-m_r \geq 2$,~$g_{\rm in} = g_{\rm out} = t^{m_1-m_r-1} \geq t \geq 2$, as required.

Now suppose $\bm{k} := (3;4)$.
In the above discussion, $r=4$, $H \cong K_{2,2}$, so $(t_1,t_2)=(3,2)$. 
Let $J$ be a copy of $H$.
If $\{1,4\}\notin J$, then $J(1)=J(4)=\{2,3\}$ while if $\{1,4\} \in J$, then $J(1)=\{i,4\}$ and $J(4)=\{j,1\}$, where $\{i,j\}=\{2,3\}$.
Let $\bm{m}=(m_1,\ldots,m_4)$ be such that $m_1-m_4 \geq 2$. 
Here, $g_{\rm out} = 2^{m_1-m_4-1} \geq 2$ and $g_{\rm in} = 3^{m_1-m_4-1}\left(\frac{2}{3}\right)^{m_i-m_j}$. If $m_i - m_j \leq 1$ then this is at least $2$; otherwise $2 \leq m_i-m_j \leq m_1-m_4$ and it is at least $\frac{4}{3}$.

Suppose finally that $\bm{k}=(3;6)$.
Here Theorem~\ref{36thm} implies that $r=8$ and $H$ is the disjoint union of two copies of $K_4$. So $(t_1,t_2)=(4,3)$.
Let $J$ be a copy of $H$.
If $\{1,8\} \in J$, then $m_{J(1)}-m_{J(8)} = m_8-m_1$, so $g_{\rm in} = 4^{m_1-m_8-1} \geq 4$.
It remains to check $g_{\rm out}$.
So suppose that $\{1,8\} \notin J$.
Then $J(1), J(8)$ are disjoint sets of size three, and $3m_8 \leq m_{J(1)},m_{J(8)} \leq 3m_1$.
Thus there is an integer $0 \leq j \leq 6(m_1-m_8)$ such that
$
m_{J(1)}-m_{J(8)} = 3(m_1-m_8) - j
$.
So
$$
g_{\rm out} = \left(\frac{81}{64}\right)^{m_1-m_8} \cdot \left(\frac{4}{3}\right)^j \cdot \frac{1}{3}.
$$
%If $m_{J(1)}-m_{J(8)} \leq 1$, then $g_{\rm out} = \frac{3}{4}\cdot 3^{m_1-m_8} \cdot \frac{1}{3} \geq \frac{9}{4}$.
%Otherwise, $g_{\rm out} = \left(\frac{3}{4}\right)^{3(m_1-m_8)-j}3^{m_1-m_8-1}$ for some $0 \leq j \leq 3(m_1-m_8)$,
%where $m_2+\ldots + m_7 = 3(m_1-m_8)-j$.
If $m_1-m_8 \geq 5$ then $g_{\rm out} \geq 1.08$, as required.
If $m_1-m_8=4$ and $j \geq 1$ then $g_{\rm out} \geq 1.14$. %If $j=0$, then $\bm{m}=(d,d,d,d,d-4,d-4,d-4,d-4)$.
If $m_1-m_8=3$ and $j\geq 2$ then $g_{\rm out} \geq 1.2$, and if $m_1-m_8=2$ and $j\geq 3$ then $g_{\rm out} \geq 1.2$.
There are only a few remaining possibilities, namely, writing $\bm{m}=(m_1,\ldots,m_1)+\bm{b}$,
\begin{center}
\begin{tabular}{l|l|l|l}
$m_1-m_8$ & $j$ & $\bm{b}$ & $\bm{b}'$             \\
\hline
$4$       & $0$ & $(0,0,0,0,-4,-4,-4,-4)$ & $(-2,-2,-2,-2,-2,-2,-2,-2)$\\
$3$       & $0$ & $(0,0,0,0,-3,-3,-3,-3)$ & $(-1,-1,-1,-1,-2,-2,-2,-2)$\\
$3$       & $1$ & $(0,0,0,0,-2,-3,-3,-3)$ & $(-1,-1,-1,-1,-1,-2,-2,-2)$\\
             &        & $(0,0,0,-1,-3,-3,-3,-3)$ &\\
$2$       & $0$ & $(0,0,0,0,-2,-2,-2,-2)$ & $(-1,-1,-1,-1,-1,-1,-1,-1)$\\
$2$       & $1$ & $(0,0,0,0,-1,-2,-2,-2)$ & $(0,-1,-1,-1,-1,-1,-1,-1)$\\
             &        & $(0,0,0,-1,-2,-2,-2,-2)$ &\\
$2$       & $2$ & $(0,0,0,0,-1,-1,-2,-2)$ & $(0,0,-1,-1,-1,-1,-1,-1)$\\
             &        & $(0,0,0,-1,-1,-2,-2,-2)$ &\\
             &        & $(0,0,-1,-1,-2,-2,-2,-2)$ &
\end{tabular}
\end{center}
In all cases, we compare by direct calculation to $\bm{\overline{m}} := (m_1,\ldots,m_1)+\bm{b}'$, where $\bm{b}'$ has its sum of entries equal to that of $\bm{b}$, and all entries as equal as possible. Indeed,
$$
\frac{\perf(\bm{\overline{m}})}{\perf(\bm{m})} = \frac{\sum_{A' \cup B' = [8]}\prod_{ij \in \binom{A'}{2}}4^{b_i'b_j'} \prod_{ij \in \binom{B'}{2}}4^{b_i'b_j'} \prod_{i \in A', j \in B'}3^{b_i'b_j'}}{\sum_{A \cup B = [8]}\prod_{ij \in \binom{A}{2}}4^{b_ib_j} \prod_{ij \in \binom{B}{2}}4^{b_ib_j} \prod_{i \in A, j \in B}3^{b_ib_j}}
$$
depends only on $\bm{b},\bm{b}'$.
We did this calculation in python (\texttt{dcheck.py}). In all cases, the ratio is at least $9$.
As before this implies that the unique supersolution has all entries as equal as possible.

Finally, we need to check that $(3;6)$ has the strong extension property.
So let $9$ be a new vertex and let $\phi'$ be an extension of $\phi$ with the property that
$$
\frac{1}{8}\sum_{i \in [8]}\log t_i = \frac{3}{4}+\frac{1}{2}\log 3 \quad\text{i.e.}\quad t_1 \ldots t_8 = 2^{6}3^4,
$$
where $t_i := \max\{|\phi'(i9)|,1\}$. 
Since each $\phi^{-1}(c) \cong K_{4,4}$ and $\phi'^{-1}(c)$ is triangle-free, in this graph $9$ has at most $4$ neighbours, and $\phi'^{-1}(c) \subseteq K_{5,4}$.
Let $G_c$ be the copy of $K_{5,4}$ containing $\phi'^{-1}(c)$, such that $G_c-\{9\} \cong K_{4,4}$, so we obtain $G_c$ from $\phi^{-1}(c)$ by adding $9$ to one of the two vertex classes of $\phi^{-1}(c)$.
One can easily check all such attachments using a computer, trying each of the $2^6$ choices (we implemented this using python (\texttt{7ext.py})).
This reveals that $\{G_c:c \in [6]\}$ has $\prod_{i \in [8]}|\{c \in [6]: i9 \in G_c\}|= 2^6 3^4$ if and only if $\phi^{-1}(c) \equiv K_{5,4}$ and there is some $i \in [8]$ such that $9$ lies in exactly the same vertex class as $i$ in each $\phi^{-1}(c)$; that is, $9$ is a strong clone of $i$ (in other words, $9$ corresponds to a copy of an existing row in the $\pm 1$-matrix representing the optimal solution).
Thus $\phi'^{-1}(c) = G_c$.
We have proved that $(3;6)$ has the strong extension property.

Therefore we can apply Theorem~\ref{exact1} with Theorem~\ref{36thm} to see that $T_8(n)$ is the unique extremal graph for sufficiently large $n$.
Moreover, letting $0 \leq j \leq 7$ be such that $n=8N+j$ for $N \in \mathbb{N}$, we see that the unique supersolution $\bm{m}$ has $j$ terms equal to $N+1$ and $8-j$ equal to $N$, and
$\perf(\bm{m})$ equals $C_j' \cdot 4^{12N} \cdot 3^{16N}$ where $C_j'$ is a constant depending only on $j$, which can be determined by calculating $|\pat(8,\bm{u})|$, by counting certain distinct Hadamard matrices (see the end of the proof of Theorem~\ref{trianglethm}). This gives the claimed value of $F(n;(3;6))$.
\epf

\bpf[Proof of Theorem~\ref{trianglethm}]
Let $\bm{k}:=(3;7)$.
By Theorem~\ref{37thm}, every element of $\opt^*(\bm{k})$ is of the form $(8,\phi,\bm{u})$, where $\bm{u}$ is uniform and $|\phi(ij)|=4$ for all $ij \in \binom{[8]}{2}$.

We need to show that $\bm{k}$ has the strong extension property.
So let $9$ be a new vertex and let $\phi'$ be an extension of $\phi$ with the property that
$$
\frac{1}{8}\sum_{i \in [8]}\log t_i = \frac{7}{4} \quad\text{i.e.}\quad t_1 \ldots t_8 = 2^{14}
$$
where $t_i := \max\{|\phi'(i9)|,1\}$.
So every $t_i$ equals $1$, $2$, or $4$, and in fact the multiset of values is either
$\{4,4,4,4,4,4,4,1\}$ or $\{4,4,4,4,4,4,2,2\}$.
Since each $\phi^{-1}(c) \cong K_{4,4}$ and $\phi'^{-1}(c)$ is triangle-free, in this graph $9$ has at most $4$ neighbours, and $\phi'^{-1}(c) \subseteq K_{5,4}$.
Let $G_c$ be the copy of $K_{5,4}$ containing $\phi'^{-1}(c)$, such that $G_c-\{9\} \cong K_{4,4}$, so we obtain $G_c$ from $\phi^{-1}(c)$ by adding $9$ to one of the two vertex classes of $\phi^{-1}(c)$.
One can easily check all such attachments using a computer, trying each of the $2^7$ choices (we implemented this using python (\texttt{7ext.py})).
This reveals that $\{G_c:c \in [7]\}$ has $\prod_{i \in [8]}|\{c \in [7]: i9 \in G_c\}|= 2^{14}$ if and only if $\phi^{-1}(c) \equiv K_{5,4}$ and there is some $i \in [8]$ such that $9$ lies in exactly the same vertex class as $i$ in each $\phi^{-1}(c)$; that is, $9$ is a strong clone of $i$ (in other words, $9$ corresponds to a copy of an existing row in the Hadamard matrix representing the optimal solution).
Thus $\bm{k}$ has the strong extension property.

Since $|\phi|$ only takes the value $4$, the same argument as at the beginning of the proof of Proposition~\ref{recoverthm} implies that the perfect colouring problem is soluble (with constant $c:= 4$) and for all integers $n$, the unique supersolution $\bm{m}_n=(m_{n,1},\ldots,m_{n,8})$ of the perfect colouring problem has $|m_{n,i}-m_{n,j}| \leq 1$ for all $ij$.
Corollary~\ref{exact2} implies that the unique extremal graph is $K_{\bm{m}_n}(n) \cong T_8(n)$.
Let $C$ be the number of Hadamard matrices with first column and last row consisting of $1$-s. Then $F(n;\bm{k}) = (C+o(1))\cdot 4^{t_8(n)}$.
Note that one can calculate $C$ if desired. Indeed, two Hadamard matrices are \emph{equivalent} if one can be obtained from the other by a sequence of permutations and/or negations of rows and/or columns. It is known that, up to equivalence, there is one Hadamard matrix $H$ of order $8$.
 So $C$ is the number of distinct matrices that can be obtained from $H$ by permuting the first $7$ rows and the last $7$ columns and making the first column and the last row consist of all $1$-s.
\epf

\section{The $(k+1,k)$-extremal graphs}\label{2colsec}

In this section, we consider the simplest case when $\bm{k}$ has the weak extension property, namely $\bm{k} = (k+1,k)$ for $k \geq 3$.
For small $k$, we determine the $(k+1,k)$-extremal graph, which turns out to have a part of size $O(k)$, and the size of this part depends on the value of $n$ modulo $k-1$.
The proof relies heavily on Theorem~\ref{stabilitycomp}, the full strength of which we have not yet required.

\begin{definition}
\rm
For all integers $k \geq 3$ and $n \geq 4$,
\begin{itemize}
\item for $q \geq 2$, let $\mathcal{J}_q(n)$ be the family of complete $k$-partite graphs with parts of size $m_1 \geq \ldots \geq m_{k-1} \geq \ell$ where $m_1-m_{k-1} \leq q$, $2 \leq \ell \leq q$ and $m_1+\ldots+m_{k-1}+\ell=n$;
\item let $\mathcal{J}_q^*(n)$ be the set of graphs in $\mathcal{J}_q(n)$ which are $(k+1,k)$-extremal;
\item for $\ell \geq 1$, let $J^\ell(n)$ denote the complete $k$-partite graph with parts of size $m_1 \geq \ldots \geq m_{k-1} \geq \ell$ where $m_1-m_{k-1} \leq 1$ and $m_1+\ldots + m_{k-1}+\ell=n$.
\end{itemize}
\end{definition}

In this section we prove the following theorem, which includes Theorem~\ref{2colapproxthm}.

\begin{theorem}\label{kk+1thm}
For all integers $k \geq 3$, there exists $n_0>0$ such that for all $n>n_0$, every $n$-vertex $(k+1,k)$-extremal graph is in $\mathcal{J}_{2(k-1)}(n)$. Moreover, $F(n;k+1,k) = O_k(1) \cdot 2^{t_{k-1}(n)}$, where the factor $O_k(1)$ is at least $2$.
\end{theorem}

Recall from~(\ref{FQeq}) that $F(k+1,k)=Q(k+1,k)$.
It is easy to solve Problem $Q^*$ for $s=2$ colours (see~\cite[Lemma 1.8]{stability2}): the unique solution is $(k-1,\phi,\bm{u})$ where $\bm{u}$ is uniform and $\phi(ij)=[2]$ for all $ij$.
Thus when $k$ is fixed and $n \to \infty$,
$$
Q(k+1,k) = \textstyle{\frac{k-2}{k-1}}\quad\text{and}\quad \log F(n;k+1,k) = \left(\textstyle{\frac{k-2}{k-1}}+o(1)\right)\binom{n}{2}.
$$
So, while~(\ref{FQeq}) determines $ \log F(n;k+1,k)/\binom{n}{2}$ asymptotically, Theorem~\ref{kk+1thm} determines $F(n;k+1,k)$ up to a multiplicative constant.

The next proposition is numerical.

\begin{proposition}\label{numerical}
Let $k \geq 3$ and $0 \leq j \leq k-2$ be integers. 
Given $m_1 \geq \ldots \geq m_{k-1} \geq \ell$, let $\bm{m} := (m_1,\ldots,m_{k-1})$ and
\begin{align*}
h(\bm{m};\ell) &:= \left(\prod_{ii' \in \binom{[k-1]}{2}}2^{m_im_{i'}}\right) 2^{\ell(m_1+\ldots+m_{k-1})} \left( 2^{-m_1}+ \ldots +2^{-m_{k-1}}\right)^{\ell},\\
h^*(\bm{m};\ell) &:= \left(\prod_{ii' \in \binom{[k-1]}{2}}2^{m_im_{i'}}\right) 2^{\ell(m_1+\ldots+m_{k-1})} \left( \sum_{i \in [k-1]}2^{-m_i} - \sum_{ii' \in \binom{[k-1]}{2}}2^{-m_i-m_{i'}}\right)^{\ell}.
\end{align*}
Let
$$
f(k,j,\ell) := \frac{h(\bm{m}^*;\ell)}{2^{t_{k-1}(n)}}\quad\text{and}\quad f(k,j) := \max_{0 \leq \ell \leq k-1}f(k,j,\ell),
$$
where $m^*_1 \geq \ldots \geq m^*_{k-1} \geq \ell$, $m^*_1 - m^*_{k-1} \leq 1$ and $m^*_1+\ldots + m^*_{k-1}+\ell=n$ (so $m^*_1,\ldots,m^*_{k-1},\ell$ are the part sizes of $J^\ell(n)$).
Let $\bm{m}^*:=(m^*_1,\ldots,m^*_{k-1})$.
\begin{enumerate}[label=(\roman*),ref=(\roman*)]
\item\label{2coli} Let $n \equiv j  ~(\!\!\!\!\mod k-1)$ with $0 \leq j < k-1$ and $n \geq k^3$.
Consider the following way of colouring $J \in \mathcal{J}_{n}(n)$ with parts of size $m_1 \geq \ldots \geq m_{k-1} \geq \ell$: for each $x$ in the $\ell$-part, pick $i_x \in [k-1]$ and colour every edge between $x$ and the $m_{i_x}$-part with colour $1$. All other edges are coloured arbitrarily. The number of colourings
is between $h^*(\bm{m};\ell)$ and $h(\bm{m};\ell)$.
\item\label{2coli.5}
If $0 \leq \ell \leq k-1$, then
$$
f(k,j,\ell) =  \begin{cases} 
%2^{\ell(j-\ell)+\binom{j-\ell}{2}-\binom{j}{2}}
2^{-\binom{\ell}{2}}\left(k-1-\frac{j-\ell}{2}\right)^\ell &\mbox{if } 0 \leq \ell \leq j \\
%2^{-\ell(\ell-j)+\binom{\ell-j}{2}-\binom{j}{2}}
2^{j-\binom{\ell+1}{2}}\left(k-1+\ell-j\right)^\ell &\mbox{if } j < \ell \leq k-1.
\end{cases}
$$
\item\label{2colii} Let $n \equiv j ~(\!\!\!\!\mod k-1)$ with $n \geq k^3$. For $3 \leq k \leq 10$, the maximum of $h(\bm{m};\ell)$ over all (part sizes of) $J \in \mathcal{J}_{2(k-1)}(n)$ is attained by $J^{\ell(k,j)}(n)$, where $\ell(k,j)$ is the $(k,j)$-th entry in the table below:

\rm
\begin{tabular}{r|rlllllllllll}
    & $j=0$ & 1 & 2 & 3 & 4 & 5 & 6 & 7 & 8 & \\ \hline
$k=3$  & 2 & 2 &   &   &   &   &   &   &   &   &    &    \\
4  & 3 & 2 & 2 &   &   &   &   &   &   &   &    &    \\
5  & 3 & 3 & 2 & 3 &   &   &   &   &   &   &    &    \\
6  & 3 & 3 & 3 & 3 & 3 &   &   &   &   &   &    &    \\
7  & 3 & 3 & 3 & 3 & 3 & 3 &   &   &   &   &    &    \\
8  & 3 & 3 & 3 & 3 & 4 & 4 & 3 &   &   &   &    &    \\
9  & 3 & 3 & 3 & 3 & 4 & 4 & 4 & 4 &   &   &    &    \\
10 & 4 & 3 & 3 & 3 & 4 & 4 & 4 & 4 & 4 &   &    &    \\
\end{tabular}
\it 

and the value is separated from the value for any other $J \in \mathcal{J}_{2(k-1)}(n)$ by a factor $1+\Omega(1)$.
%Moreover,
%$$
%F(n;k+1,k) = (1+o(1))f(k,j) \cdot 2^{t_{k-1}(n)}.
%$$
%is the $(k,j)$-th entry in the table below:
%
%\rm
%\begin{tabular}{r|rlllllllllll}
 %   & $\ell=0$ & 1 & 2 & 3 & 4 & 5 & 6 & 7 & 8 & 9 & 10 & \!\!\!11 \\ \hline
%$k=3$  & 2 & $\frac{3^2}{2^2}$ &   &   &   &   &   &   &   &   &    &    \\
%4  & $\frac{3^3}{2^2}$ & $2^2$ & $\frac{3^2}{2}$ &   &   &   &   &   &   &   &    &  \\
%5  & $\frac{7^3}{2^6}$ & $\frac{3^3}{2^2}$ & $2^3$ & $2^3$ &   &   &   &   &   &   &    &   \\
%6  & $2^3$ & $\frac{7^3}{2^5}$ & $\frac{3^3}{2}$ & $\frac{5^3}{2^3}$ & $\frac{3^6}{2^6}$ &   &   &   &   &   &    &  \\
%7  & $\frac{3^6}{2^6}$ & $2^4$ & $\frac{7^3}{2^4}$ & $3^3$ & $\frac{11^3}{2^6}$ & $\frac{5^5}{2^3}$ &   &   &   &   &    &   \\
%8  & $\frac{5^5}{2^3}$ & $\frac{3^6}{2^5}$ & $2^5$ & $\frac{7^3}{2^3}$ & $\frac{7^4}{2^6}$ & $\frac{13^4}{2^{10}}$ & $\frac{11^3}{2^6}$ &  &   &   &    &
%\end{tabular}
%\it
\item\label{2coliii} For all $0 \leq j \leq k-2$ we have $f(k,j) \geq \max\{(k+1)^2/8,2\}$.
\end{enumerate}
\end{proposition}

\bpf[Proof sketch]
For~\ref{2coli}, the assertion about $h$ is clear.
Further, $h$ overcounts those colourings with all colour-$1$ edges to at least two other parts.
By a Bonferroni-type inequality, $h^*$ is a lower bound.
%It remains to prove~(ii).
%To determine $\mathcal{J}_k^*(n)$, one needs to calculate at most $k(k+1)$ possible values.
%Indeed, fix $\ell$, and suppose that we want to determine this set for $n \equiv \ell \mod k-1$.

Let us turn to the remaining claims.
Take any admissible $(\bm{m};\ell)$.
Define the integer $b$ by $n =b(k-1)+j$ and $0 \leq j < k-1$.
Define $b_i := m_i-b$.
Then $b_1 - b_{k-1} \leq 2(k-1)$.
Moreover, $b_1+\ldots + b_{k-1} + \ell = \sum_i m_i - (k-1)b + \ell = n - (k-1)b = j$.
Then
\begin{equation}\label{fdceq}
h(\bm{m};\ell) = 2^{\binom{k-1}{2}b^2+(k-2)bj} \cdot 2^{\ell\sum_{i}b_i} \left( 2^{-b_1}+\ldots + 2^{-b_{k-1}}\right)^\ell \prod_{ii'} 2^{b_ib_{i'}}.
\end{equation}
First we prove~\ref{2colii}. For this, one needs to compare the values of
\begin{equation}\label{tildeh}
\tilde{h}(\bm{b};\ell) := 2^{\ell\sum_{i}b_i}(2^{-b_1}+\ldots+2^{-b_{k-1}})^\ell \prod_{ii'}2^{b_ib_{i'}}
\end{equation}
among all tuples $\bm{b} := (b_1,\ldots,b_{k-1})$ and $\ell$ with sum $j$ and $|b_i-b_{i'}|,\ell \leq 2(k-1)$.
We implemented this in python (\texttt{smallpart.py}), and the optimal value of $\ell$ (which is indeed unique) for $3 \leq k \leq 10$ was recorded in the table.
In all cases, we also had $b_1-b_{k-1} \leq 1$, so the optimal $m_1,\ldots,m_{k-1}$ are as equal as possible. Thus the maximum is attained by some $J^\ell(n)$.
For future reference, we note that we could have defined $b_i$ by subtracting any $b' \leq b$ from $m_i$, and then to compare values of $h(\bm{b};\ell)$, one must compare values of $\tilde{h}(\bm{b};\ell)$ over tuples $\bm{b}$ with sum $j+(k-1)(b-b')$.
%Substituting into~(\ref{fdceq}) gives the stated expression $f(k,j,\ell)$. This completes the proof of~(ii).

Note that
$$
t_{k-1}(n)
=t_{k-1}((k-1)b+j)
=\binom{k-1}{2}b^2+(k-2)bj+\binom{j}{2}.
$$
For $0 \leq \ell \leq j$,~(\ref{fdceq}) implies that the function $f(k,j,\ell)$ equals $2^{-\binom{j}{2}}\cdot\tilde{h}((1,\ldots,1,0,\dots,0);\ell)$ where there are $j-\ell$ entries of value $1$.
For $j < \ell \leq k-1$, it equals $2^{-\binom{j}{2}}\cdot\tilde{h}((0,\ldots,0,-1,\ldots,-1);\ell)$ where there are $\ell-j$ entries of value $-1$.
This proves~\ref{2coli.5}.
%Every colouring described in~\ref{2coli} is valid. Writing $n_1,\ldots, n_{k-1}$ for the size of the large parts, the number of choices is
%The relation before $h$ and $f$ follows by calculation.
%For~\ref{2coli.5} we take $\bm{b}:=(1,\ldots,1,0,\ldots,0)$ with $j-\ell$ $1$-entries when $j-\ell \geq 0$,
%and we take $\bm{b}=(0,\ldots,0,-1,\ldots,-1)$ with $\ell-j$ $-1$-entries when $0 \leq \ell-j \leq k-1$.

For~\ref{2coliii}, we bound $f(k,j,2)$ from below. For $j \geq 2$, we have
$$
f(k,j,2) = \frac{1}{2}\cdot (k-j/2)^2 \geq \frac{1}{2} \cdot (k-(k-2)/2)^2 = \frac{(k+2)^2}{8}.
$$
We have $f(k,0,2) = (k+1)^2/8$ and $f(k,1,2) = k^2/4$.
For $k \geq 11$ (i.e.~not in the table), the smallest of these is $(k+1)^2/8 \geq 18$.
\epf

The next lemma states that whenever every vertex and pair of vertices in a near-extremal complete partite graph $G$ has almost optimal contribution to $F(G;\bm{k})$, then $G$ is only optimal if it lies in $\mathcal{J}_{2(k-1)}(n)$.

\begin{lemma}\label{close}
There exist $\eps > 0$ and $n_0 \in \mathbb{N}$ such that the following holds.
Let $G$ be a graph on $n \geq n_0$ vertices which is complete multipartite, and such that $\log F(G;k+1,k) \geq \left(\frac{k-2}{k-1}-\eps\right)\binom{n}{2}$
and
\begin{align}\label{43def}
&\log F(G;k+1,k) - \log F(G-x;k+1,k) \geq (\textstyle{\frac{k-2}{k-1}}-2\eps)n\quad\forall \ x \in V(G),\\
\nonumber &\log F(G;k+1,k) - \log F(G-y-z;k+1,k) \geq (\textstyle{\frac{k-2}{k-1}}-2\eps)(n+n-1)\\
\nonumber &\hspace{6cm}\forall\text{ distinct }y,z \in V(G).
\end{align}
If $G \not\in \mathcal{J}_{2(k-1)}(n)$, then there is a graph $G'$ of order $n$ with $F(G';k+1,k) > F(G;k+1,k)$.
\end{lemma}

\bpf
As it is easy to show (or see~\cite[Lemma~1.8]{stability2}), $(k+1,k)$ has the extension property.
Let $0 < \delta' \ll \delta$.
Let $n_0,\eps>0$ be the parameters returned by Theorem~\ref{stabilitycomp} applied with parameter $\delta'$.
We may assume that $0 < 1/n_0 \ll \eps \ll \delta'$. Let $G$ be a graph on $n \geq n_0$ vertices as in the lemma.
Let $X_1,\ldots,X_t$ be the partition of $V(G)$, where $|X_1| \geq \ldots \geq |X_t|$.
By Theorem~\ref{stabilitycomp}, there is a set $\mathcal{P}$ of at least $(1-2^{-\eps n})F(G;k+1,k)$ valid colourings of $G$ such that, for each $\chi \in \mathcal{P}$, the following hold. 
\begin{itemize}
\item[(i)]  There is a coarsening $Z_0,Z_1,\ldots,Z_{k-1}$ of $X_1,\ldots,X_t$ such that $\sum_{i \in [k-1]}|\,|Z_i|-\frac{n}{k-1}\,| < \delta'n$.
\item[(ii)] Each of $Z_1,\ldots,Z_{k-2}$ is a part of $G$ and $Z_{k-1}$ contains at most two parts of $G$.
\item[(iii)] For all $x \in V(G)$, there is $i(x) \in [k-1]$, so that we have $\chi(xy) = 1$ for all $y \in N_G(x) \cap Z_{i(x)}$,
and for all $j \in [k-1]\setminus\{i(x)\}$, it holds that
\begin{equation}\label{eqvxgood}
|\chi^{-1}(c)[x,X]|=\textstyle{\frac{1}{2}}|X|\pm\delta'|Z_j|\quad\text{for all }c \in [2] \text{ and parts }X \subseteq Z_j. 
\end{equation}
Moreover, for all $j \in [k-1]$ and distinct $y,z \in V(G)$ with $i(y)=i(z) \neq j$, it holds that
\begin{equation}\label{eqpairgood}
|N_{\chi^{-1}(c)}(y,X) \cap N_{\chi^{-1}(c')}(z,X)| = \textstyle{\frac{1}{4}}|X| \pm \delta'|Z_j|\quad\text{for all }c,c' \in [2]\text{ and parts }X \subseteq Z_j.
\end{equation} 
Also, if $x \in Z_i$ for $i \in [k-1]$ then $i(x) = i$.
\end{itemize}

Fix $\chi^* \in \mathcal{P}$ and choose $Z_0,\ldots,Z_{k-1}$ and $i : V(G) \mapsto [k-1]$ as above.
(Recall that these may be different for different colourings.)
If there is a part in $Z_{k-1}$ of size at most $\delta n$, move it into $Z_0$.
Then~(i) holds with parameter $2\delta$,~(iii) holds with parameter $\delta'$, every part in $Z_{k-1}$ has size at least $\delta n$, and all but one part in $Z_0$ has size at most $\delta'n$.
Let $p^* \leq 2$ be the number of parts in $Z_{k-1}$ and $q^*$ the number of parts in $Z_{0}$ (which both depend on $\chi^*$).

\begin{claim}\label{1stclaim}
Suppose that $Z_0 \neq \emptyset$.
Then $Z_{k-1}$ contains exactly one part.
Moreover, if $x,x' \in Z_0$ lie in different parts of $G$, then for all $\chi \in \mathcal{P}$ we have $\chi(xx') = 2$ and $i(x) \neq i(x')$.
In particular, $Z_0$ contains at most $k-1$ parts.
\end{claim}

\bcpf
Let $\{ B_{p'} : p' \in [p^*] \}$ be the part(s) of $G$ contained in $Z_{k-1}$.
Suppose first that $p^*=2$ and $q^* \geq 1$.
Let $x \in Z_0$ and $\chi \in \mathcal{P}$ be arbitrary. %, and suppose without loss of generality that $i(x) = 1$.
Then $\chi(xv)=1$ for all $v \in Z_{i(x)}$, and for each $i \in [k-1]\setminus\{i(x)\}$, the vertex and pair locally good conditions~(\ref{eqvxgood}) and~(\ref{eqpairgood}) hold.
Without loss of generality, assume $|B_1| \geq |B_2|$. Thus e.g.~$|B_1| \geq n/(2k)$ and $|B_2| \geq \delta n$.
By~(\ref{eqvxgood}) there is $y \in B_2$ such that $\chi(xy)=1$ (when $i(x) \neq k-1$, since if $i(x)=k-1$ this is obvious).
Since $i(y)=k-1$ we also have $\chi(yz)=1$ for every $z \in B_1$.
Let $J := \chi^{-1}(1)$.
For $i \in [k-2]$, let $A_i := N_J(x,Z_i) \cap N_J(y,Z_i)$ and let $A_{k-1} := N_J(x,B_1) \cap N_J(y,B_1)$.
Then~(\ref{eqpairgood}) implies that $|A_i| \geq (\frac{1}{4}-\delta')|Z_i| \geq |Z_i|/5 \geq n/(5k)$ for all $i \in [k-2]$.
Also, by~(\ref{eqvxgood}), $|A_{k-1}| = |N_J(x,B_1)| \geq (\frac{1}{2}-\delta')|B_1| \geq |Z_{k-1}|/5 \geq n/(5k)$.
Lemma~\ref{adjust} implies that $J[A_i,A_j]$ is $(\sqrt{\delta'},\frac{1}{2})$-regular for all distinct $i,j \in [k-1]$.
By Lemma~\ref{embed}, $J$ contains a $K_{k-1}$ with one vertex in each of $A_1,\ldots,A_{k-1}$.
Together with $x,y$, this gives rise to a $1$-coloured copy of $K_{k+1}$ in $G$, a contradiction.
So if $p^*=2$, then $q^* = 0$.

Suppose now that $p^*=1$ and let $x,x' \in Z_0$ lie in different parts of $G$.
Similar arguments to those above show that $x,x'$ form the $2$-element set in a copy of the $k$-partite graph $K_{1,\ldots,1,2}$ in $\chi^{-1}(1)$.
Therefore (since the edge between them must have some colour) $\chi(xx') = 2$, as required.
Suppose now, for a contradiction, that $i(x) = i(x')$.
Without loss of generality, suppose their common value is $1$.
Then $U_i := N_{\chi^{-1}(2)}(x,Z_i) \cap N_{\chi^{-1}(2)}(x',Z_i)$ satisfies $|U_i| = (\textstyle{\frac{1}{4}} \pm \sqrt{\delta'})|Z_i|$ for all $2 \leq i \leq k-1$.
By Proposition~\ref{badrefine}, $\chi^{-1}(2)[U_i,U_{i'}]$ is $((\delta')^{1/3},\frac{1}{2})$-regular for all distinct $2 \leq i,i' \leq k-1$.
By Lemma~\ref{embed}, there are $z_i \in U_i$ for $2 \leq i \leq k-1$ such that
$x,x',z_2,\ldots,z_{k-1}$ form a copy of $K_{k}$ in $\chi^{-1}(2)$, a contradiction.
Thus $i(x) \neq i(x')$.
This completes the proof of the claim.
\ecpf

\medskip
\noindent
%From now on we will write $C_i=Z_i$ if $|Z_i| \leq \sqrt{\delta}$, and by a slight abuse of notation will change the meanings of $\ell,p$ slightly. Now, 
%There are four cases: (1) $p=1$ and $\ell=0$; (2) $p=2$ and $\ell=0$; (3) $p=1$ and $\ell=1$; (4) $p=1$ and $\ell \geq 2$.
%We will consider (1)--(4) and give an upper bound for $|\mathcal{P}|$ in each case.
%Write $\mathcal{P}_{p',\ell'}$ or $\mathcal{P}_{p',\ell'}(d_1,\ldots,d_{k-2+p'+\ell'})$ for the set $\mathcal{P}$ in the case when $G$ has $p'=p$ and $\ell'=\ell$ (so $G$ has $k-2+p'+\ell'$ parts), and $|X_i| = d_i$ for all $i$.

The total number of parts is $t = k-2+p^*+q^*$.
Suppose that $t=k-1$. Then 
$G$ is $K_k$-free.
Thus every $2$-edge colouring is $(k,k+1)$-valid and so
$|\mathcal{P}| \leq 2^{t_{k-1}(n)}$ and hence $F(G;\bm{k}) \leq (1-2^{-\eps n})^{-1}\cdot 2^{t_{k-1}(n)}$.

Suppose $t \geq k$ and every part of $G$ has size at least $\delta n$.
Then every $\chi \in \mathcal{P}$ has the same associated coarsening $Z_1=X_1,\ldots,Z_{k-2}=X_{k-2},Z_{k-1}=X_{k-1} \cup X_k$ in order to satisfy $\sum_{i \in [k-1]}|\,|Z_i|-\frac{n}{k-1}\,| \leq \delta' n$.
The number of perfect colourings is the same as if $X_{k-1},X_k$ were merged into a single part since edges between these parts are coloured with $1$, so
$$
|\mathcal{P}| \leq 2^{t_{k-1}(n)},\quad\text{so again}\quad F(G;\bm{k}) \leq (1+o(1)) \cdot 2^{t_{k-1}(n)}.
$$
We will see later, due to
Proposition~\ref{numerical}, that there is $G'$ with $|\mathcal{P}(G')| \geq 1.9|\mathcal{P}(G)|$.

So from now on we will assume that $t \geq k$ and there are $q \geq 1$ parts of size less than $\delta n$.
Claim~\ref{1stclaim} implies that these parts lie in $Z_0$ so $p^*=1$ and $q^*=q=t-(k-1)$.
Thus we have $|X_i| < \delta n$ for all $i \geq k$, so also $|X_1| \geq \ldots \geq |X_{k-1}| \geq \frac{n}{k-1}-(k-1)\delta n$.
%Recall that $q = t-k+2-p$ is the number of parts in $Z_{0}$, and write $L_i := X_{k-1+i}$ for all $i \in [q]$. 
%So $|L_1|$ could be as large as $\delta n$, 
We will say that $X_1,\ldots,X_{k-1}$ are the \emph{large parts}, of sizes $m_1,\ldots,m_{k-1}$ respectively, and $L_1,\ldots,L_q$ are the \emph{small parts}, of sizes $\ell_1,\ldots,\ell_q$ respectively. Also, let $m := m_1+\ldots + m_{k-1}$ and $\ell := \ell_1 + \ldots + \ell_q$, and $\bm{m} := (m_1,\ldots,m_{k-1})$ and $\bm{\ell} := (\ell_1,\ldots,\ell_q)$.
For some colourings $\chi$, one of the $L_i$ may be a part of $Z_{k-1}$, while in others all $L_i$ will be parts of $Z_0$.
Suppose there are $\chi$ and $i$ such that $Z_{k-1}=X_{k-1} \cup L_i$.
Then Claim~\ref{1stclaim} implies that $q=0$, so then $t=k$.
So for all $\chi$, either $(Z_{k-1},Z_0)=(X_{k-1} \cup L_i,\emptyset)$, or $(Z_{k-1},Z_0)=(X_{k-1},L_i)$.
We may assume that we always have the second partition, as the number of perfect colourings is a function of $(Z_1,\ldots,Z_{k-1},Z_0)$ and there are more choices with the second.
Thus we may assume that each perfect colouring $\chi$ gives rise to the same partition $(Z_1,\ldots,Z_{k-1},Z_0)=(X_1,\ldots,X_{k-1},L_1 \cup \ldots \cup L_q)$, where $q \leq k-1$.

Let $P_q(k-1)$ be the set of ordered partitions of $[k-1]$ into $q$ parts (some of which may be empty) with the following property.
Suppose that exactly $a$ different values appear in $\bm{\ell}$, so $\ell_1 = \ldots = \ell_{i_1} > \ell_{i_1+1} = \ldots = \ell_{i_2} > \ldots > \ell_{i_{a-1}+1} = \ldots = \ell_q$ for some $a \in [q]$.
For each unordered partition $\{I_1,\ldots,I_q\}$, labelled so that $\min\{I_1\}<\ldots < \min\{I_q\}$, and any permutation $\sigma : [q]\to[q]$ that has $\sigma(i_{s-1}+1)< \ldots < \sigma(i_{s})$ for all $1 \leq s \leq a+1$ (where $i_0:=0$ and $i_{a+1}:=q$), we put $(I_{\sigma(1)},\ldots,I_{\sigma(q)})$ into $P_q(k-1)$. 
So if $\ell_1>\ldots>\ell_q$, then $P_q(k-1)$ consists of all ordered partitions of $[k-1]$ into $q$ parts, while if $\ell_1=\ldots=\ell_q$, then for each unordered partition it contains exactly one ordering.
%For example, if $\ell=2$ and $\ell_1>\ell_2$, then $P_2(3)$ consists of $(1,23),(2,13),(3,12),(23,1),(13,2),(12,3)$, where~e.g.~$12$ is shorthand for $\{1,2\}$;
%but if $\ell_1=\ell_2$, then $P_2(3)$ consists of $(1,23),(2,13),(3,12)$.

Every colouring in $\mathcal{P}(G)$ is obtained as follows: First, colour edges between pairs in $Z_1,\ldots,Z_{k-1}$ arbitrarily and colour edges between the $L_i$ with colour $2$.
By Claim~\ref{1stclaim}, $i(x) \neq i(x')$ for $x,x'$ in different small parts, so choose $(I_1,\ldots,I_q) \in P_q(k-1)$ so that $i(x) \in I_j$ for all $x$ in $L_j$. Now for each $j\in[q]$ and $x$ in $L_j$, colour all edges between $x$ and $Z_{i(x)}$ 
with colour $1$ and all edges between $x$ and the remaining $Z_i$ arbitrarily.
The number of choices in this colouring procedure is at most
\begin{align*}
f(\bm{m};\bm{\ell}) := &\prod_{ij\in\binom{[k-1]}{2}}2^{m_im_j}\cdot 2^{m\ell}\cdot g(\bm{m};\bm{\ell}),
\quad\text{where}\quad g(\bm{m};\bm{\ell}) := \sum_{\substack{(I_1,\ldots,I_q)\\\in P_q(k-1)}}\prod_{j \in [q]}\left(\sum_{i' \in I_j}2^{-m_{i'}}\right)^{\ell_j}.\\
\end{align*}
We have $|\mathcal{P}(G)| \leq f(\bm{m};\bm{\ell})$ and would like an almost matching lower bound.
Every colouring from the procedure above is valid %.
%Indeed, suppose there is a $K_{k+1}$ in colour $1$.
%It contains at least two vertices in small parts since there are $k-1$ large parts. But every edge between small parts is coloured $2$, a contradiction.
%Suppose there is a $K_k$ in colour $2$, which contains $a \geq 1$ vertices $x_1,\ldots,x_a$ in small parts.
%Since these vertices are in different parts, $i(x_1),\ldots,i(x_a)$ are distinct.
%Then none of other vertices in the $K_k$ can lie in any of $Z_{i(x_1)},\ldots,Z_{i(x_a)}$, so there are only $k-1-a$ allowed parts, a contradiction.
%So all these colourings are valid and 
so we need to check which colourings have been counted more than once.
Suppose $\chi$ is a perfect colouring arising from two different choices.
Let $(I_1,\ldots,I_q)$ and $(I_1',\ldots,I_q')$ be the respective choices of partitions, which are necessarily different.
So there are $i \in [q]$ and $s \in [k-1]$ such that $s \in I_i'\setminus I_i$.
But then in $\chi$ every edge between $L_i$ and $\bigcup_{i'\in I_i \cup \{s\}}Z_{i'}$ is coloured with colour $1$.
In other words, the number of choices for colours of edges between $L_i$ and $Z_1\cup\ldots \cup Z_{k-1}$ has decreased by a multiplicative factor of $2^{m_s\ell_i}$.
The number of such colourings is therefore at most 
$$
%\prod_{ij}2^{m_im_j}\cdot 2^{dc} \cdot \sum_{(I_1,\ldots,I_\ell)}\sum_{i \in [\ell]}\sum_{s \in [k-1]\setminus I_\ell}\prod_{j \in [\ell]\setminus\{i\}}\left(\sum_{i'\in I_j}2^{-m_{i'}}\right)^{\ell_j} \cdot \left(\sum_{i'' \in I_i}2^{-m_{i''}-m_s}\right)^{\ell_i}
f(\bm{m};\bm{\ell}) \cdot q \cdot (k-1) \cdot 2^{-m_{k-1}\ell_q} \leq 2^{-n/k}\cdot f(\bm{m};\bm{\ell}).
$$ 
Thus
\begin{equation}\label{approxPeq}
(1-2^{-n/k})f(\bm{m};\bm{\ell}) \leq |\mathcal{P}(G)| \leq f(\bm{m};\bm{\ell}). 
\end{equation}
Now we can state that we are done in the case $t=k-1$ and the case $t=k$ and every part has size at least $\delta n$.
Indeed, the previous inequality, the fact that
$
f(\bm{m};(\ell)) = h(\bm{m};\ell)
$
and Proposition~\ref{numerical} imply that
\begin{equation}\label{2moreeq}
F(n;k+1,k) \geq (2+o(1))\cdot 2^{t_{k-1}(n)}.
\end{equation}
Moreover, for $3 \leq k \leq 10$, among graphs in $\mathcal{J}_{2(k-1)}(n)$, the unique graph with the largest number of valid colourings is $J^{\ell(k,j)}(n)$ where $\ell(k,j)$ is in the table in Proposition~\ref{numerical}\ref{2colii}.

We return to the only remaining case when $t \geq k$ and there is at least one part of size less than $\delta n$.
Write $f(\bm{m};\ell)$ and $g(\bm{m};\ell)$ for $f(\bm{m};(\ell))$, $g(\bm{m};(\ell))$ respectively.
In the next series of claims we will compare $|\mathcal{P}(G)| = (1+o(1)) f(\bm{m};\bm{\ell})$ with $|\mathcal{P}(G')| = (1+o(1)) f(\bm{m}';\bm{\ell}')$ where $G'$ is another complete partite graph with slightly different part sizes $\bm{m}',\bm{\ell}'$, in order to gradually pin down what $\bm{m}$ and $\bm{\ell}$ must be.
The first of these claims states that merging small parts does not decrease the number of perfect colourings.
(This is not yet sufficient for us to conclude that extremal graphs always have a single small part, as for this the number of perfect colourings would need to increase by a factor $1+\Omega(1)$.)

\begin{claim}\label{better}
If $q \geq 2$, then
$g(\bm{m};\ell) \geq g(\bm{m};\ell_1,\ldots,\ell_q) + 2^{-m_{k-1}\ell}$, so $f(\bm{m};\ell)\geq f(\bm{m};\ell_1,\ldots,\ell_q)$.
\end{claim}

\bcpf
If we obtain $G'$ from $G$ by merging $L_1,\ldots,L_q$ to obtain a single part $L$, then the number of perfect colourings increases.
Indeed, let $x \in L_1$ and $y \in L_2$, say. In $G$, the number of valid choices of $(i(x),i(y))$ is strictly less than the number in $G'$, since $i(x) \neq i(y)$ in $G$, while there is no such restriction in $G'$.
More precisely, the $(I_1,\ldots,I_q)$ term in $g$ corresponds to all the colourings that come from choosing $i(x) \in I_i$ for all $i \in [q]$ and $x \in L_i$.
In $G'$ we can choose $i(x)=k-1$ for all $x \in L$, giving the new term $2^{-m_{k-1}\ell}$.
\ecpf

In the next three claims, we assume that $G$ has one small part of size $\ell$.
Write $m_1\geq \ldots \geq m_{k-1}$ for its other parts and let $a_i := m_i-m_{k-1}$ for all $i \in [k-1]$.
By~(\ref{approxPeq}), as in the proof of Proposition~\ref{numerical}, to compare the number of perfect colourings of $G$ with another complete partite graph $G'$ with a single small part of size $\ell'$ and large parts of size $m_1',\ldots,m_{k-1}'$, it suffices to compare their associated functions $\tilde{h}(\bm{a};\ell),\tilde{h}(\bm{a}',\ell')$ defined in~(\ref{tildeh}), where $\bm{a}=(a_1,\ldots,a_{k-1})$, and $\bm{a}'=(a_1',\ldots,a_{k-1}')$ is chosen so that $\sum_i a_i +\ell = \sum_i a_i' + \ell'$.
(Note that in Proposition~\ref{numerical}, we defined $b_i := m_i-\lfloor\frac{n}{k-1}\rfloor$ whereas here it is convenient to define $a_i = m_i-m_{k-1}$ so that $a_1 \geq \ldots \geq a_{k-1}=0$.)
The next claim concludes the proof in the case when $\ell$ is large and not every $a_1,\ldots,a_{k-2}$ is comparable with $\ell$.

\begin{claim}\label{smallclaim}
Suppose $q=1$.
If $a_{k-2} \leq (\ell-1)/2$ and $\ell \geq \min\{13,2k-2\}$ then there exists an $n$-vertex graph $G'$ with $|\mathcal{P}(G')| > 1.01|\mathcal{P}(G)|$.
\end{claim}

\bcpf
Subtract $1$ from $\ell$ and add $1$ to $a_{k-2}$, letting
$\bm{a}^{(1)} := (a_1,\ldots,a_{k-2}+1,a_{k-1})$.
We let $G'$ be the complete partite graph with large parts of size $(m_{k-1},\ldots,m_{k-1})+\bm{a}^{(1)}$ and one small part of size $\ell-1$. Then
\begin{align}\label{al-1}
(1+o(1))\frac{|\mathcal{P}(G')|}{|\mathcal{P}(G)|} = \frac{\tilde{h}(\bm{a}^{(1)};\ell-1)}{\tilde{h}(\bm{a};\ell)} = 
\frac{2^{\ell-a_{k-2}-1}(2^{-a_1}+\ldots + \frac{1}{2}\cdot2^{-a_{k-2}} + 1)^{\ell-1}}{(2^{-a_1}+\ldots + 2^{-a_{k-2}} + 1)^{\ell}}.
\end{align}
Now, for $y>0$ and $0 < x \leq 1$, let
$$
f(x,y) := \frac{(y+\frac{x}{2}+1)^{\ell-1}}{(y+x+1)^\ell},\quad\text{then}\quad
\frac{\partial f}{\partial x} = -\frac{2(y+\frac{x}{2}+1)^{\ell}((\ell+1)(y+1)+x)}{(y+x+1)^{\ell+1}(2(y+1)+x)^2},
$$
so $f$ is decreasing in $x$ and hence $f(x,y) \geq f(1,y) = (y+\frac{3}{2})^{\ell-1}/(y+2)^\ell$.
Also
$$
\frac{\partial f(1,y)}{\partial y}= -\frac{(y-\frac{\ell}{2}+2)(y+\frac{3}{2})^{\ell-2}}{(y+2)^{\ell+1}}
%\quad\text{and}\quad
%\frac{\partial^2 f(1,\frac{\ell}{2}-2)}{\partial y^2} = \frac{\ell}{4}(1-\ell)
$$
so $f(1,y)$ is increasing from $y=0$ to $y=\frac{\ell}{2}-2 \geq k-3$.
Thus for $0 \leq y \leq k-3$ we have $f(1,y) \geq f(1,0) = (\frac{3}{2})^{\ell-1}/2^\ell = \frac{1}{2}\left(\frac{3}{4}\right)^{\ell-1}$.
Since $2^{-a_1}+\ldots+2^{-a_{k-3}} \leq k-3$, and $\ell \geq 13$, we see that~(\ref{al-1}) is at least
$$
2^{\ell-a_{k-2}-1} \cdot \frac{1}{2}\left(\frac{3}{4}\right)^{\ell-1} \geq \frac{1}{2}\left(\frac{\sqrt{2}\cdot 3}{4}\right)^{\ell-1} > 1.013,
$$
completing the proof of the claim.
\ecpf

\begin{claim}\label{bigclaim}
Suppose $q=1$. If $a_{k-2} \geq \ell/2$ and $\ell \geq \max\{14,2k-2\}$ then there exists an $n$-vertex graph $G'$ with $|\mathcal{P}(G')| > 1.05|\mathcal{P}(G)|$.
\end{claim}

\bcpf
Let $A := \{a_1,\ldots,a_{k-2},\ell\}$ and $b := \sum_{u \in A}u$. Then, recalling $a_{k-1}=0$,
%Let $b_i := a_i$ for $i \in [k-2]$ and $b_{k-1} := \ell$, and $b := \sum_{i \in [k-1]}b_i$. Then
\begin{align*}
%& 2^{\ell\sum_{i}a_i}\prod_{ij}2^{a_ia_j}(2^{-a_1}+\ldots+2^{-a_{k-1}})^\ell
\tilde{h}(\bm{a};\ell) \stackrel{(\ref{tildeh})}{=} \prod_{uu' \in \binom{A}{2}}2^{uu'}(2^{-a_1}+\ldots + 2^{-a_{k-2}} + 1)^{\ell} \leq 2^{t_{k-1}(b)}(2^{-a_1}+\ldots + 2^{-a_{k-2}} + 1)^{\ell}.
\end{align*}
Note that
$$
t_{k-1}(n) = t_{k-1}((k-1)m_{k-1}+b)=\binom{k-1}{2}m_{k-1}^2+(k-2)m_{k-1}b+t_{k-1}(b).
$$
Thus $(1+o(1))|\mathcal{P}(G)| = h(\bm{m};\ell) \leq 2^{t_{k-1}(n)}\cdot (2^{-a_1}+\ldots + 2^{-a_{k-2}} + 1)^{\ell}$.
But
$$
(2^{-a_1}+\ldots + 2^{-a_{k-2}} + 1)^{\ell} \leq ((k-2)2^{-a_{k-2}}+1)^{\ell} \leq \left(\left(\frac{\ell}{2}-1\right)\cdot2^{-\ell/2}+1\right)^\ell.
$$
This is less than $1.9$ for $\ell \geq 14$.
Equation~(\ref{2moreeq}) implies the existence of $G'$ with $|\mathcal{P}(G')| = (2+o(1))\cdot 2^{t_{k-1}(n)}$, and further
$2/1.9 > 1.052$.
\ecpf

The next claim implies that $a_1$ and hence all of $a_1,\ldots,a_{k-1}$ cannot be much bigger than $\ell$.

\begin{claim}\label{dparts}
Suppose $q=1$. If $a_1 \geq \ell+2$, then there is an $n$-vertex graph $G'$ with $|\mathcal{P}(G')| > 1.9|\mathcal{P}(G)|$.
\end{claim}

\bcpf
Suppose that $a_1 \geq \ell+2$.
Let $\bm{a}^{(2)} := (a_1-1,a_2,\ldots,a_{k-2},a_{k-1}+1)$.
Then, recalling $a_{k-1}=0$, 
$$
\frac{\tilde{h}(\bm{a}^{(2)};\ell)}{\tilde{h}(\bm{a};\ell)}
= 2^{a_1-1}\left(\frac{2\cdot 2^{-a_1}+2^{-a_2}+\ldots +2^{-a_{k-2}}+\frac{1}{2}\cdot 2^{-a_{k-1}}}{2^{-a_1}+\ldots + 2^{-a_{k-1}}}\right)^\ell \geq 2^{a_1-1-\ell} \geq 2,
$$
since we can lower bound the numerator by $\frac{1}{2}\sum_{i \in [k-1]}2^{-a_i}$.
\ecpf

\begin{claim}\label{onepart}
There is a positive function $p$ of $k$ such that the following holds.
If $a_1,\ldots,a_{k-1},\ell = O_k(1)$ and $q \geq 2$, then the graph $G'$ obtained by merging $L_1,\ldots,L_q$ has $|\mathcal{P}(G')|>(1+p(k))|\mathcal{P}(G)|$.
\end{claim}

\bcpf
It suffices to show that there is $a=a(k)$ such that $f(\bm{m};\ell) \geq (1+a)f(\bm{m};\ell_1,\ldots,\ell_q)$.
We have $g(\bm{m};\bm{\ell}) = 2^{-m_{k-1}\ell}g(\bm{a};\bm{\ell})$, where $\bm{a}=(a_1,\ldots,a_{k-1})$ satisfies $a_i=m_i-m_{k-1}$, so $0 \leq a_i = O_k(1)$ for all $i$.
Since also $\ell = O_k(1)$, we have $g(\bm{a};\bm{\ell})=p_0(k)$ for some positive function $p_0$ of $k$.
We have $g(\bm{a};\ell)\geq g(\bm{a};\bm{\ell}) + 1$ by Claim~\ref{better}, so
$$
\frac{f(\bm{m};\ell)}{f(\bm{m};\ell_1,\ldots,\ell_q)} = \frac{g(\bm{m};\ell)}{g(\bm{m};\ell_1,\ldots,\ell_q)} \geq \frac{p_0(k)+1}{p_0(k)},
$$
as required.
\ecpf

\begin{claim}\label{atleast2}
If $q=1$ and $\ell=1$ and $k \geq 4$, then the graph $G'$ obtained by moving one vertex from $Z_1$ to $L=L_1$ has $|\mathcal{P}(G')|>1.9|\mathcal{P}(G)|$. 
\end{claim}

\bcpf
Let $\bm{a}^{(3)}:= (a_1-1,a_2,\ldots,a_{k-1})$. Then
\begin{align*}
\frac{\tilde{h}(\bm{a}^{(3)};2)}{\tilde{h}(\bm{a};1)} &\geq \frac{2^{a_1-2}\left(2\cdot 2^{-a_1}+2^{-a_2}+\ldots + 2^{-a_{k-1}}\right)^2}{ 2^{-a_1}+2^{-a_2}+\ldots + 2^{-a_{k-1}}} \geq 2^{-1} + 2^{a_1-a_1-2}+\ldots + 2^{a_1-a_{k-1}-2}\\
&\geq \frac{1}{2}+\frac{k-1}{4} \geq \frac{5}{4},
\end{align*}
where the second inequality follows by expanding $(2^{-a_1}+\sum_{i \in [k-1]}2^{-a_i})^2$.
\ecpf

To complete the proof of Lemma~\ref{close}, assume that $G \notin \mathcal{J}_{2(k-1)}(n)$.
We will find an $n$-vertex graph $G'$ with $F(G';k+1,k)>F(G;k+1,k)$.
Let $m_1,\ldots,m_{k-1}$ be the large parts and $\ell_1,\ldots,\ell_q$ be the small parts of $G$.
Let $a_i := m_i-m_{k-1}$ for all $i \in [k-1]$.
Obtain $G_0$ by merging $L_1,\ldots,L_q$ and let $\ell := \ell_1+\ldots +\ell_q$.
%Suppose $m_1-m_{k-1} \geq k+1$.
%Then Claims~\ref{better},~\ref{dparts} and~(\ref{approxPeq}) imply that there is an $n$-vertex graph $G'$ such that $|\mathcal{P}(G)| < (1+o(1)|\mathcal{P}(G_0)| < \frac{1}{1.8}\cdot|\mathcal{P}(G')|$, as required.
 
By Claim~\ref{better}, we have that $|\mathcal{P}(G_0)| \geq (1+o(1))|\mathcal{P}(G)|$.
Suppose $\ell \geq \max\{14,2k-2\}$.
Then Claims~\ref{smallclaim} and~\ref{bigclaim} imply that there is an $n$-vertex graph $G'$ such that $|\mathcal{P}(G)|<(1+o(1))|\mathcal{P}(G_0)|<(1+o(1))\frac{1}{1.01}\cdot|\mathcal{P}(G')|$, as required.
Thus we may suppose instead that $\ell < \max\{14,2k-2\}$.

Suppose further that $a_1 \geq \ell+2$.
Then Claim~\ref{dparts} implies that there is an $n$-vertex graph $G'$ such that $|\mathcal{P}(G)| < (1+o(1))|\mathcal{P}(G_0)| < \frac{1}{1.8}\cdot|\mathcal{P}(G')|$.
So we may suppose that $\max\{14,2k-2\} \geq \ell+1 \geq a_1 \geq \ldots \geq a_{k-2} \geq a_{k-1}=0$. 

Suppose $q \geq 2$. Then Claim~\ref{onepart} implies that $|\mathcal{P}(G_0)|>(1+p(k))|\mathcal{P}(G)|$ for some positive function $p$.
So we may suppose that $q=1$.
Suppose now $\ell=1$.
Claim~\ref{atleast2} furnishes us with the required $G'$ when $k \geq 4$.

If $k \geq 8$, then we have that $14 \leq 2k-2$ and of course $k \geq 4$, so in fact Claim~\ref{atleast2} implies that $q=1$, $\ell \geq 2$ and $2k-2 \geq \ell,a_1,\ldots,a_{k-1}$, which is a contradiction since $G \notin \mathcal{J}_{2(k-1)}(n)$.
Thus we have $k \leq 7$ and $14 \geq \ell,a_1,\ldots,a_{k-1}$.
Here one can use the same program as in the proof of Proposition~\ref{numerical}\ref{2colii} with tweaked parameters to obtain that in every case there is $G'$ with large part differences $a_1',\ldots,a_{k-1}'$ such that $a_1' \leq 1$ and $\ell'$ is as in the table there, such that $\tilde{h}(\bm{a}';\ell') > \tilde{h}(\bm{a};\ell) + \Omega(1)$, since $G \neq G'$.
\epf

We derive Theorem~\ref{kk+1thm} from Lemma~\ref{close} using the same approach as~\cite{PY}.

\bpf[Proof of Theorem~\ref{kk+1thm}]
We want to show that
\begin{equation}\label{Jkn}
|\mathcal{P}(J)| = O_k(1) \cdot 2^{t_{k-1}(n)}\quad\text{for all } J \in \mathcal{J}_{2(k-1)}(n),
\end{equation}
which implies the same equality holds for $F(J;k+1,k)$ by the arguments in the proof of Lemma~\ref{close}.
For this, start with $T_{k-1}(n)$ with vertex partition $Y_1,\ldots,Y_{k-1},Y_k=\emptyset$.
Consider moves where each time we move one vertex between parts.
One can obtain $J$ in $O(k^2)$ moves, between large parts, or from a large part to the small part $Y_k$.
Each move changes $F$ by a multiplicative factor of $O_k(1)$ (we already did the calculations in Claims~\ref{dparts} and~\ref{atleast2}).
This proves~(\ref{Jkn}).

Apply Lemma~\ref{close} to obtain $\eps,n_0$. We may assume that $1/n_0 \ll \eps \ll 1$.
Let $N := n_0^2$.
Let $G$ be a $(k+1,k)$-extremal graph on $n \geq N$ vertices.
Suppose, for a contradiction, that $G \not\in \mathcal{J}_{2(k-1)}(n)$.
Let $G_n := G$.
If $G_n$ is not a complete multipartite graph, apply Lemma~\ref{small} to obtain a graph $H_n$ on $n$ vertices which is $(k+1,k)$-extremal and complete multipartite, with a part of size $1$ and let $G_n := H_n$.
%Otherwise, let $H_n := G_n$.
Observe that in both cases $G_n \not\in \mathcal{J}_{2(k-1)}(n)$.
We iteratively apply the following procedure.
Let $G_m$ be the current graph on $m \geq n_0+2$ vertices with $\log F(G_m;k+1,k) \geq \left( \frac{k-2}{k-1}-\eps \right)\binom{m}{2}$, and apply Lemma~\ref{close}.
If~(\ref{43def}) fails for some $x \in V(G_m)$, we let $G_{m-1} := G_m-x$, decrease $m$ by $1$, and repeat.
Similarly, if~(\ref{43def}) fails for some pair $y,z \in V(G_m)$, we let $G_{m-2} := G_m-y-z$, decrease $m$ by $2$, and repeat.
Note that
$$
\log F(G_{m-1};k+1,k) \geq \log F(G_m;k+1,k) - \left(\frac{k-2}{k-1} - 2\eps \right)m \geq \left(\frac{k-2}{k-1}-2\eps \right)\left(\binom{m}{2}-m\right).
$$
If~(\ref{43def}) holds for all $x \in V(G)$ and pairs $y,z \in V(G)$, but $G_{m} \not\in \mathcal{J}_{2(k-1)}(m)$, then we replace $G_m$ by the graph $G'$ returned by Lemma~\ref{close} and repeat the step (without decreasing $m$).
Recall that $G'$ has strictly more valid colourings than $G_m$.
Note that for every $m$ for which $G_m$ is defined we have
$$
\log F(G_m;k+1,k) \geq \log F(G;k+1,k) - \left(\frac{k-2}{k-1} - 2\eps \right)\left(n+(n-1)+\ldots+(m+1)\right).  
$$
It follows that we never reach $m \leq n_0$ for otherwise, when this happens for the first time, we get that
\begin{align*}
\log F(G_m;k+1,k) &\geq \left(\frac{k-2}{k-1}-\eps\right)\binom{n}{2} - \left(\frac{k-2}{k-1}-2\eps\right)\left(\binom{n+1}{2}-\binom{m+1}{2}\right)\\
&\geq \eps\binom{n}{2} + \left( \frac{k-2}{k-1}-2\eps \right)\binom{m}{2} > \binom{m}{2},
\end{align*}
i.e.~$F(m;k+1,k) > 2^{\binom{m}{2}}$, a contradiction.
Thus we stop for some $m \geq n_0+2$, with $G_m \in \mathcal{J}_{2(k-1)}(m)$.
We cannot have $m=n$ for otherwise any $J_n \in \mathcal{J}_{2(k-1)}(n)$ has more colourings than $G$ by Lemma~\ref{close}.
So $m \leq n-1$.
Almost every colouring of $G_m := J_m \in \mathcal{J}_{2(k-1)}(m)$ is perfect. Thus
\begin{align*}
1+\log |\mathcal{P}(J_m)| &> \log F(J_m;k+1,k)\\
&\geq \log F(G;k+1,k) - \left(\frac{k-2}{k-1} - 2\eps \right)\left(n+(n-1)+\ldots+(m+1)\right).  
\end{align*}
Since $t_{k-1}(\ell)-t_{k-1}(\ell-1)=\lfloor \frac{k-2}{k-1}\ell\rfloor$,~(\ref{Jkn}) implies that there is some constant $C_k>0$ such that $|\mathcal{P}(J_\ell)| \geq C_k \cdot 2^{\frac{k-2}{k-1}\ell-1}|\mathcal{P}(J_{\ell-1})|$ for all $\ell \geq n_0$ and all $J_\ell \in \mathcal{J}_{2(k-1)}(\ell)$, $J_{\ell-1}\in\mathcal{J}_{2(k-1)}(\ell-1)$. Thus
$$
\log F(G;k+1,k) \geq \log|\mathcal{P}(J_n)| \geq \frac{k-2}{k-1}\left(n+\ldots + (m+1)\right) + (C_k-1)(n-m) + \log|\mathcal{P}(J_m)|.
$$
Together with the previous displayed equation, this is a contradiction to $m \leq n-1$.
Thus $G \in \mathcal{J}_{2(k-1)}(n)$.
\epf

\section{Concluding remarks}\label{concludesec}

\subsection{The uniform problem}

Consider $\bm{k}=(k;s)$. In this paper we obtained an exact result for $(3;7)$ and a new proof for $(3;6)$.
For general $(k;s)$, it seems reasonable to expect every $(r,\phi,\ba)$ to be such that $(k-1)|r$, $\ba$ is uniform and every $\phi^{-1}(c) \cong T_{k-1}(r)$. If this were so, then for distinct $c,c'\in[s]$, $\phi^{-1}(c) \cup \phi^{-1}(c')$ would have at most $(k-1)^2$ parts. This would give rise to a constraint in Problem $L$ which is the analogue of~(\ref{constraint4}), namely
$$
(s-1+s-2)d_2 + (s-1+s-2+s-3)d_3 + \ldots + (s-1+s-2+\ldots + s-s)d_s \leq \binom{s}{2}\frac{(k-1)^2-1}{(k-1)^2}.
$$
Perhaps this constraint is always valid.
As in the case $(3;s)$, it alone does not seem enough to give a realisable solution for $s \geq 8$, but for $(5;4)$ and $(5;5)$ it does give rise to a realisable solution, with $\bm{d} = (0,\frac{3}{4},\frac{3}{16})$ and $(0,0,\frac{15}{16},0)$ respectively. The corresponding solution to Problem $Q^*$ has $r=16$, and colourings come from $\mathbb{F}^2_4$.
Thus we make the following conjecture.

\begin{conjecture}\label{5conjec2}
There exists $n_0>0$ such that for all $n \geq n_0$, $T_{16}(n)$ is the unique $(5;4)$-extremal graph, and the unique $(5;5)$-extremal graph.
\end{conjecture}

Its proof would probably follow from the following analogue of Lemma~\ref{RB}, for $b$ not too small.

\begin{problem}\label{5conjec}
Find $b>0$ as large as possible such that the following holds for all sufficiently large $n$.
Let $R,B$ be two $K_5$-free graphs on the same vertex set of size $n$ with $|R|+|B| \geq (\frac{3}{2}-b) n^2/2$.
Then $|R \cup B| \leq \frac{15}{16} n^2/2 + o(n^2)$.
\end{problem}

Suppose $R,B$ are $K_k$-free graphs on the same vertex set of size $n$, with $|R|,|B| \geq (\frac{k-2}{k-1}-d)n^2/2$ and $|R \cup B| = (\frac{(k-1)^2-1}{(k-1)^2}+d')n^2/2$.
The strong stability theorem of F\"uredi~\cite{furedi} implies that at most $dn^2/2$ edges need to be removed from each of $R,B$ to make them $(k-1)$-partite. Thus removing at most $dn^2$ edges from $R \cup B$ yields a $(k-1)^2$-partite graph, so $d' \leq 2d$.
We need $d=\Omega(1)$ and $d'=o(1)$, for example we proved $(d,d')=(\frac{1}{40},o(1))$ in Lemma~\ref{RB}.

To prove Lemma~\ref{RB}, of which
Problem~\ref{5conjec} is an analogue, we used the fact that $R(3,3)-1 -  (3-1)^2=1$, whereas $R(4,4)-1-(4-1)^2 = 8$ and $26 \leq R(5,5)-1-(5-1)^2 \leq 31$. Thus there could be many clique sizes larger than $(k-1)^2$ present in $R \cup B$, and moreover they could have many different red/blue colourings.
Therefore the proof method of Lemma~\ref{RB} seems unlikely to generalise.

\subsection{The two colour problem}
Solving the two colour problem completely seems difficult. 
Recall that $\bm{k}=(k,\ell)$ has the weak extension property but not the strong extension property, so our stability theorem Theorem~\ref{stabilitysimp} applies -- which shows there is a large family of almost extremal graphs -- but our exact theorem Theorem~\ref{exact1} does not.
The method we used to prove Theorem~\ref{2colapproxthm} about $\bm{k}=(k+1,k)$ could be used to prove further two colour results, but the number of possible coloured graph structures one needs to consider in an analogue of Lemma~\ref{close} increases with $k-\ell$.

It seems plausible that the $(5,3)$-extremal graph may be $T_4(n)$, with $F(n;(5,3)) = (3+o(1))\cdot 2^{t_2(n)}$, where almost all valid colourings come from taking one of the $\binom{4}{2}/2$ part-respecting bipartitions and allowing every cross-edge to be either colour, while edges within one part are given colour $1$. Perhaps for $\bm{k}=(\ell(k-1)+1,k)$, the $\bm{k}$-extremal graph is $T_{\ell(k-1)}(n)$, with $F(n;\bm{k})=(\frac{1}{(k-1)!}\binom{\ell(k-1)}{\ell,\ldots,\ell}+o(1)) \cdot 2^{t_{\ell(k-1)}(n)}$, and for $\bm{k}=(\ell(k-1)+j,k)$ for $j=2,\ldots,k-1$, extremal graphs contain at least one small part.

\typeout{}
\bibliography{ER}

\medskip

{\footnotesize \obeylines \parindent=0pt

\begin{tabular}{lll}
Oleg Pikhurko           &\ & Katherine Staden\\
Mathematics Institute	and DIMAP	        &\ & School of Mathematics and Statistics \\
University of Warwick 	  			&\ & The Open University   	 \\
Coventry                             			&\ & Milton Keynes  \\
CV4 2AL				  		&\ & MK7 6AA	\\
UK					&\ & UK	\\
\end{tabular}
}

\begin{flushleft}
{\it{E-mail addresses}:
\tt{o.pikhurko@warwick.ac.uk, katherine.staden@open.ac.uk}.}
\end{flushleft}

\end{document}